\documentclass[12pt]{memo-l}
\usepackage{amssymb,amscd}
\usepackage[mathscr]{eucal}
\newtheorem{thm}{Theorem}[subsection]%[section]
\newtheorem{lem}[thm]{Lemma}

\newtheorem{prop}[thm]{Proposition}

\newtheorem{cor}[thm]{Corollary}

\newtheorem{defin}[thm]{Definition}
\numberwithin{equation}{subsection}

\DeclareRobustCommand{\SkipTocEntry}[4]{}
\makeatletter
\def\l@section{\@tocline{1}{8pt plus 1pt}{0pt}{}{}}
\makeatother

\renewcommand{\o}{\operatorname}
\newcommand{\C}{{\mathbb C}}
\renewcommand{\P}{{\mathbb P}}
\newcommand{\A}{{\mathbb A}}
\newcommand{\bH}{{\mathbb H}}
\newcommand{\Jac}{{\mathsf{Jac}}}
\newcommand{\End}{{\mathsf{End}}}
\renewcommand{\H}{\o{H}}
\newcommand{\R}{{\mathbb R}}
\newcommand{\Z}{{\mathbb Z}}
\newcommand{\J}{{\mathbb J}}
\newcommand{\Id}{{\mathbb I}}
\newcommand{\Os}{{\mathcal O}}
\newcommand{\Aa}{{\mathcal A}}
\newcommand{\Ee}{{\mathcal E}}
\newcommand{\zz}{\mathcal{Z}}
\newcommand{\bb}{\mathcal{B}}
\newcommand{\M}{{\mathcal M}}
\newcommand{\la}{\langle}
\newcommand{\ra}{\rangle}

\newcommand{\hH}{{\mathscr H}}

\newcommand{\bO}{\boldsymbol{\Omega}}
\newcommand{\balpha}{\boldsymbol{\alpha}}
\newcommand{\tbalpha}{\widetilde\balpha}
\renewcommand{\Im}{\o{Im}}
\renewcommand{\Re}{\o{Re}}
\newcommand{\Hol}{\mathsf{Hol}(E)}
\newcommand{\Holu}{\mathsf{Hol}_u(E)}
\newcommand{\hol}{\mathsf{hol}}
\newcommand{\her}{\mathsf{Her}}
\newcommand{\Her}{\mathsf{Her}(E)}
\newcommand{\ac}{{\mathfrak A}(E)} % all linear connections
\newcommand{\fcu}{{\mathcal F}_u(E)} % flat unitary connections
\newcommand{\fcuu}{{\mathcal F}_u(E)} % flat unitary connections - no h
\newcommand{\fcl}{{\mathcal F}_l(E)} % flat linear connections
\newcommand{\g}{{\mathfrak g}}

\newcommand{\ff}{{\mathfrak F}} % flat structure
\newcommand{\gau}{{\mathcal G}_u(E)} % unitary gauge group
\newcommand{\gal}{{\mathcal G}_l(E)} % linear gauge group
\newcommand{\Hig}{\mathsf{Higgs}(E)}
\newcommand{\Hom}{\mathsf{Hom}}
\newcommand{\Mor}{\mathsf{Mor}}
\newcommand{\hpg}{\Hom(\pi,G)/G}
\newcommand{\hmg}{\Hom(\pi,G)}

\newcommand{\Uo}{\o{U}(1)}
\newcommand{\GLr}{\o{GL}(r,\C)}

\newcommand{\GL}{\o{GL}}
\newcommand{\U}{\o{U}}

\newcommand{\PGLtR}{\o{PGL}(2,\R)}

\newcommand{\Sp}{\o{Sp}}

\newcommand{\SOth}{\o{SO}(3)}

\newcommand{\Map}{\mathsf{Map}}
\newcommand{\db}{\overline{\partial}}
\renewcommand{\d}{\partial}
\newcommand\ka{K\"ahler}
\newcommand\kan{K\"ahlerian}

\newcommand\tZ{\tilde{Z}}
\newcommand\tS{\tilde{\Sigma}}
\newcommand\tX{\tilde{X}}

\newcommand\ca{{\mathcal A}}
\newcommand\cb{{\mathcal B}}
\newcommand\cc{{\mathcal C}}
\newcommand\Ss{{\mathcal S}}  
\newcommand\Iso{{\mathsf{Iso}}}

\newcommand\Obj{{\mathsf{Obj}}}
\newcommand\hk{{hyperk\"ah\-ler~}}
\newcommand\Hk{{Hyperk\"ahler~}}

\newcommand\bV{\bar{V}}
\newcommand\pp{\mathsf{p}}
\newcommand\qq{\mathsf{q}}
\newcommand\vv{\mathsf{v}}
\newcommand\xx{\mathsf{x}}
\newcommand\yy{\mathsf{y}}

\newcommand{\xmapsto}[1]{|\hspace{-5pt}\xrightarrow{~{#1}~}}
\newcommand\bA{\mathbf{A}}
\newcommand\bB{\mathbf{B}}
\newcommand\ba{\mathsf{a}}
\newcommand\bbb{\mathsf{b}}

\begin{document}
\frontmatter
\title[Rank One Higgs Bundles]
{Rank One Higgs Bundles and Representations of 
Fundamental Groups of Riemann Surfaces}
%\author[Goldman and Xia]{William M. Goldman and Eugene Z. Xia}
%\address{
%Department of Mathematics, University of Maryland,
%College Park, MD 20742 USA ({\it Goldman}),
%Department of Mathematics,
%National Cheng Kung University,
%Tainan 701, Taiwan
%({\it Xia}) \bigskip}
%\email{
%wmg@math.umd.edu {\it (Goldman)}, 
%ezxia@ncku.edu.tw  {\it (Xia)}}
\author{William M. Goldman}
\address{Department of Mathematics, University of Maryland,
College Park, MD 20742 USA}
\email{wmg@math.umd.edu}
\author{Eugene Z. Xia}
\address{Department of Mathematics,
National Cheng Kung University,
Tainan 701, Taiwan}
\email{ezxia@ncku.edu.tw}
\date{\today}

\begin{abstract}
This expository article details the theory of rank one Higgs bundles
over a closed Riemann surface $X$ and their relation to
representations of the fundamental group of $X$.  We construct an
equivalence between the deformation theories of flat connections and
Higgs pairs. This provides an identification of moduli spaces arising
in different contexts. The moduli spaces are real Lie groups.  From
each context arises a complex structure, and the different complex
structures define a \hk~structure. The twistor space, real forms, and
various group actions are computed explicitly in terms of the Jacobian
of $X$. We describe the moduli spaces and their geometry in terms of
the Riemann period matrix of $X$.

This is the simplest case of the theory developed by Hitchin, Simpson
and others. We emphasize its formal aspects that generalize
to higher rank Higgs bundles over higher dimensional \ka~ manifolds.

\end{abstract}

\thanks{Goldman gratefully acknowledges partial support by NSF grants
DMS-9504764, DMS-9803518, DMS-0103889 and a Semester Research Award from
the General Research Board of the University of Maryland in Fall 1998.
Xia gratefully acknowledges partial support by National Science
Council Taiwan grant NSC 91-2115-M-006-022, 93-2115-M-006-002.}
\maketitle

\begin{center}
{\it To the memory of Hsieh Po-Hsun}
\end{center}

\newpage
\setcounter{page}{5}
\tableofcontents 

\mainmatter

\thispagestyle{plain}
\addtocontents{toc}{\SkipTocEntry}
\chapter*{Rank One Higgs Bundles}
        
\section*{Introduction}
The set of equivalence classes of representations of
the fundamental group $\pi$ of a closed Riemann surface $X$ into a Lie
group $G$ is a basic object naturally associated to $\pi$ and
$G$. Powerful analytic techniques have been employed by
Hitchin, Simpson, Corlette and Donaldson 
et al~\cite{Hitchin1, Simpson1,Simpson2,Simpson3}
to understand the global topology and geometry of this object. 
Rank one Higgs bundles provide a toy model with
a more explicit form than general Higgs bundles.
This paper expounds this case, emphasizing its formal aspects to isolate
and clarify the main ideas and motivate its generalization to higher rank.

We consider three moduli spaces of seemingly different objects: 
complex characters of the fundamental group $\pi$ of a closed and
orientable surface $\Sigma$, 
flat connections on a trivial line bundle over $\Sigma$,
and pairs consisting of a holomorphic line bundle over a Riemann surface
and a holomorphic $1$-form on $X$.
These three moduli problems or {\em deformation theories\/} arise
in different contexts, but are nonetheless equivalent. 
Simpson~\cite{Simpson3} names these deformation theories 
Betti, de Rham, and Dolbeault respectively.
In rank one, the Betti moduli space is the 
collection of ordered $2k$-tuples
of nonzero complex numbers ($k$ is the genus of $\Sigma$).
The de Rham moduli space is the quotient 
\begin{equation*}
\H^1(\Sigma,\C)/\H^1(\Sigma,\Z) \cong \H^1(\Sigma,\C/\Z).
\end{equation*}
The Dolbeault moduli space is the cotangent bundle $T^*\Jac(X)$ of
the Jacobian of the Riemann surface $X$.

This paper assumes a basic knowledge of topology (for example, 
as in Fulton~\cite{Fulton}) and differential, and complex geometry.  
Our perspective is differential geometric and
gauge-theoretic. We assume basic facts about symplectic geometry and
Hamiltonian flows (as in Weinstein~\cite{Weinstein}).  
We hope this exposition of the simplest part of deep
ideas in algebraic geometry will be useful to topologists and
differential geometers interested in moduli spaces of representations
of surface groups.

The paper is organized as follows.  \S 1 is a brief introduction of
groupoids and their equivalences.  \S 2 constructs the Betti, de Rham,
groupoids and moduli spaces and the various structures on these moduli
spaces. The main player in the Betti groupoid is a surface group
$\pi$, whereas the main player in the de Rham groupoid is a closed
smooth surface $\Sigma$ with fundamental group $\pi$.  \S 3 develops
the Dolbeault group\-oid; here $\Sigma$ is given a conformal
structure, making it a Riemann surface $X$ over which we consider
various holomorphic objects.  The inherent complex structure $J$ on the Betti
and de Rham moduli spaces are isomorphic, but $J$ differs from the
inherent complex structure $I$ on the Dolbeault moduli space.  In \S 5,
these two different structures give rise to a  \hk structure 
on the underlying moduli space.
\S 6 explicitly constructs  the {\em twistor space,\/}  
a holomorphic object containing both moduli spaces.
All of these constructions apply to the cotangent bundle of an
arbitrary principally polarized abelian variety (not necessarily the
Jacobian of a curve). In particular they can all be expressed in terms
of the {\em Riemann period matrix\/} $\Pi$ of $X$.

Our story ends (\S 7) by explicitly describing the moduli spaces in
terms of the period matrix. We revisit the Betti moduli
space which appears as a product of the real torus 
\begin{equation*}
\Hom(\pi,\Uo)\cong T^{2k} 
\end{equation*}
with the real symplectic vector space 
\begin{equation*}
\Hom(\pi,\R^+) \cong  \H^1(\Sigma,\R) \cong \R^{2k}. 
\end{equation*}
The new structure
as the (Dolbeault) moduli space identifies the torus
$\Hom(\pi,\Uo)$ as the quotient of $\C^k$ with the lattice $\Z^k + \Pi\Z^k$
corresponding to the columns of $\Pi$. The conformal structure of the
Riemann surface $X$ determines a complex structure $\J_\Pi$ on $\R^{2k}$
which is easily expressed \eqref{eq:JPi} in terms of $\Pi$. The resulting
$\C^*$-action on the moduli space, which determines the full structure
of the Dolbeault moduli space and is equivalent to the {\em Hodge structure\/}
on the moduli space. 

Although many features of the rank one case generalize to higher rank,
several important technical issues are altogether absent in rank one.
We warn the reader of these simplifications.  First, since
the structure group is abelian, the moduli spaces are naturally groups
(or torsors over groups). This strong structure is missing for higher
rank Higgs bundles.  Endomorphism bundles of line bundles are trivial,
but endomorphism bundles are generally nontrivial. In rank one, moduli
spaces exist for all objects. One need not remove from consideration
``unstable'' objects in the sense of Geometric Invariant Theory.
Furthermore the main differential equation ---- Hitchin's self-duality 
equation --- decouples, leaving the Higgs field and the holomorphic 
structure completely independent of each other.
This arises from the direct product decomposition 
$\C^* \cong \Uo \times \R^+$, which does not generalize to higher rank.

We do not discuss the \hk moment map and quotient construction in
this paper, although all the moduli spaces
described here can be constructed as complex-symplectic and \hk quotients.
The constructions in rank one are particularly simple and familiar.
As we do not need this machinery, we do not discuss them here.

For other points of view and related topics, see
\cite{Arapura,ABCKT,Corlette2,Donaldson3,
Simpson1,Simpson2,Simpson3,Simpson4,Simpson5,Toen}
and references cited therein.

We thank Nigel Hitchin, Taejung Kim, Weiping Li, Michael Thaddeus and
Mike Wolf for their critical reading of this manuscript and numerous helpful
suggestions. We are also grateful to Steve Bradlow, 
Robert Bryant, Kevin Corlette,
Simon Donaldson, Lawrence Ein, Elisha Falbel, Charlie Frohman, 
Oscar Garcia-Prada, Steve Kudla,
Fran\c cois Labourie, John Loftin, John Millson, Niranjan Ramachandran,
Jonathan Rosenberg, 
Carlos Simpson, Domin\-go Toledo, Richard Wentworth, and Scott Wolpert 
for helpful conversations.
\bigskip

\noindent{\em Notation and terminology.}
To emphasize the different contexts, we reserve 
$\Sigma$ for the smooth surface, and $X$ for the Riemann surface
``diffeomorphic to  $\Sigma$.'' That is, $X$ is $\Sigma$ with
a conformal structure (which for us is the Hodge $\star$-operator
on $1$-forms). For constructions involving the differential structure,
whether we use $\Sigma$ or $X$ is a decision on the context.
For example, both $\Aa^*(\Sigma)$ and $\Aa^*(X)$ are correct notations
for the de Rham algebra of smooth differential forms on the surface.

Differential forms and cohomology classes are complex-valued,
unless otherwise stated.  Tensor products of modules are over $\Z$ unless
otherwise stated.  If $V$ is a complex vector space, then we denote the 
real vector space underlying $V$  by $V_\R$.
The complex vector space $V$ then consists of the pair $(V_\R,I)$ where
$I:V_\R\longrightarrow V_\R$ is the {\em complex structure,\/} the
$\R$-linear automorphism $I$ satisfying $I^2 = -1$ corresponding to 
scalar multiplication by $\sqrt{-1}$.
The complex vector space $\bV$ {\em complex-conjugate\/} to $V$ is 
$(V_\R,-I)$ where $-I$ is the {\em opposite\/} complex structure.
Thus an {\em anti-linear\/}  map $V\longrightarrow W$ 
into a complex vector space $W$ is a linear map $\bV \longrightarrow W$.

The multiplicative groups of nonzero complex and real numbers are denoted
by $\C^*$ and $\R^*$ respectively. The multiplicative group of positive
real numbers is denoted $\R^+$. The multiplicative group of unit complex
numbers is denoted $\Uo$. We denote the $k\times k$ identity matrix by
$\Id_k$.

By a {\em vector bundle\/} we mean a {\em smooth complex vector bundle.\/}
If $E\longrightarrow M$ is a vector bundle and $N$ is a submanifold
with inclusion map $f:N\hookrightarrow M$, we denote by
$E|_N$ the {\em restriction\/} of $E$ to $N$, that is, the vector bundle
over $N$ defined as pullback $f^*E$ of $E$ by $f$.

\section{Equivalences of deformation theories}
A moduli problem seeks to classify a class of objects up to an equivalence
relation, often defined by a group of transformations of the set of
objects. 

A {\em deformation theory\/} (or {\em transformation
groupoid\/}) $(S, G)$ consists of a category $\cc$ defined by a group
action as follows.
Let  $\alpha:G\times S\longrightarrow S$ be
a left action of a group $G$ on a set $S$.
The {\em deformation theory\/} $(S, G)$ 
consists of the category $\cc$ whose objects form a set
$\Obj(\cc) = S$ with morphisms 
\begin{equation*}
x \stackrel{g}\longrightarrow y.
\end{equation*}
corresponding to triples $(g,x,y)\in G\times S\times S$
such that $\alpha(g,x) = y$. 

The identity element $e\in G$ determines, for each object
$x\in S$ the identity morphism 
\begin{equation*}
x \stackrel{e}\longrightarrow x.
\end{equation*}
The inverse of the morphism
\begin{equation*}
x \stackrel{g}\longrightarrow y
\end{equation*}
is
\begin{equation*}
y \stackrel{g^{-1}}\longrightarrow x
\end{equation*}
and the composition of morphisms
\begin{equation*}
x \stackrel{g}\longrightarrow y\stackrel{h}\longrightarrow z
\end{equation*}
equals
\begin{equation*}
x \stackrel{hg}\longrightarrow z.
\end{equation*}
In particular every morphism is an isomorphism.

The {\em moduli set\/} corresponding to such a deformation theory is
the set $\Iso(\cc)$ of isomorphism classes of objects. The {\em
isotropy group\/} of an object $x\in\Obj(\cc)$ is the set
$\Mor(x,x)$ consisting of morphisms $x\longrightarrow x$, which has the
structure of a group.
An {\em equivalence of categories \/} is a functor
$F:\ca\longrightarrow\cb$ such that there exists a functor
$H:\cb\longrightarrow\ca$ and natural transformations from the
compositions $F\circ H$ and $H\circ F$ to the identity functors of
$\cb$ and $\ca$ respectively. (See Jacobson~\cite{Jacobson} or
Gelfand-Manin~\cite{GelM}, p.28 for discussion of this notion and
Goldman-Millson~\cite{GoldmanMillson} for an application closely
related to this one.) An equivalence of categories induces a bijection
$\Iso(\ca)\longrightarrow \Iso(\cb)$, although in general $\Obj(\ca)$
and $\Obj(\cb)$ may be enormously different. For example, each
groupoid arising from a group $G$ operating on itself by left-multiplication
is equivalent to the groupoid with one object and one morphism.

Equivalent deformation theories yield equivalent moduli sets.
However the finer notion of equivalence has
further implications--- for example isotropy groups of
corresponding points in the moduli spaces are isomorphic.

Often the sets $\Obj(\ca)$ admit additional algebraic or geometric
structures, which induce additional structures on $\Iso(\ca)$. 
For the examples discussed here, these moduli sets are Lie
groups, and the equivalences of deformation theories induces
isomorphisms of (real) Lie groups.

Equivalent deformation theories may have
different structures. An equivalence of a deformation theory $\ca$
with another deformation theory may provide additional structures to
$\Iso(\ca)$.  For example, $\Hom(\pi,\Uo)$ inherits the structure of a
complex abelian variety from every Riemann surface with fundamental
group $\pi$.

The following criterion is a useful tool for proving that a functor 
is an equivalence of categories.
A functor $F:\ca\longrightarrow\cb$ is an {\em equivalence\/} if and only if:
\begin{description}
\item[$\bullet$\quad Surjective on isomorphism classes]
The induced map
\begin{equation*}
F_*:\Iso(\ca)\longrightarrow\Iso(\cb)
\end{equation*}
is surjective;
\item[$\bullet$\quad Full]
For $x,y\in \Obj(\ca)$, the map
\begin{equation*}
F(x,y):\Mor(x,y)\longrightarrow \Mor(F(x),F(y))
\end{equation*}
is surjective;
\item[$\bullet$\quad Faithful]
For $x,y\in \Obj(\ca)$, the map
\begin{equation*}
F(x,y):\Mor(x,y)\longrightarrow \Mor(F(x),F(y))
\end{equation*}
is injective.
\end{description}

\section[Betti and de Rham moduli spaces]{The Betti and de Rham deformation theories and their moduli spaces}

This section describes the Betti and de Rham deformation theories.  
Fix a structure group $G$.  Although much of what is here generalizes to the
case that $G$ is a linear algebraic group, we restrict in this paper
to the cases of $G$ being the group $\C^*$ of nonzero complex numbers,
the group $\Uo$ of unit complex numbers, or the group $\R^*$ of
nonzero real numbers.  In what follows, $\Sigma$ is a compact smooth
oriented surface with fundamental group $\pi$.

\begin{itemize}
\item {\em The Betti groupoid\/} whose objects are
representations $\pi\longrightarrow G$, 
with morphisms $G$;
\item {\em The de Rham groupoid\/} whose objects are
flat connections on a trivial complex line bundle over $\Sigma$, 
with morphisms gauge transformations;
\end{itemize}

\subsection{The Betti groupoid}

Denote by $\Hom(\pi, G)$ the set 
of representations from $\pi$ to $G$.  The group $G$
acts on representations by conjugation.
The {\em Betti groupoid\/} is the category
having $\Hom(\pi,G)$ as the set of objects and morphisms
\begin{equation*}
g:\rho_1 \longrightarrow \rho_2
\end{equation*}
where $g\in G$, $\rho_1,\rho_2\in\hmg$ and
\begin{equation*}
\rho_2 = \iota_g\circ \rho_1 
\end{equation*}
where $\iota_g:G\longrightarrow G$ is the inner automorphism defined
by conjugation by $g$.  The Betti groupoid is $(\hmg, G)$.

The surface group $\pi$ admits a presentation
\begin{equation}\label{eq:presentation}
\la A_1,B_1,\dots,A_k,B_k \mid [A_1,B_1] \dots [A_k,B_k] = 1\ra 
\end{equation}
where $[A,B]$ denotes $ABA^{-1}B^{-1}$. The map
\begin{align*}
\hmg & \longrightarrow G^{2k} \\
\rho & \longmapsto (\rho(A_1),\rho(B_1), \dots, \rho(A_k),\rho(B_k))
\end{align*}
embeds $\hmg$ as the Zariski-closed subset of $G^{2k}$
defined by
\begin{equation}\label{eq:surfacerelation}
[\alpha_1,\beta_1] \dots [\alpha_k,\beta_k] = 1.
\end{equation}

Since $G$ is abelian, it 
acts trivially on $\hmg$. Furthermore \eqref{eq:surfacerelation}
is trivially satisfied and 
\begin{equation}\label{eq:abelian}
\hpg \cong \hmg\cong G^{2k}. 
\end{equation}
Hence the moduli space $\Hom(\pi, \C^*)$  identifies with
$(\C^*)^{2k}$ with a natural complex structure $J$.

Each $\rho \in \Hom(\pi, G)$ provides $\g$ with a $\pi$-module
structure via the adjoint action of $G$ on $\g$.  Since $G$ is
abelian, this $\pi$-module structure is trivial.  
The tangent space of $\Hom(\pi, G)$ at $\rho$ identifies with
the group cohomology $\H^1(\pi, \g)$.
Cup product 
\begin{equation*}
\H^1(\pi, \g) \times \H^1(\pi, \g) \longrightarrow \C,
\end{equation*}
defined by
\begin{equation}\label{eq:BettiOmegaJ}
\Omega(u, v) := \sum_{i = 1}^g u(A_i) v(B_i) - u(B_i) v(A_i),
\end{equation}
provides a closed nondegenerate exterior $2$-form on $\Hom(\pi,G)$. 
When $G = \C^*$, this closed exterior $2$-form $\Omega$ is holomorphic
with respect to $J$. Thus  $(J,\Omega)$ 
defines a {\em complex-symplectic structure\/} 
on  $\Hom(\pi,\C^*)$.

Each element $C \in \pi$ defines a function 
\begin{align}\label{eq:function}
f_C :\Hom(\pi,G) & \longrightarrow \C \notag\\
\rho &\longmapsto \rho(C). 
\end{align}
which determines a Hamiltonian flow on $\Hom(\pi, G)$ 
(Goldman~\cite{invariantfunctions}) as follows. 
Since $G$ is abelian, the function $f_C$ depends only
on the homology class of $C$. By applying an automorphism of $\pi$,
we may assume that $C = (A_1)^n$ for some $n\ge 0$. 
The corresponding flow 
\begin{equation*}
\Phi_t:\Hom(\pi,G)  \longrightarrow \Hom(\pi,G)  
\end{equation*}
for $t\in\R$ is defined on the generators as follows:
\begin{equation}\label{eq:flow}
  \Phi_t(\rho):\begin{cases} 
A_1 &\longmapsto  \rho(A_1) \\
B_1 &\longmapsto  \rho(B_1)e^{i n t} \\
A_j &\longmapsto  \rho(A_j) \text{~\quad if $j>1$}\\
B_j &\longmapsto  \rho(B_j) \text{~\quad if $j>1$}
\end{cases}
\end{equation}

%%%

$\Hom(\pi,\C^*)$ is a complex Lie group under pointwise
multiplication of homomorphisms. 
Namely, if $\rho_1,\rho_2$ are homomorphisms, then
\begin{align}
\pi & \longrightarrow \C^* \label{eq:pointwise}\\
\gamma &\longmapsto \rho_1(\gamma) \rho_2(\gamma) \notag 
\end{align}
is a homomorphism, defining a group structure on $\Hom(\pi,\C^*)$
isomorphic to $(\C^*)^{2k}$. 

Two real forms of $\C^*$ are $\Uo$ and $\R^*$. 
These subgroups are the fixed point sets of anti-involutions of
$\C^*$ defined by
\begin{equation*}
\iota_U: z \longmapsto (\bar{z})^{-1}, 
\qquad \iota_\R: z \longmapsto \bar{z} 
\end{equation*}
respectively. These induce real structures on $\Hom(\pi,\C^*)$ whose sets
of fixed points are $\Hom(\pi,\Uo)$ and $\Hom(\pi,\R^*)$
respectively. Note that the composition of the anti-involutions is the
involution
\begin{equation*}
\iota_U\circ \iota_\R : z \longmapsto z^{-1}.
\end{equation*}

\subsection{The de Rham groupoid}\label{sec:dRgpd}
The second deformation theory concerns flat connections on a smooth
complex line bundle $E$ over a smooth surface $\Sigma$.
The morphisms are gauge transformations, that is,
smooth maps $\Sigma\longrightarrow G$.
A flat connection defines a notion of a smooth section being locally 
constant, or {\em parallel.\/}  The connection is described as a differential
operator which vanishes precisely on locally constant sections.

We restrict to the case when $E$ is the trivial line bundle. 
As the Betti moduli space, 
the de Rham moduli space is a complex Lie group of dimension $2k$,
where $k$ is the genus of $\Sigma$. It furthermore enjoys an invariant
complex-symplectic structure as well as two real structures 
(anti-automorphisms) corresponding to the real forms of $\C^*$.

We first review exterior differential forms,
connections, and gauge transformations.
For background on fiber bundles and their differential geometry, we refer
the reader to Steenrod~\cite{Steenrod} 
and Kobayashi-Nomizu~\cite{KobayashiNomizu}.

\subsubsection*{Exterior differential calculus}
Let $\Aa^*(\Sigma)$ denote the {\em de Rham algebra of $\Sigma$,\/} that is,
the differential graded algebra of $\C$-valued smooth exterior 
differential forms on $\Sigma$. The operations are {\em wedge product\/}
\begin{equation*}
\Aa^k(\Sigma)\times \Aa^l(\Sigma) \longrightarrow \Aa^{k+l}(\Sigma)
\end{equation*}
and {\em exterior derivative\/}
\begin{equation*}
d:\Aa^k(\Sigma) \longrightarrow \Aa^{k+1}(\Sigma)
\end{equation*}
satisfying:
\begin{itemize}
\item
$d(\xi\wedge\eta) = d\xi\wedge\eta + (-1)^k \xi\wedge d\eta$
if $\xi\in\Aa^k(\Sigma)$; 
\item $d\circ d = 0$.
\end{itemize}
Complex-conjugation $\eta\longmapsto\bar{\eta}$ defines an
anti-involution of the complex vector space $\Aa^*(\Sigma)$. Its fixed-point
set is the subalgebra $\Aa^*(\Sigma;\R)$ of real differential forms.
Its $-1$-eigenspace consists of the subspace
$\Aa^*(\Sigma;i\R)$ of purely imaginary differential forms.

Let $E$ be a (smooth complex) vector bundle over $\Sigma$.
Let $\Aa^k(\Sigma;E)$ denote the collection of $E$-valued exterior
$k$-forms over $\Sigma$; then the graded vector space
\begin{equation*}
\Aa^*(\Sigma;E) = \bigoplus_{k\ge 0} \Aa^k(\Sigma;E)
\end{equation*}
is a graded module over $\Aa^*(\Sigma)$.

Let $E$ be a trivial line bundle.
A {\em trivialization\/} $\tau$ is a nonvanishing section of $E$.
Then $\Aa^k(\Sigma;E)$ identifies with $\Aa^k(\Sigma)$ via
\begin{align}\label{eq:section}
\Aa^k(\Sigma) &\longrightarrow \Aa^k(\Sigma;E) \notag\\
\eta & \longmapsto \eta \tau.
\end{align}

\subsubsection*{Gauge transformations}

A {\em linear gauge transformation\/} of $E$ is a smooth bundle automorphism
$\xi:E\longrightarrow E$ covering the identity map of $\Sigma$.  
Thus for each $x\in S$, the gauge transformation $\xi$ acts on 
the fiber $E_x$ by scalar multiplication of some scalar $g(x)\in\C^*$. 
In terms of $\tau$,
\begin{equation*}
\xi(\tau) = g \cdot \tau 
\end{equation*}
but the smooth map $g:\Sigma\longrightarrow\C^*$ is independent of
$\tau$ since $\C^*$ is abelian.  The group $\gal$ of linear gauge
transformations of $E$ identifies with the space $\Map (\Sigma,\C^*)$
of smooth maps $\Sigma\longrightarrow\C^*$.  Similarly, the
$\Uo$-gauge group is the subgroup $\Map (\Sigma,\Uo)$ of smooth
$\Uo$-valued maps.

The fundamental groups of $G=\Uo$ and $\C^*$ are infinitely cyclic, 
so a smooth map $g\in\Map(\Sigma,G)$ induces a homomorphism
\begin{equation*}
\pi_1(g) : \pi \longrightarrow \Z,
\end{equation*}
determining an element of 
\begin{equation*}
\Hom(\pi,\Z) \cong \H^1(\Sigma,\Z).  
\end{equation*}
The resulting group homomorphism
$\Map(\Sigma,G) \longrightarrow \H^1(\Sigma,\Z)$
induces an isomorphism of the group of connected components 
\begin{equation*}
\pi_0\big(\Map(\Sigma,G)\big) \cong \H^1(\Sigma,\Z).
\end{equation*}
Its kernel is the identity component $\Map(\Sigma,G)^0$ consisting
of null-homotopic maps $\Sigma\longrightarrow G$. Since
\begin{align*}
\Z \hookrightarrow \C & \stackrel{\Ee}\longrightarrow \C^* \\
z & \longmapsto \exp(2\pi i z),
\end{align*}
is a universal covering space, $\Map(\Sigma,\C^*)^0$ identifies with
\begin{equation*}
\Map(\Sigma,\C) = \Aa^0(\Sigma)  
\end{equation*}
via the exact sequence
\begin{equation*}
\Map(\Sigma,\C) \stackrel{\exp_*}\longrightarrow 
\Map(\Sigma,\C^*) \stackrel{S}\longrightarrow \H^1(\Sigma,\Z).
\end{equation*}
Similarly the identity component of $\Map(\Sigma,\Uo)$ identifies
with 
\begin{equation*}
\Map(\Sigma,i\R) = \Aa^0(\Sigma,i\R). 
\end{equation*}

\subsubsection*{Connections on vector bundles}
A {\em connection\/} on $E$ is an operator 
\begin{equation*}
D:\Aa^0(\Sigma;E)\longrightarrow \Aa^1(\Sigma;E), 
\end{equation*}
such that
\begin{equation}
D(f s) = f D(s) + df \wedge D(s).
\end{equation}
($D(s)$ is called the {\em covariant differential\/} of $s$ with respect
to $D$.)
Such a map extends to  
\begin{equation*}
D:\Aa^p(\Sigma;E)\longrightarrow \Aa^{p+1}(\Sigma;E) 
\end{equation*}
by enforcing the identity
\begin{equation*}
D(\alpha\wedge\beta) = D(\alpha) \wedge\beta + (-1)^k \alpha\wedge D(\beta)
\end{equation*}
for any $\alpha\in\Aa^k(\Sigma)$, $\beta\in\Aa^l(\Sigma;E)$.
We denote the space of all connections on $E$ by $\ac.$

A trivialization $\tau$ determines a connection
\begin{equation*}
D_0 : \Aa^0(\Sigma;E)\longrightarrow \Aa^1(\Sigma;E):
\end{equation*}
by the rule
\begin{equation}
D_0(f\tau) := df \wedge \tau.
\end{equation}
(where $f\in\Aa^0(\Sigma)$ is a smooth function determining a section 
$f\tau$ by \eqref{eq:section}).
This is the unique connection for which $\tau$ is {\em parallel\/}
(has covariant differential zero).

With respect to this connection, an arbitrary connection $D$ has the form
\begin{equation*}
D = D_0 + \eta
\end{equation*}
where $\eta\in\Aa^1(\Sigma)$
and $\eta$ acts by exterior
multiplication. 
That is, the covariant differential of a section
$s = f \tau \in\Aa^0(\Sigma;E)$  with respect to $D$ is
\begin{equation*}
D(s) =  df \wedge \tau  + \eta\wedge f \tau \in \Aa^1(\Sigma;E).
\end{equation*}
In particular the 1-form $\eta$ is given by
\begin{equation*}
D(\tau) = \eta \wedge \tau.
\end{equation*}

The {\em curvature\/} of a connection $D$ is the $\End(E)$-valued 
exterior 2-form 
$F(D) \in \Aa^2(\Sigma,\End(E))$
such that for any section $s$ of $E$,
\begin{equation*}
D(D(s)) = F(D) \wedge s.
\end{equation*}
If $E$ is a line bundle, then $\End(E)$ is canonically trivial, so
$F(D)$ identifies with an ordinary exterior 2-form.

With respect to the trivialization $\tau$, the curvature is:
\begin{equation*}
F(D) = d\eta\in\Aa^2(\Sigma).
\end{equation*}
Hence a connection $D$ is {\em flat,} that is, its curvature vanishes, 
if and only if $d\eta = 0$.
Thus the space $\fcl$ of flat connections identifies with
the subspace $Z^1(\Sigma)\subset \Aa^1(\Sigma)$ of closed 1-forms.
A line bundle $(E,D)$, where $D\in\fcl$ is a {\em flat line bundle.\/}

\subsubsection*{Gauge transformations and the de Rham groupoid} 
Let $\tau$ be a trivialization and
let $D_0$ be the corresponding connection.  If $\xi\in\gal$
corresponds to a map $g\in\Map(S,\C^*)$, then the action of
$\xi$ on a connection $D_0 + \eta$ is given by:
\begin{equation} \label{eq:gauge action}
\xi\cdot (D_0 + \eta)  := D_0 + \eta + g^{-1} dg 
\end{equation}
where $\eta\in\Aa^1(\Sigma)$.

To see this, let $\xi\cdot\tau$ be the trivialization obtained by transforming
$\tau$ by $\xi$. That is, $\xi\cdot\tau = g\tau$ (scalar multiplication).
An arbitrary section of $E$ is given by a scalar multiple
$s = f \tau$. Let $D = D_0 + \eta$ be an arbitrary connection. The
effect of the transformed connection $\xi\cdot D$ on $s$ is given by:
\begin{align*}
(\xi\cdot D)(s) & = D (f g\tau) = df \wedge g \tau + f dg\wedge \tau + f
g D(\tau) \\ 
& = df \wedge (\xi\cdot\tau) + f g^{-1} dg \wedge (g\tau) + 
f g \eta \wedge \tau \\ 
& = \big(df + f (g^{-1} dg + \eta)\big) \wedge \big(\xi\cdot\tau\big) 
\end{align*}
(since $\C^*$ is abelian).
In particular this action is independent of $\tau$.

Naturality of curvature under gauge transformations
\begin{equation*}
F(\xi\cdot D) = \xi^*(F(D))
\end{equation*}
follows, in our context, from the closedness of exact forms and the
triviality of the action of $\gal$ on $\Aa^2(\Sigma)$:
\begin{equation*}
F(\xi\cdot D) = d(\eta + g^{-1}dg ) = d\eta = F(D) = \xi^*F(D).
\end{equation*}
Hence the $\gal$-action preserves curvature.
The de Rham groupoid is $(\fcl, \gal)$.

\subsubsection*{The de Rham moduli space}
The space $\ac$ of all connections is an affine space modeled on the
space $\Aa^1(\Sigma)$ comprising smooth $\C$-valued 1-forms on $\Sigma$. 
Let $D_0$ be the connection corresponding to the trivialization $\tau$.
Every connection on $E$ is of the form $D_0 + \eta$, where 
$\eta\in\Aa^1(\Sigma)$. The objects of the de Rham groupoid --- 
flat connections --- comprise an affine subspace of $\ac$.

With $\tau$, the space $\fcl$ of flat connections on $E$ 
identifies with the vector space $\zz^1(\Sigma)$ of closed 1-forms.
We shall always assume a fixed trivialization $\tau_0$ for $E$,
which will provide a basepoint $D_0$ (an {\em origin\/}) for 
the affine space $\ac$ of connections. 

The gauge group $\gal$ identifies with $\Map(\Sigma,\C^*)$ and its   
action on $\ac$ is given by \eqref{eq:gauge action}. 

The $\gal$-action on $\ac$ decomposes into the action of the
identity component
$\gal^0 = \Map(\Sigma,\C^*)^0$ on $\ac$ and the action of $\pi_0(\gal)$
on the quotient. If $g\in\gal^0$ then $g = \exp f$ for some 
$f\in\Aa^0(\Sigma)$, and the action on $\ac$ is given by:
\begin{equation*}
D_0 + \eta \stackrel{g}\longmapsto  D_0 + \eta + df.
\end{equation*}
Thus the quotient $\fcl/(\gal^0)$ is an affine space whose underlying
vector space is the cohomology
\begin{equation*}
\H^1(\Sigma) := \zz^1(\Sigma)/\bb^1(\Sigma).
\end{equation*}

The action of $\pi_0(\gal)\cong \H^1(\Sigma;\Z)$ corresponds to the
action of $\H^1(\Sigma,\Z)$ on $\H^1(\Sigma)$ via the monomorphism
\begin{equation}\label{eqn:deRhamQuotient}
i^* : \H^1(\Sigma,\Z)\longrightarrow \H^1(\Sigma)
\end{equation}
induced by the
coefficient inclusion $\Z \stackrel{i}{\hookrightarrow} \C$.
In summary:

\begin{prop}\label{prop:dRms}
The de Rham moduli space $\fcl/\gal$ identifies with the cokernel of the
map $\H^1(\Sigma,\Z)\longrightarrow \H^1(\Sigma)$ induced by
$\Z\hookrightarrow \C$.
\end{prop}
\noindent
Therefore the moduli space inherits 
a complex structure $J$ arising from the operation
\begin{equation}\label{eq:J}
J(\eta) := i\eta 
\end{equation}
on $1$-forms.
Cup product 
\begin{equation*}
\H^1(\Sigma) \times \H^1(\Sigma) \longrightarrow \C,
\end{equation*}
defined on the level of $1$-forms by
\begin{equation}\label{eq:OmegaJ}
\Omega(\alpha, \beta) := \int_\Sigma \alpha \wedge \beta,
\end{equation}
induces a nondegenerate exterior $2$-form $\Omega$ on  $\fcl/\gal$.
Since $\Omega$ is parallel on this affine space, it is closed. 
Furthermore, this closed exterior $2$-form $\Omega$ is holomorphic
with respect to $J$. Thus$(J,\Omega)$, 
defines a complex-symplectic structure on $\fcl/\gal$.

Every $C\in\pi$ defines a function
\begin{align*}
f_C: \fcl & \longrightarrow \C \\ 
\omega &\longmapsto \int_C \omega
\end{align*}
corresponding to the function \eqref{eq:function} on the Betti moduli
space, and with a Hamiltonian flow corresponding to \eqref{eq:flow}.

The identity
\begin{equation*}
\Omega(J\alpha,\beta) =  \Omega(\alpha,J\beta) =  i\Omega(\alpha,\beta) 
\end{equation*}
means that $\Omega$ is complex-bilinear with respect to $J$ (that is, has
type $(2,0)$).

\subsubsection*{Group structure}
The de Rham moduli space $\H^1(\Sigma)/\H^1(\Sigma,\Z)$ is a Lie group
as a quotient of a vector space by a discrete subgroup.
Addition corresponds to tensor product of flat line bundles as follows.
The trivialization $\tau_0$ of $E$ determines a trivialization 
$\tau_0\otimes\tau_0$ of $E$ which we henceforth identify with $\tau_0$.
Let $(E,D_1)$ and $(E,D_2)$ be two flat line bundles with 
local sections $s_1, s_2$ respectively.
In terms of the trivialization, write $s_i = f_i \tau_0$ 
where $s_i\in\Aa^0(\Sigma)$ are functions.
The {\em tensor product flat bundle\/}
$(E,D_1)\otimes(E,D_2)$ is the unique flat bundle $(E,D_1\otimes D_2)$
for which $s_1\otimes s_2$ is parallel with respect to $D_1\otimes D_2$,
where $s_1$ and $s_2$ are local sections of $E$  parallel with respect
to  $D_1$  and $D_2$ respectively.

Write 
\begin{align*}
D_1 & = D_0 + \eta_1 \\ D_2 & = D_0 + \eta_2 \\ D_1\otimes D_2 & = D_0 + \eta
\end{align*}
where $\eta_1,\eta_2,\eta\in\zz^1(\Sigma)$. 
Then $s_i$ is $D_i$-parallel if and only if
\begin{equation*}
0 = D_i s_i =  d f_i + \eta_i f_i
\end{equation*}
for $i=1,2$ whence 
\begin{align*}
(D_1\otimes D_2) (s_1 \otimes s_2) & =  
D_1(s_1) \otimes s_2 + s_1\otimes  D_2(s_2) \\ & =
\big((d f_1 + \eta_1 f_1) f_2 +
f_1 (d f_2 + \eta_2 f_2)\big) \tau_0 \\ & =
D_0 (s_1 \otimes s_2) + (\eta_1 + \eta_2) s_1\otimes s_2.
\end{align*}
Thus the $1$-form $\eta$ corresponding to the tensor product
$D_1\otimes D_2$ equals the sum $\eta_1 + \eta_2$.

\subsubsection*{Real structure}

Corresponding to the real structures $\iota_\R$ and $\iota_U$ on $\C^*$
are real structures on  the de Rham moduli space. 
Their fixed point sets are {\em real forms\/} of the de Rham moduli
space, corresponding to sub-deformation theories for the real forms
$\Uo$ and $\R^*$ of $\C^*$.
On the Lie algebra $\C$ of $\C^*$ the anti-involutions $\iota_U,\iota_\R$ 
respectively induce maps
\begin{align}
\eta &\overset{\iota_U}\longmapsto -\bar{\eta} \notag\\
\eta & \overset{\iota_\R}\longmapsto \bar{\eta}\label{eq:iotacc}
\end{align}
which are related by the composition
\begin{equation*}
\iota_U\circ\iota_\R:\eta \longmapsto -\eta.
\end{equation*}
These anti-involutions on the Lie algebras are coefficient maps
for anti-involutions on the moduli spaces, denoted by $(\iota_U)_*$
and $(\iota_\R)_*$ respectively.
The fixed point set of $(\iota_U)_*$ is the 
moduli space $\fcu/\gau$ of flat unitary connections.
The fixed point set of $(\iota_\R)_*$ 
is the moduli space of flat real connections.

Using the real structure, the complex-symplectic form decomposes
into real and imaginary parts
\begin{equation}\label{eq:OmegaRI}
\Omega = \Omega^{\Re} + i \Omega^{\Im}, 
\end{equation}
each of which is a {\em real symplectic structure \/} on the
de Rham moduli space.  They satisfy:
\begin{align*}
\Omega^{\Re}(J\alpha,J\beta) & = -\Omega^{\Re}(\alpha,\beta) \\
\Omega^{\Im}(J\alpha,J\beta) & = -\Omega^{\Im}(\alpha,\beta) \\
\Omega^{\Im}(\alpha,\beta) & = \Omega^{\Re}(J\alpha,\beta).
\end{align*}

\subsection{Equivalence of de Rham and Betti groupoids}\label{sec:holonomy}
Let $G$ be $\C^*$.  The cases of $\R^*, \R^+$ and $\Uo$ are similar.  
Let $D$ be a flat connection on the trivial line bundle $E$
over $\Sigma$. Let $\tau_0$ denote a trivialization, and write
\begin{equation*}
D \tau_0 = \eta \tau_0 
\end{equation*}
for a closed 1-form $\eta$.  Over any smooth path
\begin{equation*}
[a,b]\stackrel{\gamma}\longrightarrow\Sigma 
\end{equation*}
{\em parallel transport\/} defines a linear map between the fibers
\begin{equation*}
\Pi^{-1}(\gamma(a)) \longrightarrow \Pi^{-1}(\gamma(b))
\end{equation*}
as follows.
The induced line bundle $\gamma^*E$ over $[a,b]$ 
has an induced flat connection 
$\gamma^*D$ as well as an induced trivialization $\tau_0\circ\gamma$.
Thus
\begin{equation*}
\gamma^*\eta = g(t) dt 
\end{equation*}
for a unique function $g$ on $[a,b]$. Similarly, every section
$s$ is
\begin{equation*}
s(t) = f(t) (\tau_0\circ\gamma) 
\end{equation*}
for a unique function $f(t)$.
The section $s$ is parallel if and only if:
\begin{equation*}
0 = (\gamma^*D) s = \big(  f'(t) + g(t) f(t) \big) \; dt\, \otimes\, 
(\tau_0\circ\gamma),
\end{equation*}
that is,
\begin{equation*}
0 =  f'(t) + g(t) f(t). 
\end{equation*}
This differential equation has unique solution
\begin{equation*}
f(t) = \exp \bigg( -\int_a^t g(s) ds\bigg)  f(a)  
\end{equation*}
for an arbitrary initial condition $f(a)$.
Starting from 
\begin{equation*}
s(a) = f(a) \tau_0\big(\gamma(a)\big),  
\end{equation*}
the corresponding section $s(t) = f(t) \tau_0(\gamma(t))$ 
is the parallel transport
of $s(a)$ along $\gamma$. Flatness of $D$ implies that the parallel
transport from $\Pi^{-1}\gamma(a)$ to $\Pi^{-1}\gamma(b)$ along $\gamma$
depends only on the homotopy class of $\gamma$ relative to its endpoints.

If $\gamma$ is a {\em based loop,\/} that is, $\gamma(a)=\gamma(b) = x_0$,
then holonomy defines a homomorphism
\begin{equation}\label{eq:holonomy}
\hol_p(D) \; :\;
\pi_1(\Sigma;x_0) 
\longrightarrow G.
\end{equation}
(Compare \cite{KobayashiNomizu}.)

Let $\xi \in \gal$ be a gauge transformation.  Evaluation of $\xi$ at $x_0$
gives a homomorphism of groups 
$\gal \rightarrow G$,
for which
\eqref{eq:holonomy} is equivariant.
Together these maps give the {\em holonomy functor\/} 
between the de Rham and the Betti groupoid 
\begin{thm}\label{thm:one1}
The holonomy functor 
\begin{equation*}
\hol: (\fcl, \gal) \longrightarrow (\Hom(\pi,G),G)
\end{equation*}
is an equivalence of deformation theories.
\end{thm}
\begin{proof}

Let $\Sigma$ be a smooth manifold with basepoint $x_0 \in\Sigma$. 
Let $\pi = \pi_1(\Sigma, x_0)$ be the corresponding fundamental group
and $\tS \stackrel{\Pi}\longrightarrow\Sigma$ the corresponding universal 
covering space.
Corresponding to a representation $\rho\in\Hom(\pi,G)$ 
is a {\em flat line bundle\/} 
$\C_\rho\longrightarrow\Sigma$, defined as follows.
The group $\pi$ acts on the total space $\tS\times \C$ of the trivial 
line bundle over $\tS$ by deck transformations on the first factor and
via $\rho$ on the second factor:
\begin{equation}\label{eq:piaction}
(\tilde s, x) \stackrel{\gamma}\longmapsto 
(\gamma\cdot \tilde s, \rho(\gamma)x) 
\end{equation}
where $x\in G$.                                 
The quotient $(\tS\times \C)/\pi$ is the total space of a smooth
line bundle 
\begin{equation*}
\C_{\rho} \stackrel{\Pi}\longrightarrow \Sigma  
\end{equation*}
which carries a natural {\em flat structure\/}
$\ff_\rho$, a foliation of $\C_\rho$ transverse to the fibration such that
the restriction of $\Pi$ to each leaf of $\ff_\rho$ is a 
covering space of $\Sigma$. The leaves are the projections
of the leaves $\tS \times \{v\}$ of the product foliation of $\tS\times\C$
under the quotient map
\begin{equation*}
\tS \times \C \longrightarrow \C_\rho. 
\end{equation*}
The representation $\rho$ is the {\em holonomy representation\/} of 
$\C_{\rho}$.

The exterior derivative
\begin{equation*}
\tilde{d}:\Aa^k(\tS) \longrightarrow \Aa^{k+1}(\tS)
\end{equation*}
induces the trivial connection $\tilde{D}$ 
on the bundle $(\tS\times \C)$
\begin{equation*}
\tilde{D}:\Aa^k(\tS, \tS \times \C) \longrightarrow \Aa^{k+1}(\tS, \tS \times \C)
\end{equation*}
by
considering sections as $\C$-valued functions and then taking
exterior derivative.  The $\pi$-action via $\rho$
is equivariant with respect to $\tilde{D}$.  Since flatness is a local
condition, $\tilde{D}$ descends to a flat connection $D$ on 
$(\tS\times \C)/\pi$

A local section $s$ is {\em parallel\/} if and only if $D(s) = 0$ 
if and only if its image lies in a leaf
of $\ff_\rho$. Equivalently, a lift of $s$ to $\tS$ has constant projection under
the Cartesian projection 
\begin{equation*}
\tS \times \C \longrightarrow \C
\end{equation*}
defining the foliation.

The {\em covering homotopy theorem\/} for fiber bundles implies that
any continuous path $\rho_t\in\hmg$, where $0\le t\le 1$, 
induces an isomorphism of line bundles
\begin{equation*}
\C_{\rho_0}  \longrightarrow \C_{\rho_1}.
\end{equation*}
Since $\hmg$ is an $\R$-algebraic set, its connected components are
path-connected and the topological type of the bundle $\C_\rho$ 
depends only on the connected component of $\hmg$ containing $\rho$.
If $\hmg$ is connected, then any representation can be connected to
the trivial representation by a continuous path, and therefore
defines the trivial bundle.  The existence of such a $D$ for each
$\rho \in \hmg$ shows that $\hol_p$ is surjective
on isomorphism classes.

Next we show $\hol_p$ is faithful and full.  Recall
that \eqref{eq:holonomy} is equivariant with respect to 
homomorphism
$\gal \rightarrow G$.  Hence $\gal$-equivalent connections give rise
to $G$-equivalent representations. 
Suppose 
\begin{equation*}
D_1, D_2 \in \fcl.
\end{equation*}
There are two cases, depending on whether or not
$D_1, D_2 $ are 
$\gal$-equivalent.
If $D_1, D_2$ are not
$\gal$-equivalent, then 
$\Mor(D_1, D_2)$
is empty.
Otherwise, there exists $\xi \in \gal$ corresponding to a map 
$g\in\Map(X,G)$ such that
\begin{equation*}
D_2 = D_1 +  g^{-1}d g
\end{equation*}
and if $\xi_1$ is another gauge transformation corresponding to a map
$g_1\in\Map(X,G)$ such that
\begin{equation*}
D_2 = D_1 + g_1^{-1} d g_1,
\end{equation*}
then $g_1 = g c$, where $c$ is a constant map. 
We denote the subgroup of $\gal$ corresponding to constant maps
$X\longrightarrow G$ by $G$.
Hence
$\Mor(D_1, D_2)$ corresponds to the 
coset $g\cdot G$.

Now suppose that $D_1,D_2\in\fcl$ and $\hol_p(D_i)=\rho_i$ for 
$i=1,2$. If $D_1$ is not $\gal$-equivalent to $D_2$, then 
then $\Mor(\rho_1,\rho_2)=\emptyset$.
Otherwise, $\rho_1 = g \rho_1 g^{-1} = \rho_2$ where $g\in G$.
$\Mor(\rho_1,\rho_2) = g G g^{-1} = G$.  Then we have an
isomorphism
\begin{equation*}
g^{-1} : \Mor(D_1, D_2) \longrightarrow \Mor(\rho_1,\rho_2).
\end{equation*}

In both cases, $\hol_p$ induces an isomorphism
\begin{equation*}
\Mor(D_1, D_2) \longrightarrow \Mor(\hol_p(D_1), \hol_p(D_2))
\end{equation*}
as desired. 
\end{proof}
For a careful treatment of holonomy for higher ranks and in general, see
Kobayashi-Nomizu~\cite {KobayashiNomizu} and 
Goldman-Millson~\cite{GoldmanMillson}.
Theorem~\ref{thm:one1} is stated and proved in \cite{GoldmanMillson}.

The holonomy representation of a tensor product of flat line bundles
is the product of the holonomy representations, as defined in
\eqref{eq:pointwise}. Therefore the corresponding isomorphism
of moduli spaces is an isomorphism of groups.

If $D = D_0 + \eta$ is a flat connection, then its holonomy is the
homomorphism
\begin{align*}
\pi & \longrightarrow G \\
\gamma &\longmapsto \exp \int_\gamma \eta
\end{align*}
which depends holomorphically on $D$. Thus the isomorphism of moduli spaces
is an isomorphism of complex Lie groups.  
This isomorphism is also
symplectomorphic with regard to the complex-symplectic structures $\Omega$
defined above.
%%%

\section{The Dolbeault groupoid}\label{sec:dolbgpd}
Let $X$ be a Riemann surface diffeomorphic to $\Sigma$. 
With this extra structure 
is associated a third deformation theory, 
the {\em Dolbeault groupoid\/}.
In this section we define this deformation theory and relate its
moduli space to the Jacobian $\Jac(X)$.

\subsection{Holomorphic line bundles}
We begin by reviewing  differential calculus on a Riemann surface.
A holomorphic structure is defined as an extension $D''$ of the $\db$-operator
to sections of a line bundle. For line bundles, a Higgs field is a holomorphic
$1$-form on $X$. A {\em Higgs pair\/} is a holomorphic structure 
together with a Higgs field. 

\subsubsection*{Exterior calculus on Riemann surface}
The Hodge $\star$-operator on $1$-forms 
\begin{equation*}
\star:\Aa^1(X) \longrightarrow \Aa^1(X) 
\end{equation*}
describes the conformal structure on $X$ as follows.
Let $\J:TX\longrightarrow TX$ denote the complex structure on $TX$:
for each $x\in X$, multiplication by $\sqrt{-1}$
is an $\R$-linear automorphism $\J_x$ of the (real) tangent space
$(T_xX)_\R$ with 
\begin{equation*}
\J_x\circ \J_x = -1. 
\end{equation*}
A $1$-form $\eta\in\Aa^1(X)$ defines
an $\R$-linear map $T_xX \longrightarrow \C$. The Hodge $\star$-operator
is the map on $\Aa^1(X)$ induced by $\J$, that is, 
define $\star\eta\in\Aa^1(X)$ as the composition
\begin{equation*}
T_xX \xrightarrow{\,\J^{-1}}
T_xX \stackrel{\eta}\longrightarrow \C.
\end{equation*}
(When $X$ is given a conformal Riemannian metric, this operator agrees
with the Riemannian $\star$-operator. For
$1$-forms on a $2$-manifold, $\star$ is independent of the Riemannian metric.)

Thus $\star$ defines an endomorphism of the complex vector space
$\Aa^1(X)$ with $\star\star  = -1$, so its
eigenvalues are $\pm i$. Furthermore
$\Aa^1(X)$ decomposes as a direct sum of its $i$-eigenspace and
$(-i)$-eigenspace respectively:
\begin{equation*}
\Aa^1(X)     = \Aa^{1,0}(X) \oplus \Aa^{0,1}(X)
\end{equation*}
and the summands are complex-conjugates of each other.

Composing the exterior derivative
$d:\Aa^0(X)\longrightarrow\Aa^1(X)$ with projections
of $\Aa^1(X)$ onto $\Aa^{1,0}(X)$ and 
$\Aa^{0,1}(X)$ defines operators
\begin{align*}
\d:\Aa^0(X)& \longrightarrow\Aa^{1,0}(X) \\
\db:\Aa^0(X)& \longrightarrow\Aa^{0,1}(X)
\end{align*}
with $d = \d + \db$.

A {\em holomorphic $1$-form\/} is a $(1,0)$-form which is closed;
a closed $(0,1)$-form is an {\em anti-holomorphic $1$-form.\/}
A {\em harmonic $1$-form\/} is the sum of a holomorphic
$1$-form and an antiholomorphic $1$-form.
The spaces of holomorphic $1$-forms, antiholomorphic $1$-forms,
and harmonic $1$-forms, denoted by $\hH^{1,0}(X)$, $\hH^{0,1}(X)$, and
$\hH^1(X)$ respectively, satisfy:
\begin{equation*}
\hH^1(X) = \hH^{0,1}(X)\oplus \hH^{1,0}(X).
\end{equation*}
where the two summands are complex-conjugates of each other.

Later we need the following basic facts, which can be found in, 
for example, Farkas-Kra~\cite{FarkasKra}, Forster~\cite{Forster}
(Theorem~19.9, p.156, Theorem~19.12, p.159) Griffiths-Harris~\cite{Gr}, 
Jost~\cite{Jost},\S\S 3.3,5.2, Wells~\cite{Wells}.

\begin{itemize}
\item
 Let $\eta\in\Aa^{0,1}(X)$.
Then there exists a unique antiholomorphic $1$-form $\eta_0\in\hH^{0,1}(X)$ 
and a function $f\in\Aa^0(X)$ such that
\begin{equation*}
\eta = \eta_0 + \db f 
\end{equation*}
\item
Let $\eta\in\zz^1(X)$ be a closed 1-form. Then there exists a unique
harmonic $1$-form $\eta_0$ and a function $f\in\Aa^0(X)$ such that
\begin{equation*}
\eta = \eta_0 + d f. 
\end{equation*}
\end{itemize}

These standard facts form the analytic 
foundation for the theory expounded 
here.

\subsubsection*{Holomorphic structures}  
Just as a connection on a vector bundle provides a differential criterion
characterizing sections to be locally constant, or {\em parallel,\/}
a holomorphic structure on a vector bundle over a complex manifold provides
a notion of a section to be holomorphic. 

From a preferred class of {\em holomorphic local sections,\/} one can find
an atlas of local trivializations of the vector bundle such that the
transition functions on overlapping coordinate patches are holomorphic.
We do not take that approach here, instead referring to Gunning\cite{Gunning1},
Kobayashi~\cite{Kobayashi}.

The trivialization $\tau_0$ determines a
holomorphic structure, such that $s = f \tau_0$ is a 
holomorphic section defined on an open set $U$  if and only if
$f$ is a holomorphic function on $U$.  In particular 
holomorphic local sections solve the {\em Cauchy-Riemann
equation\/}.

Holomorphic structures are conveniently
described by differential operators of degree 1
\begin{equation*}
D'' : \Aa^{p,q}(X;E)\longrightarrow \Aa^{p,q+1}(X;E)
\end{equation*}
which satisfy
\begin{equation*}
D'' (f\cdot s) = \db f\wedge s + f \cdot D''(s)
\end{equation*}
for $f\in\Aa^0(X)$.
Locally every holomorphic structure admits holomorphic sections,
and are thus equivalent to the standard one $D_0''$ (which
arises from $\db$ and the trivialization).
This follows from the solvability of the 
inhomogeneous Cauchy-Riemann equation.
See Atiyah-Bott~\cite{AtiyahBott}, \S 5 (pp.\ 554--55) 
or Kobayashi~\cite{Kobayashi}, Chapter 1, Propositions 3.5--3.6, p.9
for further discussion.

Denote by $D_0'$ and $D_0''$ the (1,0)- and (0,1)- parts of $D_0$
respectively:

\begin{align*}
D_0' & = (D_0)^{1,0}\\
D_0'' & = (D_0)^{0,1}
\end{align*}
so that 
\begin{equation*}
D_0 = D_0' + D_0''.
\end{equation*}
With respect to the trivialization $\tau_0$, 
a {\em holomorphic structure\/} on $E$ is an operator
\begin{equation*}
D'' = D_0'' + \Psi
\end{equation*}
(where $\Psi\in\Aa^{0,1}(X)$), which operates on smooth sections by:
\begin{align*}
\Aa^{0,0}(X,E) & \stackrel{D''}\longrightarrow\Aa^{0,1}(X,E) \\
      f\tau_0               & \longmapsto  (\db f + f\Psi)\tau_0.
\end{align*}
(When $\dim(X)>1$, holomorphic structures must satisfy
an extra integrability condition. 
Compare Kobayashi~\cite{Kobayashi}.) 

Denote the space of all holomorphic structures on $E$ by $\Hol$.
The linear gauge group $\gal$ acts on $\Hol$ by
\begin{equation}\label{eq:gaugehol}
D_0'' + \Psi \stackrel{\xi}{\longmapsto} D_0'' + \Psi + g^{-1}\db g.
\end{equation}
Denote the resulting deformation theory by $(\Hol,\gal)$.

\subsubsection*{Higgs fields}  
If $(E,D'')$ is a holomorphic vector bundle over $X$,
then a {\em Higgs field\/} on $(E,D'')$ is 
a (1,0)-form $\Phi$ on $X$ taking values in the
endomorphism bundle $\End(E)$. The Higgs field $\Phi$
is required to be {\em holomorphic\/} with respect to the
holomorphic structure on $T^*M\otimes\End(E)$ induced by
the complex structure on $X$ and the holomorphic structure
$D''$ on $E$. Thus a Higgs field is a holomorphic bundle
map $TX \longrightarrow \End(E)$. 
(In higher dimensions, the Higgs field is required to satisfy the
integrability condition $\Phi\wedge\Phi=0$, which automatically holds
when $X$ is 1-dimensional.)

A {\em Higgs pair\/} (or {\em Higgs bundle\/}) 
is a pair $((E,D''),\Phi)$ where
$(E,D'')$ is a holomorphic vector bundle and $\Phi$ is
a Higgs field on $(E,D'')$. 

When $E$ is a holomorphic line bundle, $\End(E)$ is trivial and
a Higgs field is a holomorphic $1$-form. Thus a rank one Higgs
pair is just a holomorphic structure $D''$ together
with a holomorphic $1$-form. 
Thus the space of Higgs pairs $(E,D'',\Phi)$ is
\begin{equation}\label{eq:hig}
\Hig = \Hol \times \hH^{1,0}(X).
\end{equation}
The linear gauge group $\gal$ acts on $\Hol$ by 
\eqref{eq:gaugehol} and on Higgs fields by conjugation.
Since multiplication in $\End(E)$ is commutative, 
$\gal$ acts trivially on $\hH^{1,0}(X)$.
Denote the resulting groupoid by $(\Hig, \gal)$.

The groupoid $(\Hig, \gal)$ contains the subgroupoid of holomorphic
line bundles $(\Hol, \gal)$ as Higgs pairs whose Higgs fields are
identically zero.

\subsection{The moduli spaces}
Now we describe the Dolbeault moduli spaces of rank one holomorphic
bundles and Higgs bundles over a compact Riemann surface $X$. 
These moduli spaces will identify with spaces
of unitary and linear representations of $\pi_1(X)$ respectively.

\subsubsection*{Moduli of holomorphic structures}
The space $\Hol$ of all holomorphic structures is an affine space
modeled on $\Aa^{0,1}(X)$.  If $D_0''$ is the holomorphic structure
corresponding to a trivialization, then every holomorphic structure is
of the form $D_0'' + \Psi$, where $\Psi\in\Aa^{0,1}(X)$.  As for
connections on a trivial line bundle $E$, we fix a trivialization and
its corresponding holomorphic structure.  The affine spaces 
have natural vector space structures, in which vector addition
corresponds to tensor product of holomorphic line bundles.

The $\gal$-action on $\Hol$ is given by
\begin{equation}\label{eq:galactionhol}
\Psi \mapsto \Psi + g^{-1} \db g, 
\end{equation}
where $\Psi\in\Aa^{0,1}(X)$ and $g\in\Map(X,\C^*)$. 

The $\gal$-action on $\Hol$ decomposes into the action of the
identity component
\begin{equation*}
\gal^0 = \Map(X,\C^*)^0  
\end{equation*}
on $\Hol$ and the action of $\pi_0(\gal)$
on $\Hol/(\gal^0)$.
If $g\in\gal^0$ then $g = \exp f$ for some 
$f\in\Aa^0(X)$, and the action on $\Hol$ is given by:
\begin{equation*}
D_0'' + \Psi \stackrel{g}\longmapsto  D_0'' + \Psi + \db f.
\end{equation*}
Thus the quotient $\Hol/(\gal^0)$ is an affine space whose underlying
vector space is the quotient $\Aa^{0,1}(X)/\db \Aa^0(X)$, 
the {\em Dolbeault cohomology group.\/} 
The Hodge decomposition into antiholomorphic forms and $\db$-exact forms
(\S\ref{sec:dolbgpd})
\begin{equation*}
\Aa^{0,1}(X)  = \hH^{0,1}(X)  \oplus \db\Aa^0(X) 
\end{equation*}
gives an isomorphism of
$\Aa^{0,1}(X)/\db\Aa^0(X)$
with the space $\hH^{0,1}(X)$ of
antiholomorphic $1$-forms on $X$, which we henceforth denote as $V$.
Thus $\Hol/(\gal^0)$ identifies with the complex vector space 
V = $\hH^{0,1}(X)$. 
The dimension of $V$ equals the genus $k$ of $X$.

The action of $\pi_0(\gal)\cong \H^1(X;\Z)$ corresponds to the
action of $\H^1(X,\Z)$ on $\hH^{0,1}(X)$ by 
translation via the composition map
\begin{equation}\label{eqn:JacQuotient}
\H^1(X,\Z) \stackrel{i^*}{\longrightarrow} 
\H^1(X) \stackrel{p_1}{\longrightarrow} \hH^{0,1}(X)  = V
\end{equation}
where $i^*$ is as in (\ref{eqn:deRhamQuotient}) and
$p_1$ is the projection to antiholomorphic $1$-forms.
The image of $\H^1(X,\Z)$ is a lattice $L\subset V$ of
rank 
\begin{equation*}
2 k =\dim (V_\R).  
\end{equation*}
The quotient 
$\Jac(X) := V/L$ is a complex torus, the {\em Jacobian\/} of $X$. 

The {\em Dolbeault isomorphism\/} identifies $\H^1(X,\Os)$ with
$\hH^{0,1}(X)$, where $\Os$ is the sheaf of germs of holomorphic
functions on $X$. Under the Dolbeault isomorphism, the above homomorphism
\begin{equation*}
\H^1(X,\Z) \longrightarrow \hH^{0,1}(X) 
\end{equation*}
identifies with the map induced by the inclusion $\Z \longrightarrow \Os$
of the sheaf of germs of locally constant $\Z$-valued functions.
Its cokernel then identifies with the kernel
of the homomorphism $\H^1(X,\Os^*)\longrightarrow \H^2(X,\Z)$ induced by the
sheaf homomorphism
\begin{align*}
\Os &\stackrel{\Ee}\longrightarrow \Os^* \\
f &\longmapsto \exp(2\pi i f),
\end{align*}
by the exactness of 
\begin{equation*}
0 \longrightarrow \Z \longrightarrow \Os  \stackrel{\Ee}\longrightarrow
\Os^* 
\end{equation*}
where $\Os^*$ is the multiplicative sheaf of germs of nonvanishing
holomorphic functions.

Just as the tensor product of flat line bundles corresponds to addition
of closed $1$-forms, the tensor product of holomorphic line bundles
corresponds to addition of antiholomorphic $1$-forms, that is, to
addition in $V$. Namely if $D_i'' = D_0'' + \Psi_i$, then the tensor product
holomorphic structure is
\begin{equation*}
D_1''\otimes D_2'' = D_0'' + \Psi_1 + \Psi_2.
\end{equation*}
This multiplicative structure agrees with the multiplicative structure
induced by the multiplicative structure on $\Os^*$.

\subsubsection*{Moduli of Higgs pairs}
Let $(D'',\Phi)\in\Hig$ be a Higgs pair.
Its equivalence class is recorded by the equivalence class of
the holomorphic bundle $(E,D'')$ in $\Jac(X)$ and the Higgs field
$\Phi$ in $\hH^{1,0}(X)$. 

\begin{prop}\label{prop:dolbms}
The moduli space $\Hig/\gal$  of Higgs pairs
identifies with the product
\begin{equation*}
\Jac(X) \times  \hH^{1,0}(X).
\end{equation*}
\end{prop}
\noindent
Complex-conjugation maps the space $\hH^{1,0}(X)$ of Higgs fields to 
\begin{equation*}
\hH^{0,1}(X) = \bV. 
\end{equation*}
Thus we identify the space of Higgs fields with
$\bV$ and Proposition~\ref{prop:dolbms} implies that $\Hig/\gal$ is the 
cokernel of
\begin{equation}\label{eqn:HiggsQuotient}
L \longrightarrow V \times \bV,  
\end{equation}
that is, $V/L \times \bV$.  

We denote by $I$ the complex structure on $\Hig/\gal$ arising from
the complex structures on $V$ and $\bV$. In terms of the parameters
$\Psi\in V$ and $\Phi\in \bV$ this complex structure is:
\begin{align}
\Psi &\stackrel{I}\longrightarrow i\Psi \notag \\ 
\Phi &\stackrel{I}\longrightarrow i\Phi. \label{eq:I}
\end{align}
Similarly $T^*\Jac(X)$ is an abelian group 
via a tensor product construction on Higgs pairs. 
Namely define the tensor product 
$(D_1'',\Phi_1)\otimes(D_2'',\Phi_2)$ 
to have holomorphic structure $D_1''\otimes D_2''$ 
and Higgs field $\Phi_1 + \Phi_2$.

\subsection{Geometric structure of the Dolbeault moduli space}
The conformal structure of $X$ induces a rich geometry
to the Dolbeault moduli space $\Hig/\gal$. This space is a complex
manifold under the complex structure $I$ induced from $X$, which is
the total space of a vector bundle over $\Jac(X)$. This vector bundle
is naturally isomorphic to the {\em cotangent bundle\/} $T^*\Jac(X)$,
which has a natural complex-symplectic structure. Furthermore the Riemannian
metric of $\Jac(X)$ defines a proper exhaustion function vanishing on
the zero-section $\Jac(X)$, and whose gradient flow defines a 
deformation-retraction to the zero-section. In rank two, 
Hitchin~\cite{Hitchin1} used a similar technique
to compute the Betti numbers of the Higgs bundle moduli space.

\subsubsection*{Cotangent bundles}
We identify $\Hig/\gal$ with the {\em cotangent bundle\/} $T^*\Jac(X)$
using a Riemannian metric on $\Hig/\gal$ associated to $X$ as follows.

The natural Hermitian form  on $\Aa^1(X)$ defined by 
\begin{equation}\label{eq:hermitianstructure}
\langle \alpha,\beta \rangle :=  \int_X \alpha \wedge \star \bar{\beta} 
\end{equation}
is nondegenerate, positive definite on $\Aa^{0,1}(X)$ and
negative definite on $\Aa^{1,0}(X)$.
Its restriction to $V =\hH^{0,1}(X)$ 
defines an isomorphism $\bV \longrightarrow V^*$ of complex vector spaces. 
The tangent space of $\Jac(X) = V/L$ at any point identifies with $V$.
By the isomorphism defined by $\langle,\rangle$,
the cotangent space of $\Jac(X)$ identifies with $\bV$.
Thus 
%the cotangent bundle of $\Jac(X)$ identifies with
\begin{equation*}
V/L \times \bV \cong  V/L \times V^* \cong  T^*\Jac(X).
\end{equation*}
The real part of $\langle,\rangle$ defines Riemannian metrics on $V$
and $V^*\cong\bV$ which in turn defines a Riemannian metric $g$ 
on $T^*\Jac(X)$.

The cotangent bundle $T^*M$ of any complex manifold $M$ enjoys a
natural complex-symplectic structure $\Omega_{T^*}$, In fact,
$\Omega_{T^*}$ is {\em exact,\/} that is, there exists a holomorphic
$1$-form $\alpha_{T^*}$ on $M$ such that 
$\Omega_{T^*} = -d\alpha_{T^*}$.
Let $\Pi:T^*M\longrightarrow M$ denote projection and
$d\Pi:T(T^*M)\longrightarrow TM$ denote its differential.  Let 
$u\in T^*M$ be a covector at $\Pi(u)\in M$.  
The value of $\alpha_{T^*}$ at $u$ is defined as:
\begin{equation*}
(\alpha_{T^*})_{(u)}: v  \longmapsto  u\big( d\Pi(v) \big)
\end{equation*}
for a tangent vector $v\in T_u(T^*M)$.
Then 
\begin{equation*}
\Omega_{T^*} = -d\alpha_{T^*} 
\end{equation*}
is a complex-symplectic structure
on $T^*M$.  This construction is natural, so any biholomorphism of $M$
preserves $\alpha_{T^*}$ and $\Omega_{T^*}$.

For a complex torus, the complex-symplectic-structure $\Omega_{T^*}$
is parallel. If $\qq = (q_1,\dots,q_n)$ are local coordinates on $M$,
with corresponding coordinates $\pp = (p_1,\dots,p_n)$ on the fibers,
then 
\begin{align}
\alpha_{T^*} &  = \sum_{j=1}^n p_i dq_i \notag \\
\Omega_{T^*} & = \sum_{j=1}^n dq_i\wedge dp_i. \label{eq:pq}
\end{align}
The space $T^*\Jac(X)$ supports a holomorphic completely integrable
Hamiltonian system.  The {\em Hitchin map\/} 
%is the map
\begin{align}
\Hig/\gal & \stackrel{\mathbf{H}}\longrightarrow  \hH^{1,0}(X) \notag \\
[(D'',\Phi)] & \longmapsto \Phi \label{eq:Hitchinmap}
\end{align}
associates to a pair $(D'',\Phi)$ its Higgs field.
In terms of \eqref{eq:pq}, this is the projection 
\begin{equation*}
(\qq,\pp) \longmapsto \pp.
\end{equation*}
The complex Hamiltonian vector field with potential $p_j$ is the
parallel vector field $\frac{\partial}{\partial q_j}$
which generates translation in the $q_j$-coordinate.
The Hitchin map is a holomorphic moment map for a holomorphic Hamiltonian
$\C^k$-action, whose orbits are sections of $T^*\Jac(X)$ corresponding
to the trivialization $\mathbf{H}$, where $k$ is the genus of $X$. 
Equivalently these parallel copies
of $\Jac(X)$ are the preimages of $\mathbf{H}$.  Although in this case
this construction is elementary, its generalization to higher rank is
highly nontrivial and revealing. See Hitchin~\cite{StableIntegrable}.

On $T^*\Jac(X) = V/L \times \bV$, this structure arises 
from the linear complex-symplectic structure  on the vector space
\begin{equation*}
V \oplus \bV = \hH^{0,1}(X) \oplus \hH^{1,0}(X)  
\end{equation*}
defined by 
\begin{equation}\label{eq:OmegaI}
\Omega_{T^*}\big(
(\Psi_1,\Phi_1),(\Psi_2,\Phi_2)\big)
:=  \langle\Psi_1,\overline{\Phi_2}\rangle  - 
\langle\Psi_2,\overline{\Phi_1}\rangle. 
\end{equation}

\subsubsection*{Functions and flows on the Dolbeault moduli space}
The {\em energy function\/}
\begin{align}\label{eq:energy}
T^*\Jac(X) & \stackrel{e}\longrightarrow \R  \notag \\
(\qq,\pp) & \longmapsto \frac12\langle \pp,\pp\rangle
\end{align}
vanishes precisely on the zero-section $\Jac(X)$, which contains
the only critical points of $e$.
These critical points are minima.  Let $\nabla e$ denote the gradient
vector field of $e$ with respect to the Riemannian metric $g$.  Then
the flow of $-\nabla(e)$ defines a deformation retraction of
$T^*\Jac(X)$ to its zero-section $\Jac(X)$, from which one can
determine the homotopy type of the moduli space.
Hitchin~\cite{Hitchin1} uses a similar technique to determine the
topology of some moduli spaces in rank two.  However, in higher rank,
not all critical points are minima, and the calculation is
considerably more complicated. For further applications see
\cite{BGG1,BGG2,Gothen1,Gothen2,Gothen2,Xia0,Xia1,Xia2}.

The vector field $I\nabla(e)$ is the Hamiltonian vector field having potential
$e$ with respect to the \ka~form on $T^*\Jac(X)$. It generates the Hamiltonian
circle action:
\begin{align}\label{eq:hamciract}
S^1  \times T^*\Jac(X) & \longrightarrow T^*\Jac(X) \\
\big(\lambda, (\qq,\pp)\big) & \longmapsto (\qq, \lambda\pp)
\end{align}
which evidently extends to a holomorphic $\C^*$-action on $T^*\Jac(X)$.
The fixed point set of each action equals the zero section
$\pp = 0$, which is precisely the critical set of $e$.
Furthermore the $\C^*$-action extends to an action of the multiplicative
semigroup $\C$; the limit as $\lambda\longrightarrow 0$ lies in the critical
set. We revisit this geometry in \S\ref{sec:functionsflows}.

\section{Equivalence of de Rham and Dolbeault groupoids}
The Betti, de Rham, and Dolbeault deformation theories of the previous 
sections are equivalent.
Hence the three moduli {\em sets\/} are isomorphic.
However these moduli sets have natural differentiable structures
(in fact they are complex Lie groups).
With this extra structure, we refer to the moduli sets as 
{\em moduli spaces.\/}
The equivalence between the de Rham groupoid and the Dolbeault groupoid
induces an isomorphism of real Lie groups  
between their moduli spaces. 

\subsection{Construction of the equivalence}\label{sec:DolDR}
Passage from de Rham to Dolbeault groupoids requires the construction 
of a holomorphic structure and a Higgs field 
for each flat connection $D$.  
To this end, we introduce a Hermitian metric $H$ on $E$
to decompose $D$ into a unitary connection and a $1$-form $\phi$.
The unitary connection determines a unique holomorphic structure $D''$.
The condition that $\phi$ determines a Higgs field
is equivalent to a {\em harmonicity\/} condition on $H$,
which follows from the standard theory of harmonic forms.

\subsubsection*{Hermitian metrics} 
A {\em Hermitian metric\/} $H$ on $E$ 
is a family of positive definite Hermitian forms $\la,\ra_H$ 
\begin{equation*}
E_x \times E_x \longrightarrow \C 
\end{equation*}
on the fibers $E_x$ smoothly varying with $x\in \Sigma$. 
For any vector space $V$, let $\her(V)$ denote the space of positive
definite Hermitian forms $V\times V\longrightarrow\C$.
The space $\Her$ of Hermitian metrics on the vector bundle $E$
is the space of sections of the $\her(\C^r)$-bundle associated to $E$.
In terms of a basis, a Hermitian form is represented by a Hermitian
matrix $h$:
\begin{equation*}
H(u,v) = u^\dag  h \bar{v}
\end{equation*}
where $\dag$ indicates transpose.                       
The action of a linear transformation $g\in\GLr$ on a Hermitian form 
defined by a matrix $h$ is then:
\begin{equation*}
g: h \longmapsto g^\dag h \bar{g}.
\end{equation*}
If $E$ is a flat vector bundle with holonomy $\phi:\pi\longrightarrow\GLr$, 
then a Hermitian metric $H\in\Her$ corresponds to a $\phi$-equivariant map
\begin{equation}\label{eq:positivefunction}
\tS  \stackrel{h}\longrightarrow \her(\C^r).
\end{equation}
In particular $\End(E)$ is a complex line bundle with a canonical
{\em everywhere  nonzero\/} section, the identity endomorphism
$\Id:E_x\longrightarrow E_x$. Hence $\End(E)$ identifies canonically
with the trivial line bundle $\C$ over $X$.

The associated bundles of Hermitian metrics are trivial, 
but not canonically trivial.

However, in the important special case when the
holonomy lies in the subgroup $\Uo$, the standard Hermitian form on $\C$
\begin{equation*}
\la z,w\ra := z\bar{w} 
\end{equation*}
defines a Hermitian metric on $E$. 
Gauge transformations act on Hermitian metrics by:
\begin{equation*}
\la u,v\ra_{\xi\cdot H} :=  \la g^{-1}\cdot u,g^{-1}\cdot v\ra_{H}
\end{equation*}
and, in rank one, the corresponding positive function $h$ transforms by:
\begin{equation*}
h \longmapsto \vert g\vert^{-2} h.
\end{equation*}
Two Hermitian metrics $H_1,H_2\in\Her$ on a complex line bundle
$E$ relate by their {\em ratio\/}, the positive function 
$h:X\longrightarrow\R^+$ defined by:
\begin{equation*}
\la u,v \ra_{H_1} = h(x) \la u,v \ra_{H_2}
\end{equation*}
where $u,v\in E_x$. 
\begin{lem} \label{lem:transitive}
The linear gauge group $\gal$ acts transitively on the space $\Her$ of
Hermitian metrics.
\end{lem}
\begin{proof}
If $h_1,h_2$ are positive functions representing Hermitian metrics on
$E$, then
\begin{equation*}
g(z) = \sqrt{h_1(z)/h_2(z)}
\end{equation*}
defines a gauge transformation $g\in\gal$ with
$g\cdot h_2 = h_1$.
\end{proof}
The unique Hermitian metric $H_0$ for which 
\begin{equation*}
\la \tau,\tau\ra_{H_0}=1  
\end{equation*}
corresponds to the constant function $1$ and its stabilizer in $\gal$
is the {\em unitary gauge group\/} $\gau$, corresponding to functions
$X\longrightarrow \Uo$.

A Hermitian metric $H$ on $E$ induces Hermitian pairings over $\Aa^*(X)$
\begin{equation*}
\Aa^k(X;E) \times \Aa^l(X;E) \longrightarrow \Aa^{k+l}(X).
\end{equation*}
A connection $D$ is {\sl unitary\/} with respect to a metric $H$
if and only if
\begin{equation*}
d \la s_1,s_2\ra_H =  \la Ds_1,s_2\ra_H + \la s_1,Ds_2\ra_H,
\end{equation*}
for sections $s_1, s_2 \in \Aa^0(X,E)$.
Equivalently we say that $H$ is {\em parallel\/} with respect to $D$.
With respect to $\tau_0$, the connection 
$D = D_0 + \eta$ is unitary with respect to  $H$ if and only if:
\begin{equation} \label{eq:exactRealPart}
\eta + \bar{\eta} = 2\Re(\eta) = h^{-1}dh = d\log h.
\end{equation}
Thus $D = D_0 + \eta \in \fcl$ is unitary with respect to some
Hermitian metric on $E$ if and only if the 1-form $\Re(\eta)$ is exact.

In particular a connection $D = D_0 + \eta$ is unitary with respect to the
Hermitian metric $H_0$ corresponding to the trivialization $\tau_0$
(that is, the Hermitian metric with $h \equiv 1$) if and only if
$\eta\in\Aa^1(X;i\R)$.                                
We denote the subset of $\fcl$ consisting of flat connections unitary
with respect to $H_0$ by $\fcu$.                           
Since $\gal$ acts transitively on $\Her$,
a connection $D \in \fcl$ is $\gal$-equivalent to a connection
in $\fcu$ if and only if $D$ is unitary with respect to some Hermitian
metric on $E$.

The following well known result relates holomorphic structures, Hermitian
structures and unitary connections:
\begin{prop}\label{prop:uniqueconnection}
Let $E$ be a complex line bundle over
$X$ with Hermitian metric $H$, and let $D''$
be a holomorphic structure on $E$. Then
there exists a unique connection $D$ on $E$ such that
\begin{itemize}
\item The $(0,1)$-part $D^{0,1}$ of $D$ equals $D''$;
\item $H$ is parallel with respect to $D$.
\end{itemize}
\end{prop}
\noindent We call 
$D$ the connection {\em compatible\/} with $D''$ and $H$.

\begin{proof}
Write $D'' = D_0'' + \mu$, where $\mu\in\Aa^{0,1}(X)$. 
Let $h$ be the equivariant positive function on $\tilde X$ 
corresponding to $H$ as in  \eqref{eq:positivefunction}.
Let $D = D_0 + \eta$ where
\begin{equation*}
\eta = \mu - \bar{\mu} + h^{-1}\d h.
\end{equation*}
Since $\eta + \bar\eta = h^{-1}dh$, \eqref{eq:exactRealPart} implies
that $D$ is unitary with respect to $H$. Since 
$\eta^{0,1} = \mu$, the holomorphic structure defined by $D$ equals
$D''$ as desired.
 
For uniqueness, let $D_1,D_2$ be two $\Uo$-connections
with the same $(0,1)$-part. 
By Lemma~\ref{lem:transitive}, we may assume $H=H_0$.
With respect to the trivialization $\tau_0$,
write $D_i = D_0 + \psi_i$ where $D_0$ is the trivial connection. Then
$\psi_1$ and $\psi_2$ are 1-forms with values in $i\R$.
Their difference $\psi_1-\psi_2$ has values in $i\R$ with Hodge type (1,0).
Since a purely imaginary 1-form of Hodge type (1,0) is necessarily zero,
$\psi_1-\psi_2$ is zero.
\end{proof}
\noindent
(Although we only proved Lemma~\ref{lem:transitive} and
Proposition~\ref{prop:uniqueconnection} for line bundles, these
results hold for vector bundles of arbitrary rank.)

\subsubsection*{The unitary case}
Suppose $D$ is a connection. Its composition
\begin{equation*}
\Aa^0(X;E) \stackrel{D}\longrightarrow
\Aa^1(X;E) \longrightarrow  \Aa^{0,1}(X;E)
\end{equation*}
with $(0,1)$-projection is a holomorphic structure, denoted $D^{0,1}$.  If 
\begin{equation*}
D = D_0 + \eta,
\end{equation*}
then the corresponding holomorphic structure is
\begin{equation*}
D^{0,1} = D_0'' + \Psi
\end{equation*}
where $\Psi = \eta^{0,1}$.
 
The correspondence between holomorphic structures on $E$ and flat
unitary connections is almost an equivalence of deformation theories.
Define a functor
\begin{equation*}
\big(\fcu,\gau\big) \stackrel{\Ss}\longrightarrow \big(\Hol,\gal\big),
\end{equation*}
which on objects is the 
map assigning to a flat $\Uo$-connection $D$ on $E$
the holomorphic structure $D^{0,1}$, and on morphisms
is the inclusion $\gau\hookrightarrow\gal$.

Unfortunately $\Ss$ is not full.
Consider the connection $D_0$ on $E$ and $D_0''$ the trivial
holomorphic structure. Then
\begin{equation*}
\Mor(D_0,D_0) \longrightarrow  \Mor(D_0'',D_0'')
\end{equation*}
is not surjective:
$\Mor(D_0,D_0)$ consists of scalar multiplication by unit complex numbers,
while $\Mor(D_0'',D_0'')$ corresponds to scalar multiplication by
$\C^*$.

To overcome this difficulty, rigidify
the groupoid $(\Hol,\gal)$ by including a
Hermitian metric as an extra piece of data.
Therefore $\gal$ acts on $\Hol \times \Her$.

Suppose $D'' \in \Hol$ and $H\in\Her$.
Let $D$ be the unique connection compatible with $D''$
and $H$ (Proposition~\ref{prop:uniqueconnection}).

\begin{defin}\label{def:adapted}
A Hermitian metric $H$ is {\em adapted\/}
to $D  ''$ if and only if 
the connection $D$ compatible with $D''$ and $H$ is flat, that is, 
$F(D)=0$. 
\end{defin}
(In Corlette~\cite{Corlette2}, such metrics are called {\em
harmonic,\/} with respect to a holomorphic line bundle. To avoid
confusion, we chose to restrict the terminology ``harmonic metric''
to the context of flat bundles and refer to the corresponding notion
for holomorphic bundles as ``adapted metric''.)

Let $\Holu$ be the subset of $\Hol\times \Her$ consisting of all
$(D'',H)$ such that $H$ is adapted to $D''$ in the sense of
Definition~\ref{def:adapted}.
Then $\gal$ acts on $\Holu$, and
we denote the corresponding deformation theory $(\Holu, \gal)$.
\begin{prop} \label{prop:equiv}
The natural map
\begin{align*}
\fcu & \longrightarrow  \Holu \\
D & \longmapsto  (D^{0,1},1)
\end{align*}
is equivariant with respect to the natural $\gal$-action. The
corresponding functor 
\begin{equation*}
(\fcuu,\gau) \stackrel{\mathcal T}\longrightarrow (\Holu,\gal)
\end{equation*}
is an equivalence of deformation theories.
\end{prop}
\begin{proof}
Let $(D'', h) \in \Holu$.
By applying a linear gauge transformation, we may assume that
$h = 1$.  
Apply Proposition~\ref{prop:uniqueconnection}, to obtain a
$\Uo$-connection $D$ satisfies $D^{0,1} = D''$.
By definition of $\Holu$, the connection $D$ is unitary, whence
$D\in\fcu$. Thus $\mathcal{T}_*$ is surjective.

Since the subgroup of $\gal$ preserving $H_0$ is 
$\gau$, the induced map
\begin{equation*}
\Mor(D_1,D_2) \stackrel{\mathcal{T}(D_1,D_2)}
\longrightarrow  \Mor(D_1'',D_2'')
\end{equation*}
is an isomorphism.
Thus $\mathcal{T}$ is faithful and full.
\end{proof}

\begin{prop} \label{prop:adaptedmetric}
For each $D'' \in \Hol$, there exists $h \in \Her$ such
that $(D'', h) \in \Holu$.
\end{prop}
\begin{proof}
Let $D'' \in \Hol$.
One must find a flat unitary connection which is equivalent
by a gauge transformation in $\gal$ to a connection whose
$(0,1)$-part equals $D''$.
Write $D'' = D_0'' + \Psi$ where $\Psi\in\Aa^{0,1}(X)$.
Now 
\begin{equation*}
\Aa^{0,1}(X) = \hH^{0,1}(X) \oplus \db\Aa^0(X)
\end{equation*}
and $\Psi = \Psi_0 + \db s$, where $\Psi_0$ is antiholomorphic.
(\S\ref{sec:dolbgpd}). 
In particular $\Psi_0$ is closed.
Let 
\begin{equation*}
D = D_0 + \Psi_0 - \bar{\Psi}_0 
\end{equation*}
and $g=\exp(s)$. Then $D$ is flat
(since $d\Psi_0 = d\bar{\Psi}_0 = 0$), unitary with respect to $1$
(since $\Psi_0 -\bar{\Psi}_0$ is purely imaginary).
Furthermore the gauge transformation
$g$ takes 
\begin{equation*}
D^{0,1} = D_0'' + \Psi_0 
\end{equation*}
to $D''$. The metric
$H_0$ corresponds to the function $\bar{g} g$.  
The desired adapted metric is $g\cdot H_0$.
\end{proof}
Together, Propositions~\ref{prop:equiv} and \ref{prop:adaptedmetric} imply:
\begin{cor} \label{cor:u1}
The induced map
\begin{equation*}
\fcu/\gau \longrightarrow  \Hol/\gal
\end{equation*}
is an isomorphism of sets.
\end{cor}
That every topologically trivial holomorphic line bundle over $X$
arises from a unitary character of $\pi_1(X)$ is a standard result;
see Weyl\cite{Weyl}, \S 18, 
Gunning~\cite{Gunning1}, \S 8a, pp.\ 131--135, or
Farkas-Kra~\cite{FarkasKra}, \S III.9, pp.\ 119--129. In the
terminology of \cite{FarkasKra}, a {\em character\/} is a homomorphism
$\pi_1(X)\longrightarrow\C^*$ and a {\em normalized character\/} is a
unitary character. A {\em multiplicative function\/}
in the sense of \cite{Weyl,FarkasKra}  corresponds to a
holomorphic section of a flat complex line bundle.

\subsubsection*{The linear case}
Now we extend the preceding theory from $\Uo$-rep\-re\-sen\-ta\-tions 
and flat unitary connections to $\C^*$-representations and flat linear
connections.

Suppose $D \in \fcl$.  With respect to $\tau_0$,
write $D = D_0 + \eta$, where $\eta \in \zz^1(X)$.
Let $\phi, \psi$ be the real and imaginary parts of $\eta$ respectively.
Since $D$ is flat,
\begin{equation*}
d \phi + d \psi = 0.
\end{equation*}
Since $d \psi$
is purely imaginary and $d \phi$ is real, this single equation
is equivalent to the pair of equations: 
\begin{equation*}
d \phi = d \psi = 0.
\end{equation*}
Thus $D$ is flat if and only if  $\phi$ and $\psi$ are closed.

For any Hermitian metric $H$ on $E$, decompose the connection $D$ into
its skew-adjoint part (a connection $D_H$ compatible with $H$) and its
self-adjoint part, a 1-form $\phi_1$ which is self-adjoint with
respect to $H$:
\begin{equation}
D = D_H + \phi_1,
\end{equation}
where 
\begin{equation*}
D_H = D_0 + \psi + \frac12 h^{-1}dh  
\end{equation*}
and 
\begin{equation*}
\phi_1 = \phi - \frac12 h^{-1}dh 
\end{equation*}
where $h$ is the positive function corresponding to $H$
as in \S\ref{sec:DolDR}.

Then $D_H$ is flat ($h^{-1}dh$ is exact) and 
compatible with $H$ (by \eqref{eq:exactRealPart}).
Since $D_H$ is flat, $\phi_1$ is a closed 1-form.
Let $\Phi = (\phi_1)^{1,0}$.

\subsubsection*{Harmonic metrics.}

We also recall basic facts about harmonic functions.
A smooth function $f$ is {\em harmonic\/}
if and only if $d\star df = 0$.
In particular holomorphic functions are harmonic.
A harmonic 1-form locally is the differential
of a harmonic function. Every harmonic function is locally the
real part of a holomorphic function.
A smooth function  $f:X\longrightarrow\R^+$ is 
{\em multiplicatively harmonic\/} if
and only if its logarithm 
\begin{equation*}
\log f:\tX\longrightarrow\R  
\end{equation*}
is a harmonic function.

\begin{defin} Let $D$ be a flat connection on the line bundle $E$ over $X$. 
A Hermitian metric $H\in\her(E)$ 
is {\em harmonic with respect to $D$\/} if and only if
the equivariant map
\begin{equation*}
\tilde X \stackrel{\tilde h}\longrightarrow \R^+    
\end{equation*}
corresponding to $H$ is a multiplicatively harmonic function on $\tilde X$,
that is, its logarithm $\log \tilde h$ is harmonic function
$X\longrightarrow\R$.
\end{defin}
This definition is independent of the trivialization of $D$ over the 
universal covering space $\tilde X\longrightarrow X$ used to define
$\tilde{h}$.

$H$ is a harmonic metric with respect to $D$ if and only if
the 1-form $\phi_1 = \phi - \frac12 h^{-1}dh$ is harmonic,
or equivalently $\db \Phi = 0$.

\begin{lem}\label{lem:harmonicmetric}
Let $D$ be a flat connection. Then there
exists a  Hermitian metric $H$ harmonic with respect to $D$.
Furthermore $H$ is unique up to multiplication by a positive constant.
\end{lem}
\begin{proof}
Write $D = D_0 + \psi + \phi$ as above, where $\phi\in\zz^1(X,\R)$. 
Then $\phi$ is cohomologous
to a unique harmonic form $\phi_1\in\hH^1(X;\R)$, that is,
\begin{equation*}
\phi = \phi_1 + df.
\end{equation*}
for some smooth function $f\in\Aa^0(X, \R)$  (\S\ref{sec:dolbgpd}).
The desired harmonic metric is $h = e^{2f}$.
Since $f$ is unique up to addition of a real constant, the metric
$H$ is unique up to multiplication by a positive constant.
\end{proof}

In higher rank, these facts generalize to deeper
technical results. The existence and uniqueness
of flat unitary connections on a topologically trivial
holomorphic line bundle (Proposition~\ref{prop:adaptedmetric})
generalizes to Donaldson's proof of
the Narasimhan-Seshadri theorem \cite{Donaldson1},  
Hitchin's proof of solutions of the self-duality 
equations in rank two \cite{Hitchin1}
and Simpson's proof in general ~\cite{Simpson1,Simpson2}. 
Lemma~\ref{lem:harmonicmetric}
generalizes to the existence of harmonic metrics on flat bundles with
reductive holonomy (Donaldson~\cite{Donaldson2} in rank two, and
Corlette~\cite{Corlette1} in general). The existence of a harmonic
metric generalizes the existence and uniqueness of harmonic maps 
of compact Riemannian manifolds into 
complete Riemannian manifolds of nonnegative curvature, due to 
Eells-Sampson~\cite{EellsSampson}, and has been further generalized
by Labourie~\cite{Labourie}).

Let $D \in \fcl$ and $H$ its harmonic metric and $D = D_H + \phi_1$,
where $D_H$ is compatible with $H$.
Let $D''$ be the (0,1)-part of $D_H$
and $\Phi$ be the (1,0)-part of $\phi_1$.
Since $\phi_1$ is harmonic, $\db \Phi = 0$. 
Thus $(E, D'', \Phi)$ is a Higgs pair.

\begin{thm} \label{thm:c}
The resulting functor
\begin{equation*}
(\fcl,\gal) \stackrel{\Ss}\longrightarrow (\Hig,\gal).
\end{equation*}
is an equivalence of deformation theories.
\end{thm}
\begin{proof}
First we show that $\Ss$ is surjective on isomorphism
classes.  Let $(D'',\Phi) \in \Hig$.  By 
Proposition~\ref{prop:adaptedmetric}, $D''$ has an adapted metric $H$
corresponding to a function 
$h : X \rightarrow \R^+$.  Let $D_H$ be the unique connection
compatible with $H$ and $(D_H)^{0,1} = D''$
(Proposition~\ref{prop:uniqueconnection}).
Write $D_H = D_0 + \psi$.  
Then
$\xi \cdot H = H_0$ for some linear gauge transformation $\xi\in\gal$.
Let 
\begin{equation*}
D = \xi \cdot (D_H + \Phi + \bar{\Phi}) = 
D_0 +  \Phi + \bar{\Phi}.
\end{equation*}
The curvature of $D$ is
\begin{equation*}
F(D) = F(\xi \cdot D_H) + d (\Phi + \bar{\Phi}).
\end{equation*}
Since $H$ is adapted to $D''$, the curvature $F(D_H) = 0$, so
$F(\xi \cdot D_H) = 0$.  Since $\Phi$ is holomorphic,
\begin{equation*}
d (\Phi + \bar{\Phi}) = 0. 
\end{equation*}
Hence $F(D) = 0$ and $D \in \fcl$.
Then $\Ss(D)$ is a Higgs pair 
differing from $(D'',\Phi)$ by the gauge transformation $\xi$.

Next we show $\Ss$ is faithful and full. Suppose 
\begin{equation*}
(D_1'', \Phi_1), (D_2'',\Phi_2) \in \Hig.
\end{equation*}
There are two cases, depending on whether or not
$(D_1'', \Phi_1), (D_2'',\Phi_2) $ are 
$\gal$-equivalent.
If $(D_1'', \Phi_1), (D_2'',\Phi_2) $ are not
$\gal$-equivalent, then 
\begin{equation*}
\Mor\big((D_1'', \Phi_1), (D_2'',\Phi_2)\big) 
\end{equation*}
is empty.
Otherwise, there exists $\xi \in \gal$ corresponding to a map 
$g\in\Map(X,\C^*)$ such that
\begin{equation*}
D_2'' = D_1'' +  g^{-1}\db g, \ \ \ \Phi_2 = \Phi_1
\end{equation*}
and if $\xi_1$ is another gauge transformation corresponding to a map
$g_1\in\Map(X,\C^*)$ such that
\begin{equation*}
D_2'' = D_1'' + g_1^{-1} \db g_1,
\end{equation*}
then $g_1 = g c$, where $c$ is a constant map. 
We denote the subgroup of $\gal$ corresponding to constant maps
$X\longrightarrow \C^*$ by $\C^*$.
Hence
$\Mor\big((D_1'', \Phi_1), (D_2'',\Phi_2)\big)$ corresponds to the 
coset $g\cdot\C^*$.

Now suppose that $D_1,D_2\in\fcl$ and $\Ss(D_i)=(D_i'',\Phi_i)$ for 
$i=1,2$. If $D_1$ is not $\gal$-equivalent to $D_2$, then 
then $\Mor(D_1,D_2)=\emptyset$ since $\C^*$ acts trivially on $\fcl$.
Otherwise, $D_1 = g \cdot D_2$ where $g\in\gal$. Then 
$\Mor(D_1,D_2) = g\cdot\C^*$.

In both cases, $\Ss$ induces an isomorphism
\begin{equation*}
\Mor(D_1, D_2) \longrightarrow \Mor(\Ss(D_1), \Ss(D_2))
\end{equation*}
as desired. 
\end{proof}

\begin{cor} \label{cor:cstar}
The induced map
\begin{equation*}
\fcl/\gal \stackrel{\Ss_*}\longrightarrow  \Hig/\gal
\end{equation*}
is an isomorphism of sets.
\end{cor}

\subsection{Higgs coordinates}
To see that the bijection $\Ss_*$ between the de Rham and 
Dolbeault moduli spaces 
is an isomorphism of real Lie groups, 
we examine in detail the constructions in the proof
of Theorem~\ref{thm:c}.  We explicitly describe how to pass from the
equivalence class of a flat connection $D = D_0 + \eta$ to an equivalence
class of a Higgs pair $(D'',\Phi)$, where $D'' = D''_0 + \Psi$.
We call $\eta\in\hH^1(X)$ 
the {\em connection coordinates\/} of the point in $\fcl/\gal$ 
and $(\Psi,\Phi)\in V\oplus\bV$ the {\em Higgs coordinates.\/}

Passing from connection coordinates $\eta\in\hH(X)$ to Higgs
coordinates $(\Psi,\Phi) \in V \oplus \bV$ involves the
following steps.

We assume that $\eta$ is harmonic, since we
are only interested in the equivalence class of $D$ under $\gal$.
Decompose $\eta$ into its imaginary part $\psi$ and its real part
$\phi$:
\begin{equation*}
\eta = \psi + \phi 
\end{equation*}
where
\begin{equation*}
\psi = \frac12 (\eta - \bar\eta), \qquad  
\phi = \frac12 (\eta + \bar\eta).
\end{equation*}
The unitary connection $D_H := D_0 + \psi$
corresponds to a holomorphic structure $D'' = D_0'' + \Psi$ where the
antiholomorphic $1$-form $\Psi$ is the $(0,1)$-part of $\psi$:
\begin{align*}
\Psi & =  \psi^{0,1} \\
\psi & = \Psi - \bar\Psi.
\end{align*}
The $(1,0)$-part $\Phi$ of $\phi$ is a holomorphic $1$-form:
\begin{align*}
\Phi & =  \phi^{1,0} \\
\phi & = \Phi + \bar\Phi.
\end{align*}
In summary, to pass from Higgs coordinates to connection coordinates:
\begin{equation}\label{eq:etaPsiPhi}
\eta = \Psi - \bar\Psi + \Phi + \bar\Phi
\end{equation}
and to pass from connection coordinates to Higgs coordinates:
\begin{align}
\Psi & = (i \Im \eta)^{0,1} \notag\\
\Phi & = (\Re \eta)^{1,0} \label{eq:cc2hc}
\end{align}

\subsubsection*{Holomorphic structures associated to flat connections}
Since
\begin{align*}
\eta  & = \psi + \phi \\
& = (\Psi - \bar{\Psi} ) + (\Phi +\bar{\Phi})  \\
& = (\Psi + \bar{\Phi} ) + (\Phi -\bar{\Psi})
\end{align*}
where $\Psi\in\hH^{0,1}(X)$ and $\Phi\in\hH^{1,0}(X)$, the Hodge components
of $\eta$ are:
\begin{align*}
\eta^{0,1} & =  \Psi + \bar{\Phi} \\
\eta^{1,0} & =  \Phi - \bar{\Psi}.
\end{align*}
The flat connection $D$ determines a holomorphic structure:
\begin{equation*}
D^{0,1} :=  D_0'' + \eta^{0,1} = D_0'' + \Psi + \bar{\Phi} = D'' + \bar{\Phi}.
\end{equation*}
The holomorphic structure $D''$ in the Higgs pair 
$(D'',\Phi)$ is:
\begin{equation*}
D'':=  D_0'' +  \Psi. 
\end{equation*}
Thus, unless $\Phi = 0$, the holomorphic structure in the Higgs pair is
{\em not\/} the holomorphic structure determined by the flat connection.

\subsubsection*{Isomorphism of Lie groups}
Now we strengthen Corollary~\ref{cor:cstar}:
\begin{prop}
The induced map
\begin{equation*}
\fcl/\gal \stackrel{\Ss_*}\longrightarrow  \Hig/\gal
\end{equation*}
is an isomorphism of real Lie groups
\end{prop}
\begin{proof}
Using harmonic decomposition, the moduli space 
$\fcl/\gal$ identifies with the quotient of $\hH^1(X)$
by the discrete subgroup corresponding to the subgroup
$\H^1(X,\Z)\hookrightarrow \H^1(X)$. 
The moduli space $\Hig/\gal$ corresponds to the quotient of
$V \oplus\bV$ by a lattice $L\subset V$ in the first factor. 
The map $\Ss_*$ arises from the $\R$-linear map
\begin{align*}
\hH^1(X) & \longrightarrow V \oplus \bV \\
\eta &\longmapsto (\Psi,\Phi)
\end{align*}
in Higgs coordinates.
The induced bijection $\Ss_*$ is an isomorphism of real Lie groups.
\end{proof}

\subsubsection*{Real structures}
The anti-involutions defining the real forms $\Uo$ and $\R^*$
were described in connection coordinates in \eqref{eq:iotacc}.
Using \eqref{eq:etaPsiPhi} and \eqref{eq:cc2hc} to pass to 
Higgs coordinates, $(\iota_U)_*$ and $(\iota_{\R})_*$ are:

\begin{align}
(\Psi,\Phi)  & \xmapsto{(\iota_U)_*} (\Psi,-\Phi)\notag \\
(\Psi,\Phi)  & \xmapsto{(\iota_\R)_*}(-\Psi,\Phi). \label{eq:iotahc}
\end{align}
The fixed point set of $(\iota_U)_*$ equals $\fcu/\gau$.
In the description of the moduli space as $T^*\Jac(X)$,
the involution $\iota_U$ 
is scalar multiplication by $-1$ on the fibers of 
$T^*\Jac(X)\longrightarrow\Jac(X)$. 

In contrast, $\iota_\R$ leaves $\Phi$ unchanged, 
but replaces the holomorphic structure $D''$
by its inverse.
In the cotangent bundle description, 
$\iota_\R$ is the lift of inversion 
\begin{align*}
\Jac(X)& \longrightarrow\Jac(X) \\
L & \longmapsto L^{-1}
\end{align*}
to $T^*\Jac(X)$
which fixes the parallelism 
\begin{equation*}
T^*\Jac(X) \longrightarrow  T_{L_0}^*\Jac(X) \cong 
V^* \cong \bV. 
\end{equation*}
The differential of inversion induces the composition of this map with
the map $\iota_U$ above given by scalar multiplication by $-1$ on the
cotangent fibers. 

The fixed points of inversion form the 2-torsion subgroup
\begin{equation*}
\Jac_2(X) \hookrightarrow \Jac(X)
\end{equation*}
corresponding to
\begin{equation*}
\bigg(\frac12 L\bigg)/L \subset V/L 
\end{equation*}
which is isomorphic to $(\Z/2)^{2k}$ where $k$ is the
genus of $X$.
Thus the fixed point set of $\iota_\R$ acting on 
$\Hom(\pi,\C^*)$ is the subset 
$\Hom(\pi,\R^*)$ of
\begin{equation*}
\Hom(\pi,\C^*)\longleftrightarrow T^*\Jac(X) 
\end{equation*}
corresponding to
\begin{equation*}
T^*_{\Jac(X)}\big|_{\Jac_2(X)} \cong
(\Z/2)^{2k} \times \bV. 
\end{equation*}

\subsection{Involutions}
The involution $\iota_U$ is holomorphic with respect to $I$, 
but anti-holomorphic with respect to $J$.
The fixed point set of $t = -1$ in the $\C^*$-action is holomorphic with
respect to $I$, but is a purely real subset with respect to $J$.  
Hence this fixed point set is $\Hom(\pi, \Uo)$ which embeds into
$\Hig/\gal$ as a
holomorphic submanifold by \eqref{eq:iotahc}.  

Similarly, the involution $\iota_\R$ is holomorphic with respect
to $I$ and antiholomorphic with respect to $J$.
The fixed point set of $\iota_\R$ is $\Hom(\pi,\R^*)$ 
which is holomorphic with respect to $I$.

The common fixed point set of $\iota_U$ and $\iota_\R$ 
is the intersection
\begin{equation*}
\Hom(\pi,\Uo) \,  \cap \, \Hom(\pi,\R^*) \, = \, \Hom(\pi,\pm\Id) 
\end{equation*}
which identifies with the 2-torsion subgroup
$\Jac_2(X)$ of the Jacobian. 
(This is also the fixed point set of the composition
$\iota_U\circ\iota_\R$.)
In particular $\Hom(\pi,\R^*)$ corresponds
to Higgs pairs $(D'',\Phi)$ where the holomorphic line bundle $(E,D'')$
has order two in the Jacobian. Its identity component $\Hom(\pi,\R^+)$
corresponds to Higgs bundles whose underlying holomorphic structure is
trivial.

The holomorphic structure of the Higgs pair corresponding to an
$\R^+$-representation is trivial.  If $\rho:\pi\longrightarrow\R^+$ is
a nontrivial real character, then there exists a global holomorphic section
to the corresponding flat line bundle but {\em no \/} global parallel
section. This holomorphic section arises as follows. Under the
identifications
\begin{equation*}
\Hom(\pi,\R^+) \xrightarrow{~\log~}
\Hom(\pi,\R) \; \cong \; \H^1(X;\R)
\end{equation*}
$\rho$ corresponds to a harmonic 1-form $\phi$. On the
universal cover, $\phi$ lifts to an exact 1-form which is
the differential of a harmonic function $f:\tX\longrightarrow \R$
satisfying:
\begin{equation}\label{eq:harmonicfunction}
f\circ\gamma - f = \log\rho(\gamma).
\end{equation}
There exists a holomorphic function  $F:\tX \longrightarrow \C$
with $\Re F = f$ which necessarily satisfies \eqref{eq:harmonicfunction}.
(The function $F$ is unique up to an additive purely imaginary constant.)
Then $\exp F$ is a holomorphic function on $\tX$ and defines a 
nonvanishing holomorphic section of the corresponding flat line bundle.

This holomorphic section is nonvanishing for the following general reason.
Since every flat line bundle has degree zero, a nonzero holomorphic
section has no zeroes. Compare the discussion in
Gunning~\cite{Gunning2}, pp.237--238.
The function $f$ arose earlier in the construction of the harmonic
metric in Lemma~\ref{lem:harmonicmetric}.

In general the direct product decomposition of Lie groups 
$\C^* \cong \Uo \times \R^+$
induces a decomposition of Betti moduli spaces
\begin{equation*}
\Hom(\pi,\C^*) \cong \Hom(\pi,\Uo) \times \Hom(\pi,\R^+). 
\end{equation*}
The corresponding decomposition of Dolbeault moduli spaces
\begin{equation*}
\Hig/\gal  \cong \Jac(X) \times \H^1(X;\R)
\end{equation*}
associates to a Higgs pair the underlying holomorphic structure
and the class of the
logarithmic differential of the
harmonic metric in $\H^1(X;\R)$.

\section{Hyperk\"{a}hler geometry on the moduli space}\label{sec:geomdefsp}

Theorem~\ref{thm:one1} and Corollaries~\ref{cor:u1} and
\ref{cor:cstar} establish bijections between the three moduli spaces
\begin{equation*}
(\C^*)^{2k}, \quad  \H^1(\Sigma)/\H^1(\Sigma,\Z), \quad  T^*\Jac(X) 
\end{equation*}
which are isomorphisms of real Lie groups. 
The first two --- Betti and de Rham --- 
were isomorphic as complex Lie groups with isomorphic
complex-symplectic structures.
We henceforth will only work with the de Rham 
and Dolbeault moduli spaces.

The complex structure $J$ on the de Rham moduli space
and the complex structure $I$ on the Dolbeault space are quite
different. Using the isomorphism $\Ss_*$, we superimpose these
structures on a single object which we call $\M$.
The interaction between these structures leads to a rich
{\em \hk geometry \/} and group actions on $\M$.
 
\subsection{The quaternionic structure}\label{sec:quaternion}
In connection coordinates $\eta\in\hH^1(X)$ and Higgs coordinates
$(\Psi,\Phi)\in V\oplus \bV$, 
the complex structure $I$ on $\M$ is:
\begin{align}\label{eq:II}
\eta & \stackrel{I}\longmapsto -\star\bar{\eta}\notag\\
(\Psi,\Phi) & \longmapsto (i\Psi,i\Phi). 
\end{align}
The formula in connection coordinates may be derived as follows.
Let $(\Psi,\Phi)\in V \oplus \bV$ be the Higgs coordinates of $\eta$. 
By definition of $I$, the Higgs coordinates of $I(\eta)$ are 
$(i\Psi,i\Phi)$ so \eqref{eq:etaPsiPhi} implies: 
\begin{align*}
I(\eta) & = (i\Psi) -\overline{(i\Psi)} + (i\Phi) +
\overline{(i\Phi)} \\
& = i\Psi + i\bar\Psi + i\Phi - i\bar\Phi \\
& = \star\Psi  -\star\bar\Psi -\star\Phi  -\star\bar\Phi \\
& = -\star\overline{(-\bar\Psi + \Psi + \bar\Phi + \Phi)} \\
& = -\star\bar\eta. 
\end{align*}
Using the equivalence $\Ss_*$ with the de Rham moduli space,
the natural complex structure on the de Rham moduli space \eqref{eq:J}
defines the complex structure $J$ on $\M$: 
\begin{align}\label{eq:JJ}
\eta & \stackrel{J}\longmapsto i\eta \notag \\
(\Psi,\Phi) & \longmapsto (i\bar{\Phi}, -i\bar{\Psi}). 
\end{align}
To derive the expression in Higgs coordinates, decompose $\eta$ as in \eqref{eq:etaPsiPhi}. Then:
\begin{align*}
J(\eta) = i\eta &= i \Psi - i \bar{\Psi} + i \Phi + i\bar{\Phi} \\
& = (i\bar\Phi) - \overline{(i\bar\Phi)} + (-i\bar\Psi) + 
\overline{(-i\bar\Psi)}
\end{align*}
so that the Higgs coordinates of 
$J(\eta)$ are $(i\bar{\Phi},-i\bar{\Psi})$.
Let $K := IJ$. In coordinates:
\begin{align}\label{eq:K}
\eta & \stackrel{K}\longmapsto i\star\bar\eta \\
(\Psi,\Phi) & \longmapsto (-\bar{\Phi}, \bar{\Psi}). \notag
\end{align}

The complex structures $I,J$ satisfy
\begin{equation*}
I^2 = J^2 = -1, \qquad IJ = -JI
\end{equation*}
and together $I,J,K$ define a left action of the 
{\em algebra of quaternions\/} 
$\bH$ on $\hH^1(X)$. Here $\bH$ is the $\R$-algebra with basis
$1,I,J,K$ and multiplication satisfying: 
\begin{equation}\label{eq:quaternionidentities}
I^2 = J^2 = K^2  = -1  
\end{equation}
\begin{equation*}
IJ = -JI  = K, \quad JK = -KJ  = I, \quad
KI = -IK  = J. 
\end{equation*}
Furthermore every purely imaginary unit quaternion
$u = u_I I + u_J J + u_K$
(where $u_I,u_J,u_K\in\R$ satisfy $u_I^2 + u_J^2 + u_K^2 = 1$)
satisfies $u^2 = -1$, obtaining a whole $S^2$ of complex structures.
Opposite complex structures are related by {\em the antipodal map:\/}
\begin{align}
S^2 & \longrightarrow S^2  \notag\\ 
u & \longmapsto - u. \label{eq:antipodal}
\end{align}

\subsection{The Riemannian metric}
The real part of the Hermitian form $\langle,\rangle$ on $V$ in
\eqref{eq:hermitianstructure} and the real part of the corresponding 
Hermitian form on $\bV$ define a Riemannian metric $g$ on $\M$:
\begin{align}\label{eq:g}
g\big((\Psi_1,\Phi_1),(\Psi_2,\Phi_2)\big) & =
\Re \langle \Psi_1,\Psi_2\rangle +
\Re \langle \Phi_1,\Phi_2\rangle \notag \\
g(\eta_1,\eta_2) & = \int_X \eta_1 \wedge \,\star\, \overline{~\eta_2}. 
\end{align}
The Riemannian metric is invariant under the 
complex structures $I,J,K$:
\begin{equation*}
g(\eta_1,\eta_2) =  
g(I\eta_1,I\eta_2) =  
g(J\eta_1,J\eta_2) =  
g(K\eta_1,K\eta_2)
\end{equation*}
and more generally $g(u\eta_1,u\eta_2)= g(\eta_1,\eta_2)$ for $u\in S^2$.
Thus, for each $u\in S^2$, the Riemannian metric $g$ is Hermitian with
respect to $u$, that is, of type $(1,1)$ with respect to $u$.

Each purely imaginary unit quaternion $u\in S^2$
determines a real-symplectic structure, or {\em \ka~ form\/} $\omega_u$:
\begin{equation}\label{eq:kahlerform}
\omega_u(\eta_1,\eta_2) :=  g(u\eta_1,\eta_2). 
\end{equation}
$\omega_u$ is the real-symplectic structure associated to 
a \ka~ structure $H$
for which the associated Riemannian structure is $g$, with respect to
the complex structure $u$:
\begin{align*}
\Im H & = \omega_u \\ \Re H & = g.
\end{align*}
The general definition of a \hk structure is a quadruple
\begin{equation*}
(g,I,J,K)  
\end{equation*}
where $g$ is a Riemannian metric and $I,J,K$ are
integrable almost complex structures satisfying the quaternion
identities \eqref{eq:quaternionidentities}, for which $g$ is
\kan~ with respect to each of them. 

A \hk structure also determines complex-symplectic structures.
For each fixed $u\in S^2$, there is a family of complex-symplectic structures
which are holomorphic with respect to $u$. Namely, extend $u$ to
a positively oriented orthonormal frame 
$(u,u_1,u_2)$, and define:
 \begin{equation}\label{eq:csympst}
\Omega_{(u_1,u_2)} :=   \omega_{u_1} + i \omega_{u_2},
\end{equation}
where the \ka~forms $\omega_{u_i}$ are defined in
\eqref{eq:kahlerform}.
Then $\Omega_{(u_1,u_2)}$ is $\C$-bilinear 
with respect to the complex structure $u$. 
Furthermore every other
positively oriented orthonormal frame 
$(u,u_1',u_2')$ is obtained by rotation:
\begin{align*}
u_1' & = \cos(\theta)u_1 - \sin(\theta)u_2 \\
u_2' & = \sin(\theta) u_1 + \cos(\theta) u_2
\end{align*}
for some $\theta\in\R$ and
\begin{equation*}
\Omega_{(u_1',u_2')} =   e^{i\theta} \Omega_{(u_1,u_2)}. 
\end{equation*}
The unit vector $u_2$ is determined by $u$ and $u_1$ as the cross-product
$u\times u_1$. Thus these complex-symplectic structures are parametrized
by orthonormal pairs $(u,u_1)$. 
Since the tangent space to $S^2$ at $u\in S^2$ is the orthogonal complement
$u^\perp\subset\R^3$, the complex-symplectic structures are parametrized
by the {\em unit tangent bundle\/} $T_1S^2$.
Left-multiplication by $u$ preserves $u^\perp$
and defines an integrable almost complex structure on $S^2$.
In \S\ref{sec:projective line}, we describe 
how this complex manifold identifies with the complex
projective line $\P^1$ and the
corresponding line bundle is the {\em holomorphic tangent bundle\/}
of $P^1$. It has degree $2$ and is denoted $\Os(2)$.
In \S\ref{sec:twistor} these complex-symplectic structures
on the moduli space unify into one single holomorphic exterior
$2$-form on the twistor space which takes values in a line bundle
induced from $\Os(2)$.

\subsection{Complex-symplectic structure}
An alternative approach to \hk geometry is to begin with a complex-symplectic
manifold and refine the complex-symplectic structure to a \hk structure.
Suppose that $(M,J,\Omega)$ is a complex-symplectic manifold, that is,
$(M,J)$ is a complex manifold and $\Omega$ is closed nondegenerate
complex-valued exterior 2-form of type $(2,0)$ with respect to $J$.
A {\em \hk structure refining $(J,\Omega)$\/}
consists of another complex structure 
$I$ which is anti-linear with respect to $J$ and for which the symmetric
form defined by 
\begin{equation}\label{eq:kametric}
g(\alpha,\beta) :=  \Re \Omega(\alpha,I\beta) 
\end{equation}
is \kan~ with respect to $J$ and $I$.
 
This approach can be succinctly described in terms of classical Lie groups
as follows. Let $\C^{2n}$ be the standard $2n$-dimensional
complex symplectic vector space.
The complex symplectic group $\Sp(2n,\C)$ is a subgroup of the 
automorphisms of $\C^{2n}$.
The stabilizer of the standard Hermitian structure on $\C^{2n}$ 
is the compact symplectic group
\begin{equation*}
\Sp(2n) = \Sp(2n,\C) \cap \o{U}(2n),
\end{equation*}
the maximal compact subgroup of $\Sp(2n,\C)$.

A reduction of structure group for the tangent bundle of a $4n$-dimensional
real manifold from $\GL(4n,\R)$ to $\Sp(2n,\C)$ consists of a pair
$(\Omega,J)$ where $J$ is an almost complex structure  and 
$\Omega$ a $J$-bilinear complex-valued nondegenerate exterior $2$-form
on $M$. 
This reduction of structure group is 
{\em integrable,\/} that is, corresponds to a complex-symplectic structure
on $M$ if and only if $J$ is integrable and $\Omega$ is closed.

A further reduction of structure group from $\Sp(2n,\C)$ to
its maximal compact subgroup $\Sp(2n)$ consists of another complex structure
$I$ which is $J$-anti-linear, and such that
\eqref{eq:kametric} defines a positive definite Hermitian metric
with respect to $I$. Such a reduction of structure group
gives tensor fields $g,I,J,K$ etc.\ which satisfy the \hk identities.
For further examples and discussion of \hk geometry, we refer the
reader to
Besse~\cite{Besse},
Dancer~\cite{Dancer},
Hitchin~\cite{Bourbaki}, 
Hitchin-Karlhede-Lindstr\"om-Ro\u{c}ek~\cite{HKLR},
Joyce~\cite{Joyce}, or
Salamon~\cite{Salamon}.

\subsubsection*{Relation to the moduli spaces}
These objects relate to the $J$-holo-morphic
complex-symplectic structure $\Omega$ 
on the de Rham moduli space in \S\ref{sec:dRgpd}.
By the definitions (\eqref{eq:OmegaJ}, \eqref{eq:I},
\eqref{eq:g}) of $\Omega, I, g$,
\begin{equation*}
g(\eta_1,\eta_2) = \Re\, \Omega(I \eta_1,\eta_2).
\end{equation*}
Thus
\begin{equation}\label{eq:OmegaIK}
\Omega = \Omega_{(-I,K)} = -\omega_I + i \omega_K 
\end{equation}
in the notation of \eqref{eq:csympst}, 
with real and imaginary parts discussed in \eqref{eq:OmegaRI}.

Similarly the $I$-holomorphic 
complex-symplectic structure $\Omega_{T^*}$ on the
Dolbeault moduli space \eqref{eq:OmegaI} is:
\begin{equation}\label{eq:OmegaTstar}
\Omega_{T^*} =  8 \Omega_{J,K}
\end{equation}
in the notation of \eqref{eq:csympst}.
We indicate the proof of \eqref{eq:OmegaTstar} as follows:
The projections onto Hodge components are:
\begin{xalignat*}{2}
\hH^1(X) & \longrightarrow \hH^{0,1}(X) 
&\qquad 
\hH^1(X) & \longrightarrow \hH^{0,1}(X)  
\\
\eta &\longmapsto  \frac12 (\eta - i\star\eta)
&\qquad 
\eta &\longmapsto  \frac12 (\eta + i\star\eta).
\end{xalignat*}
The $(0,1)-$ and $(1,0)-$ components of $\psi = 2i \Im(\eta)$
and $\phi = 2 \Re(\eta)$ respectively, are:
\begin{align*}
\Psi & = \; \frac{1-i\star}2 \;\; \frac{\eta-\bar\eta}2 \; =\;
\frac14 (\eta -\bar\eta - i\star\eta + i\star\bar\eta) \\
\Phi & = \; \frac{1+i\star}2 \;\; \frac{\eta+\bar\eta}2  \; =\;
\frac14 (\eta +\bar\eta + i\star\eta + i\star\bar\eta).
\end{align*}
Applying \eqref{eq:etaPsiPhi} to \eqref{eq:OmegaI}
\begin{multline*}
\Omega_{T^*}(\eta_1,\eta_2)  =  2 \bigg\{
(\eta_1\wedge\eta_2 -\bar\eta_1\wedge\bar\eta_2) \\
 + (\eta_1\wedge i\star\bar\eta_2 -\bar\eta_1\wedge i\star\eta_2)  
+(-i\star\eta_1\wedge\bar\eta_2 + i\star\bar\eta_1\wedge\eta_2)  \\
 + (\star\eta_1\wedge\star\eta_2 - \star\bar\eta_1\wedge\star\bar\eta_2) 
\bigg\} 
\end{multline*}
Taking real parts and imaginary parts, \eqref{eq:OmegaTstar} follows.
The \ka~form on $T^*\Jac(X)$ is $\omega_I$.

\subsection{Quaternionization}

The \hk structure identifies the vector space
$V\oplus\bV$ with the tensor product (over $\C$)
\begin{equation*}
V_{\bH} := \bH \otimes_{\C} V.
\end{equation*}
Here the quaternion algebra $\bH$ is a complex vector space 
where the action of $i$ is left-multiplication
by $I\in\bH$.
This quaternionic action extends to $V_\bH$, making $V_\bH$ a
left $\bH$-module.

Identify $V\oplus\bV$ with $V_\bH$  as follows:
Let $a:V\longrightarrow\bV$ be the anti-isomorphism given
by the identity map on $V_\R$. The $\R$-linear automorphism 
\begin{equation*}
(\Psi,\Phi) \stackrel{J}\longmapsto (-a^{-1}(\Phi),a(\Psi))
\end{equation*}
of $V\oplus\bV$ is anti-linear and satisfies $J^2 = -1$.
The given complex structures $I$ and $J$ define a 
left $\bH$-action on $V\oplus\bV$. The resulting map
\begin{equation*}
V_\bH \longrightarrow V\oplus\bV 
\end{equation*}
is an isomorphism of left $\bH$-modules.

The two complex structures $\pm I$ are special, since the $\R$-span of
$L$ is invariant under each of these structures. For those two complex
structures, $V_\bH/L$ contains the complex torus $V/L$
as a compact holomorphic submanifold. However, 
for all other complex structures $u\in S^2$, 
the resulting complex manifold is {\em Stein:\/}

\begin{thm} Let $u\in S^2$. Suppose that $u\neq \pm I$.
Then the complex manifold
$(V_\bH/L,\J_u)$ is biholomorphic to $(\C^*)^{2k}$ where
$k = \dim_\C(V)$.
\end{thm}
\begin{proof}
The lattice $L$ spans $V$ over $\R$, and $V$ is $I$-invariant.
Denote by $\C_u$ the subalgebra $\R 1 + \R u$ generated by $u$. 
Since $u\neq \pm I$, the two unit quaternions $I$ and $u$ generate
$\bH$. Thus $V$ generates $V_\bH$ over the subalgebra $\C_u\subset\bH$.
Thus the $\R$-linear map
\begin{align*}
\C \otimes_{\R} V & \longrightarrow V_\bH \\
(x + y i)\otimes v & \longmapsto  (x + y u) v
\end{align*}
is surjective. It is an isomorphism because
\begin{equation*}
\dim_\R  \big(\C \otimes_{\R} V\big) = 4k = \dim_\R  V_\bH.
\end{equation*}
As abelian groups, $V = L \otimes \R$, so 
\begin{equation*}
V_\bH \cong \C \otimes_{\R} V \cong \C \otimes L
\end{equation*}
Since $L$ is free abelian of rank $2k$,
\begin{align*}
V_\bH/L & \cong (\C\otimes L)/L \\ 
& \cong  (\C\otimes L)/(\Z\otimes L) \\ 
& \cong (\C/\Z)\otimes L \\ 
&\cong (\C^*)\otimes L \cong (\C^*)^{2k}
\end{align*}
as desired
\end{proof}

\section{The twistor space}\label{sec:twistor}
A \hk structure is equivalent to a holomorphic construction which
embodies all the complex and complex-symplectic structures, 
due to Hitchin-Karlhede-Lindstr\"om-Ro\u{c}ek~\cite{HKLR}.

For each unit purely imaginary quaternion $u\in S^2$, 
left-multiplication by $u$ describes a complex
structure on $V_\bH/L$. 
The collection of complex manifolds
$(V_\bH/L,u)$ 
unifies in a single holomorphic fibration $\Pi_Z$ over $S^2$, where
$S^2$ identifies with the complex projective line $\P^1$ 
by stereographic projection \eqref{eq:stereo}.
The corresponding holomorphic fibration 
$Z\xrightarrow{\Pi_Z} \P^1$ is 
{\em the twistor space\/} of $\M$.

In this section we explicitly describe this fibration
and use it to describe Simpson's $\C^*$-action on $\M$.
We realize $\P^1$ as a space of complex structures on $\bH$
by stereographic projection. 
The twistor space results, as well as a closed holomorphic
exterior $2$-form $\bO$ on $Z$. This exterior $2$-form
defines the holomorphic  complex-symplectic structure on the fibers of $\Pi_Z$,
taking values in the pullback $\Pi_Z^*\Os(2)$ 
of the line bundle $\Os(2)$ on $\P^1$.
Differentiably, $Z$ is the product bundle 
$V_\bH/L \times \P^1$, but holomorphically it is nontrivial.
The rational curves $l_\vv = \{[\vv]\}\times\P^1 $  define a family of
holomorphic sections of $Z$, {\em the twistor lines,\/} whose normal bundles
are isomorphic to $\Os(1)$. Finally the antipodal map $\balpha$ on $\P^1$
lifts to an anti-holomorphic involution on $Z$ which satisfies compatibility
relations with the twistor lines and $\bO$. The basic result 
(Theorem~3.3 of \cite{HKLR})
is that the collection $(\Pi_Z,\{l_{\vv}\},\bO,\tbalpha)$
is equivalent to a \hk structure. 
We describe these objects explicitly in this special case
(Theorem~\ref{thm:hktwistor}).

Finally the $\C^*$-action is described in terms of the twistor  
construction.

\subsection{The complex projective line}\label{sec:projective line}
We briefly review the geometry of $\P^1$. For more details, see
\S 1--\S 2 of Hitchin's article in~\cite{HSW}.

Points in $\P^1$ correspond to complex lines 
($1$-dimensional linear subspaces)
in $\C^2$.  Cover $\P^1$ by two affine patches $\A_1, \A_2$ with coordinate
charts:
\begin{xalignat*}{2}
\C & 
\xrightarrow{\;\psi^{(1)}} \A_1  
&\qquad 
\C & \xrightarrow{\;\psi^{(2)}} \A_2 
\\
\zeta &\longmapsto  \C \bmatrix \zeta \\ 1 \endbmatrix \subset \C^2
&\qquad 
\xi &\longmapsto  \C \bmatrix 1 \\ \xi  \endbmatrix \subset \C^2
\end{xalignat*}
where $\zeta\xi = 1$ on the intersection $\C^* \cong \A_1\cap \A_2$.

\subsubsection*{Stereographic projection}
The unit sphere $S^2$ of purely imaginary unit quaternions, with its complex
structure defined in \S\ref{sec:quaternion}, identifies with $\P^1$ as follows.
We embed $\C\subset\bH$ as the subalgebra generated by $I$.
The two affine patches $\A_1, \A_2\subset \P^1$ map affinely into the
multiplicative subgroup $\bH^*$ of nonzero quaternions by:
\begin{align*}
\A_1 &\longrightarrow \bH^* \\
\psi^{(1)}(\zeta) & \longmapsto 1 + \zeta K
\end{align*}
and 
\begin{align*}
\A_2 &\longrightarrow \bH^* \\
\psi^{(2)}(\xi) & \longmapsto \bar{\xi} +  K
\end{align*}
respectively. Use these maps to conjugate $I\in\bH$ to purely imaginary
unit quaternions in $S^2$:
\begin{align}
I^{(1)}_\zeta  &:= (1 + \zeta K) I  (1 + \zeta K)^{-1} \notag\\
I^{(2)}_\xi  &:= (\bar{\xi} +  K) I  (\bar{\xi} +  K)^{-1} 
= (\xi - K)^{-1} I (\xi - K) 
\label{eq:stereo},
\end{align}
where $\zeta,\xi\in\C$. These two maps combine to define a map
$\P^1\longrightarrow S^2$ because:
\begin{lem}
If $\xi\zeta = 1$, then $I^{(1)}_\zeta = I^{(2)}_\xi.$
\end{lem}
\begin{proof}
Since $\zeta\in\C = \R 1 + \R I$, 
\begin{equation*}
1 + \zeta K = 1 +  K\bar{\zeta} = (\bar{\zeta}^{-1} + K)\bar{\zeta} =
(\bar{\xi} + K)\bar{\zeta}
\end{equation*}
so
\begin{equation*}
(1 + \zeta K)^{-1} = (\bar\zeta)^{-1}  (\bar\xi + K)^{-1}
\end{equation*}
and 
\begin{align*}
I^{(1)}_\zeta & = (1 + \zeta K) I(1 + \zeta K)^{-1} \\ & = 
(\bar\xi + K)\bar\zeta I (\bar\zeta)^{-1}  (\bar\xi + K)^{-1} \\
& = (\bar\xi + K) I  (\bar\xi + K)^{-1} \\
& = I^{(2)}_\xi.
\end{align*}
\end{proof}
\noindent
Evidently $I^{(1)}_0 = I$ , $I^{(2)}_0 = -I$, and
$I^{(1)}(1) = I^{(2)}(1) = J$.
Explicitly, 
$I^{(1)}_\zeta, I^{(2)}_\xi$ are {\em stereographic projections:\/}
\begin{align}
I^{(1)}_\zeta  & \,=\, 
\frac{1 - \zeta\bar\zeta}{1 + \zeta\bar\zeta} I \,+\,
\frac{2\Re(\zeta)}{1 + \zeta\bar\zeta} J \,+\,
\frac{2\Im(\zeta)}{1 + \zeta\bar\zeta} K \notag \\
I^{(2)}_\xi  & \,=\, 
- \frac{1 - \xi\bar\xi}{1 + \xi\bar\xi} I \,+\,
\frac{2\Re(\xi)}{1 + \xi\bar\xi} J \,-\,
\frac{2\Im(\xi)}{1 + \xi\bar\xi} K. \label{eq:explicitstereo}
\end{align}
(This differs slightly from formulas (3.70) and (3.71) of
\cite{HKLR}. The formulas in \cite{HKLR} relate
to ours by replacing $\zeta$ with $\bar\zeta$. In our convention,
stereographic projection is a map from $\P^1$ to
$S^2$ with the orientation of $S^2$ given by the {\em outward\/} normal
vector.)

\subsubsection*{The antipodal map}
Stereographic projection relates the antipodal map 
$S^2\longrightarrow S^2$ defined in \eqref{eq:antipodal} 
to the map $\P^1\stackrel{\balpha}\longrightarrow\P^1$ induced by
the anti-linear involution of $\C^2$:
\begin{align*}
\C^2 & \stackrel{\bar\balpha}\longrightarrow \C^2 \\
\bmatrix z_1 \\ z_2 \endbmatrix & \longmapsto
\bmatrix -\bar{z}_2 \\ \bar{z}_1 \endbmatrix. 
\end{align*}
In terms of the atlas, it is given by:
\begin{align*}\label{eq:alphaatlas}
\psi^{(1)}(\zeta) &\stackrel{\balpha}\longleftrightarrow 
\psi^{(2)}(-\bar\zeta) \notag \\
\psi^{(2)}(\xi) &\stackrel{\balpha}\longleftrightarrow 
\psi^{(1)}(-\bar\xi) 
\end{align*}
for $\xi,\zeta\in\C$ and in particular
\begin{equation*}
\psi^{(1)}(0) \stackrel{\balpha}\longleftrightarrow  \psi^{(2)}(0). 
\end{equation*}
\noindent

\begin{lem}\label{lem:antipodcs}
The antipodal map takes a complex structure to its opposite. That is,
if $\zeta\in\C$ is nonzero, then
\begin{equation*}
I^{(1)}_{-\bar\zeta^{-1}} =  I^{(2)}_{-\bar\zeta} = 
- I^{(1)}_{\zeta} 
\end{equation*}
\end{lem}
\begin{proof}
\begin{align*}
I^{(1)}_{-\bar\zeta^{-1}} & =  
(1 - \bar\zeta^{-1} K) I (1 - \bar\zeta^{-1} K)^{-1} \\
& = \frac1{1 + \vert\bar\zeta^{-1}\vert^2}\;
(1 - \bar\zeta^{-1} K) I (1 + \bar\zeta^{-1} K) \\
& = \frac{\vert\zeta\vert^2}{\vert\zeta\vert^2 + 1 }\;
(1 - \bar\zeta^{-1} K)  I (1 + \bar\zeta^{-1} K) \\
& = \frac1{1 + \vert\zeta\vert^2}\;
(1 - K \zeta^{-1}) \zeta I \bar\zeta (1 + \bar\zeta^{-1} K) \\
& = \frac1{1 + \vert\zeta\vert^2}\;
(\zeta  - K ) I (\zeta +  K) \\
& = \frac1{1 + \vert\zeta\vert^2}\;
(- \zeta K - 1) KIK (1 - \zeta K) \\
& = \; - \; \frac1{1 + \vert\zeta\vert^2}\;
(1 + \zeta K) I (1 - \zeta K) \\
& = \; - \; I^{(1)}_\zeta
\end{align*}
as claimed.
\end{proof}

\subsubsection*{Line bundles}
For each $d\in\Z$ there is a holomorphic line bundle
$\Os(d)$ of degree $d$ over $\P^1$, unique up to isomorphism.
Every holomorphic line bundle over $\P^1$ is isomorphic to some $\Os(d)$.

Explicitly, the total space of $\Os(d)$ is defined in the above atlas
as the identification space of the disjoint union
\begin{equation*}
\bigg(\C \times \A_1 \bigg) 
\coprod 
\bigg(\C \times \A_2 \bigg) 
\end{equation*}
under the equivalence relation:
\begin{equation}\label{eq:Od}
\big(z,\psi^{(1)}(\zeta), z\big)  \sim 
\big(\xi^d z,\psi^{(2)}(\xi)\big),
\end{equation}
where $\zeta\xi = 1$. A section $s$ of $\Os(d)$ determines a pair
$\big(s^{(1)},s^{(2)}\big)$ of maps $\C \longrightarrow \C$ such that
$s\big(\psi^{(1)}(\zeta)\big)$ corresponds to 
$\big(\psi^{(1)}(\zeta), s^{(1)}(\zeta)\big)$ in
$\A_1\times\C$
and $s\big(\psi^{(2)}(\xi)\big)$ corresponds to 
$\big(\psi^{(2)}(\xi), s^{(2)}(\xi)\big)$
in $\A_2\times\C$.
By \eqref{eq:Od}, the pair $\big(s^{(1)},s^{(2)}\big)$ must satisfy:
\begin{equation}\label{eq:sections}
s^{(1)}(\zeta) = \xi^{-d} s^{(2)}(\xi)
\end{equation}
whenever $\zeta\xi=1$.
Conversely any pair
$\big(s^{(1)},s^{(2)}\big)$
satisfying \eqref{eq:sections} corresponds to a
section of $\Os(d)$

The trivial rank two bundle 
\begin{equation*}
\C^2 \times \P^1 \longrightarrow \P^1 
\end{equation*}
contains the (tautological) subbundle $\Os(-1)$ whose fiber over
$[v]\in\P^1$ is the line $\C v\subset\C^2$.  
The {\em hyperplane bundle\/} over $\P^1$
is its dual $\Os(1)$: the 
fiber over $\C v$ 
consists of linear functionals on $\C_v$.  A linear functional
\begin{equation*}
\C^2 \stackrel{a}\longrightarrow \C 
\end{equation*}
defines a section of $\Os(1)$ by
assigning to $[v]\in\P^1$ its restriction to $\C v$.
Every holomorphic section of $\Os(1)$ 
arises from such a linear functional on $\C^2$. 
A nonzero holomorphic section has 
a single zero of order $1$, corresponding to the
kernel of $a$. 
For example, the linear functional
\begin{equation*}
\bmatrix z_1 \\ z_2 \endbmatrix \longmapsto z_2 
\end{equation*}
determines the section $\sigma$ of $\Os(1)$ having
\begin{equation*}
s^{(1)}(\zeta) = 1,\quad 
s^{(2)}(\xi) = \xi.
\end{equation*}
This section is holomorphic on $\A_1$ and has a single simple zero
at $\psi^{(2)}(0)\in\A_2$.

The {\em canonical line bundle\/} $T^*\P^1$ has degree $-2$ and is
isomorphic to $\Os(-2)$;
the holomorphic 1-form 
\begin{equation*}
d\zeta = -\xi^{-2} d\xi 
\end{equation*}
is holomorphic on $\A_1$ with a double pole at $\zeta = \infty$. 
In the atlas, this section is given by holomorphic functions
\begin{align*}
s^{(1)}(\zeta) & = 1 \\
s^{(2)}(\xi)  & = \xi^{-2}
\end{align*}
where the local nonvanishing sections over $\A_1$ and $\A_2$ are
$d\zeta$ and $-d\xi$ respectively.

The inverse of the canonical bundle, 
the {\em tangent bundle\/} $T\P^1\cong \Os(2)$, 
has a holomorphic section 
\begin{equation}\label{eq:o2}
\delta = \frac{\partial}{\partial\zeta} = -\xi^2 \frac{\partial}{\partial\xi} 
\end{equation}
with a double zero at $\xi = 0$. 
In coordinate charts,
\begin{align*}
\delta^{(1)}(\zeta) & =  1 \\ 
\delta^{(2)}(\xi) & = \xi^2.
\end{align*}
This section corresponds to the section
$-\sigma^2 = -\sigma\otimes\sigma$ of $\Os(2)=\Os(1)\otimes\Os(1)$.
 
\subsubsection*{Antiholomorphic line bundles}
Given any holomorphic vector bundle $E\longrightarrow M$, the {\em
conjugate\/} vector bundle $\bar{E}\longrightarrow M$ is defined as
the holomorphic vector bundle with opposite complex structure on the
fibers. Fiberwise complex-conjugation defines an isomorphism of smooth
vector bundles between $E$ and $\bar{E}$ which is {\em antiholomorphic.\/}
A Hermitian structure on $\C^2$ defines Hermitian structures on any $\Os(d)$
which induces isomorphisms
\begin{equation*}
\overline{\Os(d)} \cong \Os(-d).
\end{equation*}
\begin{lem}\label{lem:tbalpha}
Suppose that $d\in\Z$ is even. Then the antipodal map
$\balpha$ of $\P^1$ lifts to an antiholomorphic isomorphism
\begin{equation*}
\Os(d)\xrightarrow{\tbalpha} \overline{\Os(d)}.
\end{equation*}
\end{lem}
\begin{proof}
In the coordinate atlas, define
\begin{align*}
\big(z, \psi^{(1)}(\zeta)\big) & \stackrel{\tbalpha}\longmapsto
\big(\bar{z}, \psi^{(2)}(-\bar\zeta)\big) \\
\big(w, \psi^{(2)}(\xi)\big) & \stackrel{\tbalpha}\longmapsto
\big(\bar{w}, \psi^{(1)}(-\bar\xi)\big). 
\end{align*}
Suppose that 
$z,\zeta\in\C$ define a point  $\big(z,\psi^{(1)}(\zeta)\big)$ 
in $\C\times \A_1$, 
whose image 
$\tbalpha\big(z,\psi^{(1)}(\zeta)\big)$ corresponds to
\begin{equation}\label{eq:first}
\big(\ (-\bar\xi)^d\bar{z},\  \psi^{(1)}(-\bar\xi)\ \big)\,\in\, \C \times\A_1,
\end{equation}
where $\zeta\xi = 1$. On the other hand, 
$\big(z,\psi^{(1)}(\zeta)\big)$ corresponds to the point
$\big(\xi^dz,\ \psi^{(2)}(\xi)\big)$ in $\C\times \A_2$.
This point maps under $\tbalpha$ to 
$\big((\overline{\xi^dz}),\  \psi^{(1)}(-\bar\xi)\big)$,
which equals \eqref{eq:first} when $d$ is even.
\end{proof}
\noindent When $d=2$, the image $(\tbalpha)^*d\zeta =
-\overline{d\xi}$ is
the image of the section defined above under the differential of $\balpha$.

\subsection{The twistor space as a smooth vector bundle}
The triviality of the vector bundle $\C^{2k}\otimes\Os(1)$ as a 
{\em smooth real vector bundle\/} 
may be seen as follows. Isomorphism classes of
bundles over $\P^1$ with structure group $G$ are classified by the
homotopy class of their {\em clutching function\/} 
$\A_1\cap\A_2\longrightarrow G$, which identifies with an element of
$\pi_1(G)$
(see for example Steenrod~\cite{Steenrod}).
Namely, a fiber bundle over a contractible space is trivial,
and writing a trivialization over $\A_1$ in terms of a trivialization
over $\A_2$ defines a map $\A_1\cap\A_2\longrightarrow G$, whose homotopy
class is a complete invariant of the isomorphism type as a smooth
fiber bundle with structure group $G$. Since $\A_1\cap\A_2$ is 
homotopy-equivalent to $S^1$, isomorphism classes of $G$-bundles
are classified by $\pi_1(G)$.

The invariant of the line bundle $\Os(1)$ is a generator $\alpha\in\pi_1(\Uo)$.
The invariant of $\C^{2k}\otimes\Os(1)$ is the image of 
\begin{equation*}
\alpha\oplus\dots\oplus\alpha\,\in\,\pi_1\big(\Uo\times\dots\times\Uo\big) 
\end{equation*}
under the homomorphism 
\begin{equation*}
\begin{CD}
\pi_1\big(\Uo\times\dots\times\Uo\big) @>>> \pi_1\big(\o{U}(2k)\big) \\
@| @| \\
\Z^{2k} @>>> \Z
\end{CD}
\end{equation*}
induced by the inclusion of diagonal matrices 
$\Uo\times\dots\times\Uo$ in $\o{U}(2k)$. The isomorphism 
$\pi_1(\o{U}(2k))\longrightarrow \Z$ is induced by 
\begin{equation*}
\det:\o{U}(2k)\longrightarrow \Uo.  
\end{equation*}
Thus the invariant of
$\C^{2k}\otimes\Os(1)$  as a complex vector bundle equals
$2k\alpha$, where $\alpha$ denotes a generator of 
$\pi_1\big(\o{U}(2k)\big)\cong\Z$.
The isomorphism type {\em as a real vector bundle\/} is the image of
$2k \alpha$ under the homomorphism
\begin{equation*}
\begin{CD}
\pi_1(\o{U}(2k))  @>>> \pi_1(\o{O}(4k)) \\
@|    @| \\
\Z @>>> \Z /2
\end{CD}
\end{equation*}
which is evidently zero.

\subsection{A holomorphic atlas for the twistor space}

The twistor space of the quaternionic vector space
$V_\bH$ % and later pass to the quotient $V_\bH/L$. 
% The twistor space for $V_\bH$ is the universal covering space
% $\tZ$ of the twistor space $Z$ for $V_\bH/L$. 
is defined as the product $\tZ = V_\bH\times \P^1$ with coordinate 
projection:
\begin{equation*}
\tZ =  V_\bH\times \P^1 \stackrel{\Pi_{\tZ}}\longrightarrow \P^1.
\end{equation*}
The complex structure on $\tZ$ is the unique complex structure 
such that $\Pi_{\tZ}$ is a holomorphic
vector bundle, and the complex structure on the fiber over $[v]\in\P^1$
is the complex structure defined by left-multiplication by the purely
imaginary unit quaternion corresponding to $[v]$. 
We shall see that the twistor space $\tZ$ of $V_\bH$ is the total space
of the holomorphic vector bundle $\C^{2k}\otimes\Os(1)$, and describe
the structures (holomorphic symplectic structures, real structure) 
on this space explicitly in terms of the atlas $\{(\A_1,\psi^{(1)}),
(\A_2,\psi^{(2)})\}$ in this and subsequent sections.

% The complex structure on $\bH$ defined by 
% left-multiplication by a purely imaginary unit quaternion $u\in S^2$ 
% induces a complex structure on $V_\bH$ and hence on 
% $V_\bH/L$.

Represent complex vectors $\vv\in\C^{2k}$ by the notation:
\begin{equation}\label{eq:cvectornotation}
\vv = \bmatrix \qq \\ \pp \endbmatrix   
\end{equation}
where the components are vectors  $\qq,\pp\in\C^k$. 
Thus a general quaternion vector in $\bH^k$ is:
\begin{equation*}
\qq + \pp J  = \xx_1 + \yy_1 I + \xx_2 J + \yy_2 K,
\end{equation*}
where
$\qq = \xx_1 + i \yy_1$ and $\pp = \xx_2 + i \yy_2$
for real vectors $\xx_1,\xx_2,\yy_1,\yy_2\in\R^k$.

For $\zeta,\xi\in\C, \vv \in \C^{2k}$, define
\begin{align*}
\C^{2k} & \xrightarrow{\;f^{(1)}_\zeta} 
  \bH^k \\
\vv = \bmatrix \qq \\ \pp \endbmatrix 
& \longmapsto 
\frac1{1 + \zeta \bar\zeta} 
(1 + \zeta K) (\qq + \pp J)
= (1 - \zeta K)^{-1} (\qq + \pp J)
\end{align*}
and
\begin{align*}
\C^{2k} & \xrightarrow{\;f^{(2)}_\xi}   \bH^k \\
\vv = \bmatrix \qq \\ \pp \endbmatrix 
& \longmapsto 
\frac1{1 + \xi \bar\xi} 
(\bar\xi + K) (\qq + \pp J)
= (\xi - K)^{-1} (\qq + \pp J).
\end{align*}
For each $\xi,\zeta\in\C$, the maps $f^{(1)}_\zeta,f^{(2)}_\xi$,
are $\R$-linear isomorphisms.

%Define a holomorphic rank $2k$ vector bundle $\tZ$ over $\P^1$ by charts: 
%\begin{align*}
%\C^{2k}\times\A_1& \xrightarrow{~F^{(1)}~} \bH^k\times\A_1  \\
%\C^{2k}\times\A_2& \xrightarrow{~F^{(2)}~} \bH^k\times\A_2
%\end{align*}
%\begin{align}\label{eq:holocharts}
%F^{(1)}(\vv, \zeta) & = \big(f^{(1)}_\zeta(\vv),\psi^{(1)}(\zeta)\big) \notag\\
%F^{(2)}(\vv, \xi)   & = \big(f^{(2)}_\xi(\vv),\psi^{(2)}(\xi)\big) 
%\end{align}

%Then we have:
%\begin{equation*}
%F^{(1)}\big(\vv,\psi^{(1)}(\zeta)\big) = 
%F^{(2)}\big(\xi\vv,\psi^{(2)}(\xi)\big).
%\end{equation*}

Hence if we define charts
\begin{align*}
\C^{2k}\times\A_1& \xrightarrow{~F^{(1)}~} \bH^k\times\A_1  \\
\C^{2k}\times\A_2& \xrightarrow{~F^{(2)}~} \bH^k\times\A_2
\end{align*}
\begin{align}\label{eq:holocharts}
F^{(1)}(\vv, \zeta) & = \big(f^{(1)}_\zeta(\vv),\psi^{(1)}(\zeta)\big) \notag\\
F^{(2)}(\vv, \xi)   & = \big(f^{(2)}_\xi(\vv),\psi^{(2)}(\xi)\big), 
\end{align}
we have
\begin{equation*}
F^{(1)}\big(\vv,\psi^{(1)}(\zeta)\big) = 
F^{(2)}\big(\xi\vv,\psi^{(2)}(\xi)\big),
\end{equation*}
thus, 
obtaining a holomorphic rank $2k$ vector bundle 
\begin{equation*}
\Pi_{\tZ} : \tZ \longrightarrow \P^1.
\end{equation*}
The fiber of $\tZ_\zeta$ over $\psi^{(1)}(\zeta)$ 
is $\bH^k$ with the complex structure $I^{(1)}_\zeta$.
We denote the fiber $\Pi_{\tZ}^{-1}\big(\psi^{(2)}(0)\big)$
by $\tZ_\infty$. 

%define holomorphic charts on $\bH^k \times \A_1$ 
%and $\bH^k \times \A_2$.
% \pagebreak
\begin{lem}\label{lem:Zatlas}
Let $\xi\zeta = 1$. Then $f^{(1)}_\zeta(\vv) = f^{(2)}_\xi(\xi \vv)$.
\end{lem}

\begin{proof}
\begin{align*}
f^{(1)}_\zeta(\vv) 
&= \frac{1}{1 + \zeta\bar\zeta}      (1 + \zeta K)      (\qq + \pp J) \\
&= \frac{\xi\bar\xi}{\xi\bar\xi + 1} (1 + \xi^{-1} K)   (\qq + \pp J) \\
&= 
\frac1{1 + \xi\bar\xi}
(1 + \xi^{-1} K) \xi\bar\xi   (\qq + \pp J) \\
&= 
\frac1{\xi\bar\xi + 1}
(1 + K \bar{\xi}^{-1}) \bar\xi\, \xi (\qq + \pp J) \\
&= 
\frac1{\xi\bar\xi + 1}
(\bar\xi + K)                    \xi (\qq + \pp J) \\ 
&=  f^{(2)}_\xi(\xi \vv)
\end{align*}
\end{proof}
\noindent
Thus $\{F^{(1)},F^{(2)}\}$ is an atlas for a holomorphic vector bundle
$\tZ$ of rank $2k$ over $\P^1$. Lemma~\ref{lem:Zatlas} implies that
this vector bundle is $\C^{2k}\otimes \Os(1)$, a direct sum of $2k$ copies
of $\Os(1)$.

Finally, the complex structure on the fiber 
$\tZ_\zeta$ is $I^{(1)}_\zeta$ and  $I^{(2)}_0$ if $\zeta=\infty$.
\begin{lem}
Let $\zeta,\xi\in\C, \zeta\xi = 1$. Then
\begin{equation*}
f^{(1)}_\zeta (i \vv) =  
I^{(1)}_\zeta f^{(1)}_\zeta(\vv), \qquad
f^{(2)}_\xi (i \vv) =  
I^{(2)}_\xi f^{(2)}_\xi(\vv).
\end{equation*}
\end{lem}
\begin{proof}
\begin{align*}
f^{(1)}_\zeta (i \vv) 
& =  (1 + \zeta\bar\zeta)^{-1} \big(1 + \zeta K\big) (i\qq + i\pp J) \\
& =  (1 + \zeta\bar\zeta)^{-1} \big(1 + \zeta K\big) I (\qq + \pp J) \\
& =  (1 + \zeta\bar\zeta)^{-1} I^{(1)}_\zeta \big(1 + \zeta K\big) 
(\qq + \pp J) \\
& = I^{(1)}_\zeta f^{(1)}_\zeta (\vv)
\end{align*}
and similarly for $f^{(2)}_\xi$.
\end{proof}

\subsection{The twistor lines} By
(\ref{eq:holocharts}) and Lemma \ref{lem:Zatlas}, 
each $\vv \in \C^{2k}$ defines a nowhere 
vanishing holomorphic section of $l_\vv$, 
called the {\em twistor line\/} corresponding
to $\vv$. The normal bundle to $l_\vv$ is isomorphic to the vector bundle 
$\Pi_{\tZ}$ itself, that is, to $\C^{2k}\otimes\Os(1)$.

Here is how the the twistor lines appear in the holomorphic
atlas. Write
\begin{equation*}
\vv_0 = \qq_0 + \pp_0 J 
\end{equation*}
with $\qq_0,\pp_0\in\C^k$. The twistor line $l_\vv\big(\psi^{(1)}(\zeta)\big)$
intersects the fiber $\tZ_\zeta$ at
\begin{equation*}
F^{(1)}\bigg(\bmatrix \qq^{(1)}(\zeta) \\ \pp^{(1)}(\zeta)\endbmatrix,
\psi^{(1)}(\zeta)\bigg),  
\end{equation*}
where $\qq^{(1)}(\zeta),\pp^{(1)}(\zeta)\in\C^{k}$ are defined by
\begin{equation}\label{eq:twistorlinepq}
\bmatrix \qq^{(1)}(\zeta) \\ \pp^{(1)}(\zeta) \endbmatrix := 
\bmatrix \qq_0 + i \zeta \bar{\pp}_0 \\ 
\pp_0 - i \zeta\bar{\qq}_0 \endbmatrix. 
\end{equation}
Similarly, the twistor line $l_\vv\big(\psi^{(2)}(\xi)\big)$
intersects $\tZ_{\xi^{-1}}$ at
\begin{equation*}
F^{(2)}\bigg(\bmatrix \qq^{(2)}(\xi) \\ \pp^{(2)}(\xi)\endbmatrix,
\psi^{(2)}(\xi)\bigg)  
\end{equation*}
where 
\begin{equation}\label{eq:twistorlinepqxi}
\bmatrix \qq^{(2)}(\xi) \\ \pp^{(2)}(\xi) \endbmatrix := 
\bmatrix \xi \qq_0 + i \bar{\pp}_0 \\ 
\xi \pp_0 - i \bar{\qq}_0 \endbmatrix.
\end{equation}
Conversely, the twistor line containing 
\begin{equation*}
F^{(1)}\bigg(\bmatrix \qq \\ \pp\endbmatrix,
\psi^{(1)}(\zeta)
\bigg) 
\in \Pi_{\tZ}^{-1}(\A_1)
\end{equation*}
intersects $\tZ_0$ at
\begin{equation*}
F^{(1)}\bigg(
\frac1{1 + \zeta\bar\zeta} 
\bmatrix \pp + i\zeta\bar\qq \\ \qq - i \zeta\bar\pp\endbmatrix,
\psi^{(1)}(0)\bigg)  
\end{equation*}
and the  twistor line containing 
\begin{equation*}
F^{(2)}\bigg(\bmatrix \qq \\ \pp\endbmatrix,\psi^{(2)}(\xi)\bigg) 
\in \Pi_{\tZ}^{-1}(\A_2)
\end{equation*}
intersects $\tZ_\infty$ at
\begin{equation*}
F^{(2)}\bigg(
\frac1{1 + \xi\bar\xi} 
\bmatrix \bar\xi\pp + i\bar\qq \\ \bar\xi\qq - i \bar\pp\endbmatrix,
\psi^{(2)}(0)\bigg).  
\end{equation*}

\subsection{The real structure on the twistor space}
The antipodal map $\P^1\stackrel{\balpha}\longrightarrow\P^1$ lifts to the
twistor space $\tZ$ preserving the fibration and the twistor lines.
In terms of the product fibration in the coordinate charts, the lift
$\tbalpha$ to $\tZ$ is given by:

\begin{equation*}
F^{(1)}\bigg( 
\bmatrix \qq \\ \pp \endbmatrix, 
\psi^{(1)}(\zeta)\bigg) 
\quad\stackrel{\tbalpha}\longleftrightarrow\quad
F^{(2)}\bigg( 
\bmatrix -i\bar\xi\bar{\pp} \\ i\bar\xi\bar{\qq} \endbmatrix, 
\psi^{(2)}(-\bar\xi)\bigg) 
\end{equation*}
where $\xi\zeta = 1$ and
\begin{equation}\label{eq:antipodpoles}
F^{(1)}\bigg( 
\bmatrix \qq \\ \pp \endbmatrix, \psi^{(1)}(0)\bigg) 
\quad\stackrel{\tbalpha}\longleftrightarrow\quad
F^{(2)}\bigg( 
\bmatrix -i\bar{\pp} \\ i\bar{\qq} \endbmatrix, 
\psi^{(2)}(0) \bigg). 
\end{equation}
Lemma~\ref{lem:antipodcs} and the antiholomorphicity of $\balpha$ imply
that $\tbalpha$ is an antiholomorphic involution of $Z$.

\subsection{Symplectic geometry of the twistor space}
The complex-symplectic structures defined in \eqref{eq:csympst}  
fit together 
giving a holomorphic exterior $2$-form on the complex manifold 
$\bO$ on $Z$, taking values in the
holomorphic line bundle $\Pi_{\tZ}^*\Os(2)$, where $\Pi_{\tZ}$ denotes
the twistor fibration.

Let $u\in S^2$. There is a family of complex-symplectic structures on 
$V_\bH$ which are holomorphic with respect to $u$, which are parametrized
by a unit vector $u_1\in T_1S^2$. 
Namely \eqref{eq:csympst} defines a complex-symplectic
structure $\Omega_{(u_1,u_2)}$, where $u_1\in u^\perp$ and 
$u_2 = u\times u_1$.
Let $u$ be $I^{(1)}_\zeta$ or $I^{(2)}_\xi$. Define $u_1$ 
(respectively $u_2$) as $J^{(1)}_\zeta, J^{(2)}_\xi$ 
(respectively $K^{(1)}_\zeta, K^{(2)}_\xi$), where:
\begin{align}
J^{(1)}_\zeta  &:= (1 + \zeta K) J  (1 + \zeta K)^{-1} \notag \\
K^{(1)}_\zeta  &:= (1 + \zeta K) K  (1 + \zeta K)^{-1} \notag \\
J^{(2)}_\xi  &:=  (\bar{\xi} +  K) J  (\bar{\xi} +  K)^{-1} \notag \\
K^{(2)}_\xi  &:=  (\bar{\xi} +  K) K  (\bar{\xi} +  K)^{-1}, \label{eq:defJK}
\end{align}
analogous to the definition \eqref{eq:stereo} of
$I^{(1)}_\zeta$ and $I^{(2)}_\xi$.

Evidently, $(J^{(1)}_0,K^{(1)}_0) = (J,K)$ and 
$(J^{(2)}_0,K^{(2)}_0) = (-J,K)$.
For $\zeta,\xi\in\C$, let
\begin{align*}
\bO^{(1)}_\zeta & := (1 + \zeta\bar\zeta) 
\Omega_{(J^{(1)}_\zeta,K^{(1)}_\zeta)} \\
\bO^{(2)}_\xi & := (1 + \xi\bar\xi) 
\Omega_{(J^{(2)}_\xi,K^{(2)}_\xi)}.
\end{align*}
Corresponding to $\psi^{(1)}(0)\in\A_1$ and
$\psi^{(2)}(0)\in\A_2$ are
\begin{equation*}
\bO^{(1)}_0 = \Omega_{(J,K)},\quad \bO^{(2)}_0 = \Omega_{(-J,K)}, 
\end{equation*}
respectively.
\begin{lem}
Let $\alpha$ and $\beta$ be tangent vectors and $\zeta,\xi\in\C$.  Then
\begin{align*}
\bO^{(1)}_\zeta (\alpha,\beta) & = \;  
\Omega_{(J,K)}\big( (1-\zeta K)\alpha,(1-\zeta K)\beta\big) \\
\bO^{(2)}_\xi (\alpha,\beta) & =  \;
\Omega_{(J,K)}\big( (\xi - K)\alpha,(\xi - K)\beta\big)
\end{align*}
\end{lem}
\begin{proof}
First recall the metric $g$ as defined in (\ref{eq:g}).  
For any $\alpha$ and $\beta$,
\begin{align*}
g\big( (1 + \zeta K)\alpha,\beta\big) & = 
g\big( \alpha,(1 - \zeta K)\beta\big) \\
g\big( (\bar\xi + K)\alpha,\beta\big) & = 
g\big( \alpha,(\xi - K)\beta\big)
\end{align*}
since $g(h\alpha,\beta)= g(\alpha,-h\beta)$ 
for any purely imaginary quaternion $h\in\bH$.
\begin{align*}
\bO^{(1)}_\zeta (\alpha,\beta)  & = (1 + \zeta\bar\zeta)\; 
\bigg(g(J^{(1)}_\zeta\alpha,\beta) + ig(K^{(1)}_\zeta\alpha,\beta)\bigg) \\
& = (1 + \zeta\bar\zeta)\; \bigg(
g\big( (1 + \zeta K) J (1 + \zeta K)^{-1} \alpha,\beta\big)  \\
& \qquad \qquad \qquad 
+  ig\big( (1 + \zeta K) K (1 + \zeta K)^{-1} \alpha,\beta\big)
\bigg) \\
& = 
g\big( (1 + \zeta K) J (1 - \zeta K) \alpha,\beta\big) 
+  ig\big( (1 + \zeta K) K (1 - \zeta K) \alpha,\beta\big)
\\
& = g\big( J (1 - \zeta K) \alpha, (1-\zeta K)\beta\big)
+ ig\big( K (1 - \zeta K) \alpha,  (1 - \zeta K)\beta\big) \\
& =  \Omega_{(J,K)}\big(( 1 - \zeta K) \alpha, (1-\zeta K)\beta\big) 
\end{align*}
and similarly for $\bO^{(2)}_\xi (\alpha,\beta).$
\end{proof}

\begin{cor} If $\zeta\xi = 1$, then
$ \bO^{(2)}_\xi  = \xi^2  \bO^{(1)}_\zeta $.
\end{cor}
\begin{proof} 
Since $\xi - K = \xi(1 - \zeta K)$ and $\Omega_{(J,K)}$ is 
holomorphic with respect to $I$,
\begin{align*}
\bO^{(2)}_\xi(\alpha,\beta) & = 
\Omega_{(J,K)} \big((\xi - K) \alpha,(\xi - K)\beta\big)  \\
& =  \xi^2   \Omega_{(J,K)} 
\big((1 - \zeta K) \alpha, (1 - \zeta K)\beta\big) \\
& = \xi^2  \bO^{(1)}_\zeta(\alpha,\beta)
\end{align*}
as desired.
\end{proof}
\noindent
Thus $\{\bO^{(1)},\bO^{(2)}\}$ defines a $\Pi_{\tZ}^*\Os(2)$-valued 
exterior 2-form $\bO$ on $\tZ$.  Denoting $\delta$ the meromorphic
section of $\Os(-2)$ in \eqref{eq:o2}, 
we see that $\bO\otimes\Pi_Z^*\delta$ is a (scalar-valued) exterior 2-form 
on $\tZ$.
\begin{cor} $\bO$ is holomorphic.
\end{cor}
\begin{proof}

$\bO$ is holomorphic in the fiber directions since
$\bO^{(1)}_\zeta (I^{(1)}_\zeta(\alpha), \beta)  = 
i\bO^{(1)}_\zeta (\alpha, \beta)$
on $\bH^k\times\A_1$.
Now
\begin{align*}
\bO^{(1)}_\zeta (\alpha, \beta) & = 
\Omega_{(J,K)} \big( (1 - \zeta K)\alpha, (1 - \zeta K)\beta\big)\\ & = 
\Omega_{(J,K)}(\alpha,\beta) \\ 
& \qquad \qquad + \, \zeta\;\big(\Omega_{(J,K)} (K\alpha,\beta) \,+\, 
\Omega_{(J,K)}(\alpha,K\beta)\big)   \\
& \qquad \qquad \qquad \qquad +\,
\zeta^2\; \Omega_{(J,K)}(K\alpha,K \beta) \\
& = \Omega_{(J,K)}(\alpha,\beta)  \; - \;
2\, \zeta\; \omega_I(\alpha,\beta) \; + 
\zeta^2\;  \Omega_{(-J,K)}(\alpha,\beta).
\end{align*}
That is,
\begin{equation}\label{eq:csymp1}
\bO^{(1)}_\zeta = \Omega_{(J,K)} - 2\zeta \omega_I + \zeta^2\Omega_{(-J,K)};
\end{equation}
therefore, $\bO^{(1)}$ is holomorphic with respect to $\zeta$.

Similarly, 
$\bO^{(2)}_\xi (I^{(2)}_\xi(\alpha), \beta)  = 
i\bO^{(2)}_\xi (\alpha, \beta)$
on $\bH^k\times\A_2$,
and 
\begin{equation}\label{eq:csymp2}
\bO^{(2)}_\xi = \Omega_{(-J,K)} - 2\xi \omega_I  + \xi^2\Omega_{(J,K)}.
\end{equation}
\end{proof}
\noindent
(Formulas \eqref{eq:csymp1},\eqref{eq:csymp2} differ slightly from
formula (3.87) of~\cite{HKLR}, due to our differing conventions;
see the discussion after \eqref{eq:explicitstereo}.)

\begin{lem} Let $\tbalpha$ be the real structure on $\tZ$ and $\bO$
the $\Os(2)$-valued holomorphic exterior 2-form defined above.
Then
\begin{equation*}
(\tbalpha)^*\bO = -\overline{\bO}.
\end{equation*}
\end{lem}
\begin{proof}
The action of $\tbalpha$ on the coefficient bundle $\Os(2)$ is defined
in Lemma~\ref{lem:tbalpha}. Now
\begin{equation*}
\bO^{(1)}_0 = \Omega_{(J,K)} = - \overline{\Omega_{(-J,K)}} =
- \overline{\bO^{(2)}_0}.
\end{equation*}
Let $\xi\zeta = 1$. Applying \eqref{eq:csymp1},
\begin{equation*}
\bO^{(1)}_{-\bar\xi} = -\bar\xi^2 \bO^{(1)}_{\zeta}.
\end{equation*}
\end{proof}
\subsection{The lattice quotient}
We have described the essential structures on the twistor space 
$\tZ \longrightarrow \P^1$ of the quaternionic vector space $V_\bH$.
That is $\tZ \cong \bH^k \times \P^1$ with the complex structure
which restricts to the fiber over $[v]\in\P^1$ by the corresponding
complex structure defined by \eqref{eq:explicitstereo}. The coordinate
charts $F^{(i)}$ defined in \eqref{eq:holocharts} 
trivialize the holomorphic vector bundle $\tZ$ over the affine 
patches $\A_1,\A_2$ of $\P^1$.

The twistor space for the \hk manifold $V_\bH/L$ is the quotient 
$Z = \tZ / L $  with the  fibration $Z \xrightarrow{\Pi_Z}\P^1$ 
induced from $\Pi_{\tZ}$. 
Denote the fiber $\Pi_Z^{-1}\big(\psi^{(1)}(\zeta)\big)$ by $Z_\zeta$
and the fiber $\Pi_Z^{-1}\big(\psi^{(2)}(0)\big)$ by $Z_\infty$, similar
to $\tZ_\zeta$ and $\tZ_\infty$. 
The complex structure on $\tZ$, the complex $\Pi_{\tZ}^*\Os(2)$-valued
exterior form $\bO$ and the real structure 
$\tbalpha$ are all invariant under translations in $L$, and therefore
induce corresponding structures on the quotient $Z = \tZ/L$.

The collection of twistor lines $l_{\vv}$ is invariant under the action
of $L$.  If $\gamma\in L$, then
\begin{equation*}
l_{\vv + \gamma} = l_{\vv} + \gamma.  
\end{equation*}
Thus the twistor lines in $\tZ$  define {\em twistor lines\/} in $Z$,
namely holomorphic sections with normal bundle isomorphic to 
$\C^{2k}\otimes\Os(1)$.

In summary we obtain the theorem of 
Hitchin-Karlhede-Lindstr\"om-Ro\u{c}ek~\cite{HKLR}:
\begin{thm}\label{thm:hktwistor}
The twistor space $Z$ is a complex manifold with 
holomorphic fibration $Z \xrightarrow{\Pi_Z}\P^1$.
Furthermore there exists a holomorphic 
$\Pi_Z^*\Os(2)$-valued exterior $2$-form $\bO$ on $Z$ inducing
a nondegenerate complex-symplectic structure on each fiber of $\Pi_Z$.
There exists a family of holomorphic sections $l_{\vv}$ of $\Pi_Z$
each with normal bundle isomorphic to $\C^{2k}\otimes\Os(1)$. 
There exists a real structure $\tbalpha$ on $Z$ covering the antipodal map 
$\P^1 \xrightarrow{\balpha}\P^1$, preserving each twistor line $l_\vv$
and $(\tbalpha)^*\bO = -\overline{\bO}$.
\end{thm}
In terms of the coordinate atlas, $L$ acts on $\tZ$ as follows.
On the fiber %$\Pi_{\tZ}^{-1}\big(\psi^{1}(0)\big)$, 
$\tZ_0$, the action of $\gamma\in L$ is:
\begin{equation}\label{eq:actionL}
F^{(1)}\bigg( \bmatrix \qq \\ \pp \endbmatrix, \psi^{(1)}(0)\bigg)
\stackrel{\gamma}\longmapsto
F^{(1)}\bigg( \bmatrix \qq + \gamma \\ \pp \endbmatrix, \psi^{(1)}(0)\bigg).
\end{equation}
Applying \eqref{eq:antipodpoles}, the action of $\gamma\in L$ on
the opposite fiber $\tZ_\infty$ is:
\begin{equation*}
F^{(2)}\bigg( \bmatrix \qq \\ \pp \endbmatrix, \psi^{(2)}(0)\bigg)
\stackrel{\gamma}\longmapsto
F^{(2)}\bigg( \bmatrix \qq  \\ \pp - i\bar\gamma \endbmatrix, 
\psi^{(2)}(0)\bigg)
\end{equation*}
Applying \eqref{eq:twistorlinepq} and \eqref{eq:twistorlinepqxi}, the
action of $\gamma\in L$ on the other fibers are:
\begin{align*}
F^{(1)}\bigg( \bmatrix \qq \\ \pp \endbmatrix, \psi^{(1)}(\zeta)\bigg)
& \stackrel{\gamma}\longmapsto
F^{(1)}\bigg( \bmatrix \qq + \gamma \\ \pp - i\zeta\bar\gamma \endbmatrix, 
\psi^{(1)}(\zeta)\bigg)  \\
F^{(2)}\bigg( \bmatrix \qq \\ \pp \endbmatrix, \psi^{(2)}(\xi)\bigg)
& \stackrel{\gamma}\longmapsto F^{(2)}\bigg( 
\bmatrix \qq + \xi\gamma  \\ \pp - i\bar\gamma \endbmatrix, 
\psi^{(2)}(\xi)\bigg)
\end{align*}

\subsection{Functions and flows}\label{sec:functionsflows}
The holomorphic $\C^*$-action on the fiber $Z_0 =
T^*\Jac(X)$ in \S\ref{sec:dolbgpd} extends to a holomorphic
$\C^*$-action on $Z$. This $\C^*$-action covers a holomorphic
$\C^*$-action on $\P^1$, where the action of $\lambda\in\C^*$ is
defined by:
\begin{align*}
\psi^{(1)}(\zeta) &\longmapsto  \psi^{(1)}(\lambda\zeta) \\
\psi^{(2)}(\xi) &\longmapsto  \psi^{(2)}(\lambda^{-1}\xi).
\end{align*}
The $\C^*$-action on $Z$ preserves the foliation by
twistor lines. Since the fiber $Z_0$ is a cross-section for
this foliation, \eqref {eq:twistorlinepq} implies:
\begin{prop}
The $\C^*$-action is described as:
\begin{equation*}
F^{(1)}\bigg( \bmatrix \qq \\ \pp \endbmatrix, 
\psi^{(1)}(\zeta) \bigg)
\stackrel{\lambda}\longmapsto
F^{(1)}\bigg( 
\bmatrix 
\frac{1 + \vert\lambda\zeta\vert^2}{1 + \vert\zeta\vert^2}\;\qq
\,+\, i\zeta \;\frac{\vert\lambda\vert^2-1}{1 + \vert\zeta\vert^2}\; \bar\pp\
 \\
\lambda\pp \endbmatrix, 
\psi^{(1)}(\lambda\zeta) \bigg)
\end{equation*}
in $\A_1$ and
\begin{equation*}
F^{(2)}\bigg( \bmatrix \qq \\ \pp \endbmatrix, 
\psi^{(2)}(\xi) \bigg)
\stackrel{\lambda}\longmapsto
F^{(2)}\bigg( 
\bmatrix 
\lambda^{-1} \;
\frac{\vert\lambda\vert^2 + \vert\xi\vert^2}{1 + \vert\xi\vert^2}
\;\qq \,+\, i\lambda^{-1}\xi \;
\frac{\vert\lambda\vert^2-1}{1 + \vert\xi\vert^2}\, \bar\pp \ \\
\pp \endbmatrix, 
\psi^{(2)}(\lambda^{-1}\xi) \bigg)
\end{equation*}
in $\A_2$.
In particular, if $\vert\lambda\vert = 1$, then
\begin{equation*}
F^{(1)}\bigg( \bmatrix \qq \\ \pp \endbmatrix, \psi^{(1)}(\zeta) \bigg)
\stackrel{\lambda}\longmapsto
F^{(1)}\bigg( \bmatrix \qq \\ \lambda\pp \endbmatrix, 
\psi^{(1)}(\lambda\zeta) \bigg).
\end{equation*}
\end{prop}
\begin{proof}
The point represented by
\begin{equation*}
F^{(1)}\bigg( \bmatrix \qq \\ \pp \endbmatrix, \psi^{(1)}(\zeta) \bigg) 
\in \Pi_{\tZ}^{-1}(\A_1)
\end{equation*}
lies on the same twistor line as
\begin{equation*}
F^{(1)}\bigg( \frac1{1+\vert\zeta\vert^2}
\bmatrix \qq - i \zeta \bar\pp \\ \pp + i\zeta\bar\qq \endbmatrix, 
\psi^{(1)}(0) \bigg) 
\end{equation*}
which maps under $\lambda\in\C^*$ to
$F^{(1)}(\vv_\lambda,\psi^{(1)}(0))$, where
\begin{equation*}
\vv_\lambda := \bmatrix \qq_\lambda \\ \pp_\lambda \endbmatrix := 
 \frac1{1+\vert\zeta\vert^2}
\bmatrix \qq - i \zeta \bar\pp \\ \lambda(\pp + i\zeta\bar\qq) \endbmatrix. 
\end{equation*}
The twistor line containing $F^{(1)}\big(\vv_\lambda,\psi^{(1)}(0)\big)$ 
intersects the the fiber $Z_{\lambda\zeta}$ at
\begin{equation*}
F^{(1)}\bigg( 
\bmatrix \qq_\lambda +  i (\lambda\zeta) \bar\pp_\lambda \\ 
\pp_\lambda - i(\lambda\zeta)\bar\qq_\lambda \endbmatrix, 
\psi^{(1)}(\lambda\zeta)\bigg).
\end{equation*}
The other claims follow similarly.
\end{proof}
As $\lambda\longrightarrow 0$, the orbit converges to a point in $\M$ which
is fixed by the $\C^*$-action and hence lies in the critical set
of the Hitchin map \eqref{eq:Hitchinmap}.

The circle action inside the $\C^*$-action is Hamiltonian for the
symplectic structure $\omega_I$. The Hamiltonian potential function
is the energy function defined in \eqref{eq:energy}:
\begin{align*}
\bH^k/L & \xrightarrow{\phi} \R \\ 
[\qq + \pp J] & \longmapsto  \frac12 \pp\cdot\bar\pp.
\end{align*}
Compare \eqref{eq:hamciract}. In contrast:
\begin{thm}\label{thm:kpotential}
The function
$\phi$ is a \ka~ potential for the metric $g$ and
the complex structure
\begin{equation*}
e^{I \theta}J = \cos(\theta) J + \sin(\theta) K 
\end{equation*}
for each $\theta \in \R$.
\end{thm}

Let $(M,g,J)$ be a \ka~ manifold, that is, $M$ is a manifold with
complex structure $J$ and Riemannian metric $g$ which is \kan~ with
respect to $J$.  Recall that a {\em \ka~ potential\/} is a function
$M\xrightarrow{f}\R$ such that 
\begin{equation*}
-\frac12 d( df\circ J) = \omega_J.
\end{equation*}
Here $df$ denotes the exterior derivative and $df\circ J$ denotes
the $1$-form (essentially the Hodge $\star$-operator applied to $df$)
obtained by composition
\begin{equation*}
TM \xrightarrow{J} TM \xrightarrow{df} \R. 
\end{equation*}
The notation $\omega_J$ was introduced in \eqref{eq:kahlerform}.
The proof of Theorem~\ref{thm:kpotential} divides into several lemmas:
\begin{lem}\label{lem:everypot}
For every $u\in S^2$, the function
\begin{align*}
\bH^k & \xrightarrow{\rho} \R \\ 
\qq + \pp J & \longmapsto  \frac12 (\qq\cdot\bar\qq + \pp\cdot\bar\pp)
\end{align*}
is a \ka~ potential for the Riemannian metric
\begin{equation*}
g =  d\qq\cdot d\bar\qq + d\pp\cdot d\bar\pp
\end{equation*}
on $\bH^k$.
\end{lem}
\begin{proof}
The rotation group $\SOth$ acts transitively on $S^2$ and isometrically
on $(\bH^k,g)$. Thus, by symmetry we may assume that $u=I$.
Now
\begin{align*}
\rho  & = \frac12 (\qq\cdot\bar\qq + \pp\cdot\bar\pp) \\
d\rho  & = \frac12 (\qq\cdot d\bar\qq + \bar\qq\cdot d\qq
+ \pp\cdot d\bar\pp + \bar\pp\cdot d\pp) \\
d\rho\circ I  & = \frac12 (-i \qq\cdot d\bar\qq + 
i \bar\qq\cdot d\qq -i \pp\cdot d\bar\pp + i \bar\pp\cdot d\pp) \\
& = \frac{i}2
(-\qq\cdot d\bar\qq + \bar\qq\cdot d\qq - 
\pp\cdot d\bar\pp + \bar\pp\cdot d\pp)
\end{align*}
and 
\begin{align*}
d (d\rho\circ I)  & = \frac{i}2
(-d \qq\cdot d\bar\qq + d\bar\qq\cdot d\qq - 
d \pp\cdot d\bar\pp + d\bar\pp\cdot d\pp) \\
& = \omega_I 
\end{align*}
as claimed.
\end{proof}
\noindent
If $u$ is a complex structure, then a function $M\xrightarrow{\psi}\R$ is
{\em $u$-pluriharmonic\/} if $d(d\psi\circ u) = 0$. 
The sum of a \ka~ potential with a pluriharmonic function is
again a \ka~ potential.
Theorem~\ref{thm:kpotential} follows from 
adding to $\rho$ (the \ka~ potential of Lemma~\ref{lem:everypot})
the pluriharmonic function 
$-\psi$ defined in the following lemma:

\begin{lem}\label{lem:circlepot}
For each $\theta\in\R$, the function
\begin{align*}
\bH^k/L & \xrightarrow{\psi} \R \\ 
\qq + \pp J & \longmapsto  \qq\cdot\bar\qq -\pp\cdot\bar\pp
\end{align*}
is $e^{I\theta}J$-pluriharmonic.
\end{lem}
\begin{proof}
From the definition of $J$,
\begin{align*}
d\pp \circ J & = d\bar\qq \\
d\qq \circ J & = -d\bar\pp \\
d\bar\pp \circ J & = -d\qq \\
d\bar\qq \circ J & = d\pp.  
\end{align*}
Now
\begin{align*}
\psi  & = \qq\cdot\bar\qq  - \pp\cdot\bar\pp \\
d\psi  & = \qq\cdot d\bar\qq + \bar\qq\cdot d\qq
- \pp\cdot d\bar\pp - \bar\pp\cdot d\pp \\
d\rho\circ e^{I\theta}  & = e^{-i\theta} \qq\cdot d\bar\qq + 
e^{i\theta} \bar\qq\cdot d\qq -e^{-i\theta} \pp\cdot d\bar\pp 
- % + i  % changed with Eugene 090705
e^{i\theta} \bar\pp\cdot d\pp \\  % changed - in exponent to + 
d\rho\circ e^{I\theta}\circ J & = e^{-i\theta} \qq\cdot d\pp + 
e^{i\theta} \bar\qq\cdot (-d\bar\pp) -e^{-i\theta} \pp\cdot (-d\qq) 
- e^{i\theta} \bar\pp\cdot d\bar\qq \\ %
& = e^{-i\theta} \qq\cdot d\pp - 
e^{i\theta} \bar\qq\cdot d\bar\pp + e^{-i\theta} \pp\cdot d\qq 
- e^{i\theta} \bar\pp\cdot d\bar\qq \\ %
& = d( e^{-i\theta} \pp\cdot\qq - e^{i\theta}\bar\pp\cdot\bar\qq) %
\end{align*}
is closed, that is, $d \big(d\rho\circ (e^{I\theta}J)\big) = 0$, as claimed.
\end{proof}
\noindent
The function $\phi$ is {\em not\/} a 
\ka~ potential for any of the other complex structures in the \hk family.

\section{The moduli space and the Riemann period matrix}
Here we explicitly compute these moduli spaces in terms of 
classical invariants.
We return to the original point of view, involving the fundamental
group and the Betti moduli space. 

\subsection{Coordinates for the Betti moduli space}
Under the direct product decomposition of Betti moduli space 
\begin{align*}\label{eq:directproduct}
\Hom(\pi,\C^*)  &\xrightarrow{\cong} \Hom(\pi,\Uo)\times\Hom(\pi,\R^+) \\
\rho & \longmapsto (\rho_u,\rho_\R),
\end{align*}
the representation $\rho$ decomposes into a unitary part $\rho_u$ 
and a positive real part $\rho_\R$ as:
\begin{equation*}
\rho_u(\gamma) = \frac{\rho(\gamma)}{\vert\rho(\gamma)\vert},\qquad
\rho_\R(\gamma) = \vert\rho(\gamma)\vert.
\end{equation*}
By the real analytic isomorphism $\R^+\xrightarrow{\log}\R$,
we identify
\begin{equation*}
\Hom(\pi,\R^+) \longrightarrow  \Hom(\pi,\R) \longrightarrow 
\H^1(X,\R) \cong \R^{2k}
\end{equation*}

Let $A_1,B_1,\dots,A_k,B_k\in\pi$ be generators
for the presentation \eqref{eq:presentation}. 
In terms of these generators,
a representation $\rho\in\Hom(\pi,\C^*)$ is determined by
an arbitrary $2k$-tuple
\begin{equation*}
\big(\rho(A_1),\dots, \rho(B_k)\big) \in (\C^*)^{2k}.
\end{equation*}
The unitary part $\rho_u$ is then given by
\begin{equation*}
\big(\frac{\rho(A_1)}{\vert\rho(A_1)},\dots, 
\frac{\rho(B_k)}{\vert\rho(B_k)\vert}\big) \in \Uo^{2k}
\end{equation*}
and the positive-real part by
\begin{equation*}
\big(\vert\rho(A_1)\vert,\dots, \vert\rho(B_k)\vert\big)\in(\R^+)^{2k}
\end{equation*}
which we identify with the real vector
\begin{equation*}
\big(\log \vert\rho(A_1)\vert,\dots, \log\vert\rho(B_k)\vert\big)
\in\R^{2k}.
\end{equation*}

With a conformal structure on $X$, $\rho_u$ corresponds to a point
in $\Jac(X)=V/L$, where $V\cong\C^k$ is a complex vector space and 
$L$ is the lattice in $\C^k$ spanned by $\Z^k$ and the columns of the
period matrix $\Pi$. The real part $\rho_\R$ 
%(the Higgs field) 
corresponds to a point 
in the real symplectic vector space $\H^1(X;\R)\cong\R^{2k}$.  Thus,
the conformal structure provides a linear $\C^*$-action generated by a complex
structure $\J_\Pi$ depending on $\Pi$.

\subsection{Abelian differentials and their periods}
The homology classes
corresponding to the generators $A_1, B_1 \dots,A_k, B_k\in\pi$
determine a symplectic basis (also denoted by $A_1, B_1 ,\dots A_k, B_k$) 
of $\H_1(X;\Z)$ satisfying:
\begin{align*}
A_j \cdot A_l =  A_j \cdot A_l & =  0 \\   
A_j \cdot B_l =  -B_j \cdot A_l & = \delta_{j,l}   
\end{align*}
with respect to the intersection form on $\H_1(X;\Z)$.
If $\phi,\psi$ are $1$-forms on $X$,
then
\begin{equation}\label{eq:cupproduct}
\int_X \phi \wedge \psi \;=\; \sum_{j=1}^k \bigg(
\int_{A_j} \phi  \int_{B_j} \psi -
\int_{B_j} \phi  \int_{A_j} \psi\bigg). 
\end{equation}

Recall that the space $\bV=\hH^{1,0}(X)$ of abelian differentials is a
complex vector space for which the Hermitian inner product defined by
\eqref{eq:hermitianstructure} is negative definite.  Choose a basis
$\{\omega_1,\dots,\omega_k\}$ of abelian differentials. %%%$\bV=\hH^{1,0}(X)$ 
Call such a basis {\em unitary\/} if the corresponding
basis $\{\bar\omega_1,\dots,\bar\omega_k\}$ of $V$
is orthonormal with respect to the {\em positive definite\/} Hermitian
inner product on $V$ induced by \eqref{eq:hermitianstructure}.

We record this basis as a column vector
of $1$-forms
\begin{equation*}
\omega = \bmatrix \omega_1 \\ \vdots \\ \omega_k \endbmatrix.
\end{equation*}
If $\{\omega_1',\dots,\omega_k'\}$ is another basis, then the
corresponding vectors $\omega$ and $\omega'$ satisfy
\begin{equation*}
\omega' = U \omega,
\end{equation*}
where $U\in\GL(k,\C)$ is an invertible $k\times k$ matrix. 
If $\omega$ is a unitary basis, then  $\omega'$ is unitary if and only
if $U\in \U(k)$ is a unitary matrix.

If $\alpha,\beta\in\bV$ are abelian differentials, 
then $\star\alpha =i\alpha$ and $\star\bar\beta = -i\beta$, 
and \eqref{eq:hermitianstructure} implies 
\begin{equation*}
\langle \alpha,\beta \rangle =  \int_X \alpha \wedge \bar{\beta}  
\end{equation*}

The {\em period matrix\/} for $\omega$ is defined as the $k\times 2k$
matrix $\bmatrix \bA & \bB \endbmatrix$, where the $k\times k$-matrices 
$\bA,\bB$ are defined by:
\begin{equation}\label{eq:defnAB}
\bA_{j,l} := \int_{A_l}\omega_j, \qquad 
\bB_{j,l} := \int_{B_l}\omega_j  
\end{equation}
respectively.
Formula \eqref{eq:cupproduct} implies that $\bA$ and $\bB$ are invertible.
If $\omega'$ is another basis related to $\omega$ by 
$\omega' = U \omega$, then the period matrix
$\bmatrix \bA' & \bB' \endbmatrix$ of $\omega'$  relates to the 
period matrix of $\omega$ by:
\begin{equation}\label{eq:periodmatrixtransformation}
\bmatrix \bA' & \bB' \endbmatrix = U \bmatrix \bA & \bB \endbmatrix.
\end{equation}
Since the $(2,0)$-form $\omega_j\wedge\omega_l$ is  identically zero, 
\begin{equation*}
0  =\; \int_X \omega_j\wedge\omega_l 
\end{equation*}
and  \eqref{eq:cupproduct} implies:
\begin{equation}\label{eq:generalperiods1}
0 = \bA \bB^\dag - \bB \bA^\dag.  
\end{equation}
Suppose that $\omega$ is a unitary basis. Then
%$\dim(X)=1$ and unitarity implies:
\begin{equation*}
- i\,\delta_{j,l}  =\; \int_X \omega_j\wedge\bar{\omega}_l 
\end{equation*}
and \eqref{eq:cupproduct} implies
\begin{equation}\label{eq:generalperiods2}
-i\Id  = \bA \bar{\bB}^\dag - \bB \bar{\bA}^\dag. 
\end{equation}
The {\em normalized basis of abelian differentials\/}
$\{\omega_1',\dots,\omega_k'\}$ is the basis such that
\begin{equation*}
\int_{A_j} \omega_l' = \delta_{j,l}.
\end{equation*}
Together with \eqref{eq:periodmatrixtransformation}, this implies:
\begin{equation*}
\omega' = \bA^{-1} \omega,
\end{equation*}
where $\bA$ is defined in \eqref{eq:defnAB}.
By 
\eqref{eq:periodmatrixtransformation}, 
the corresponding {\em period matrix\/}
%%%
is 
\begin{equation*}
\bmatrix \bA' & \bB' \endbmatrix = 
\bmatrix \Id \quad \Pi \endbmatrix 
\end{equation*}
where $\Pi = \bA^{-1} \bB$ is the matrix 
%%%
\begin{equation*}
\Pi_{i,j} :=  \int_{B_j} \omega_i'.
\end{equation*}
\eqref{eq:generalperiods1} and \eqref{eq:generalperiods2} 
imply the {\em Riemann bilinear relations:\/}
\begin{align*}
\Pi & = \Pi^{\dag}  \\
\Im(\Pi) & = \frac{\Pi-\bar\Pi}{2i} = \frac12 (\bA^{-1})
\overline{(\bA^{-1})}^\dag > 0.
\end{align*}
That is,
$\Pi$ is a $k\times k$ symmetric matrix whose imaginary part is positive
definite.

For a symmetric $k\times k$ complex matrix $\Pi$ 
with positive definite imaginary part, 
the matrix $\bmatrix \bA & \bB \endbmatrix$ corresponding
to a unitary basis can be recovered. Namely, let $\sqrt{2\ \Im(\Pi)}$ denote
the unique positive definite symmetric square-root of $2\Im(\Pi)$ and let
$U\in\U(k)$ be an arbitrary unitary matrix. Then 
\begin{align*}
\bA & \;:=\;  U \big(\sqrt{2\,\Im(\Pi)}\big)^{-1} \\
\bB & \;:=\;  U \big(\sqrt{2\,\Im(\Pi)}\big)^{-1} \, \Pi 
\end{align*}
define a period matrix $\bmatrix \bA & \bB \endbmatrix$
which is unique up to the choice of $U$.

\subsection{Flat connections}
Consider a flat connection $D$ on a trivial line bundle over $X$. Write
\begin{equation*}
D = D_0 + \psi + \phi,  
\end{equation*}
where $D_0$ is the flat connection arising from the trivialization,
$\psi$ purely imaginary and $\phi$ real 1-forms. By a gauge
transformation, one may write
\begin{align*}
\psi & =  \Psi - \bar{\Psi} \\
\phi & =  \Phi + \bar{\Phi},
\end{align*}
where $\Psi\in V$ and $\Phi\in\bV$.
The corresponding vectors $\qq,\pp\in \C^k$ satisfy
\begin{equation*}
\Psi = \sum_{j=1}^k q_j \bar{\omega}_j = \bar\omega ^\dag\cdot\qq,\qquad
\Phi = \sum_{j=1}^k p_j \omega_j = \omega^\dag\cdot\pp,
\end{equation*}
where
\begin{equation*}
\qq = \bmatrix q_1 \\ \vdots \\ q_k \endbmatrix, \qquad 
\pp = \bmatrix p_1 \\ \vdots \\p_k \endbmatrix.  
\end{equation*}
The Higgs pair consists of the holomorphic structure $D'' = D_0'' +
\qq^\dag\bar\omega$ and the Higgs field  $\Phi=\pp^\dag\omega$.
Corresponding to $\qq,\pp\in\C^k$  is the connection
$D = D_0 + \eta$, where
\begin{equation*}
\eta  =  
\big(\qq^\dag \bar{\omega} \, -\, \bar{\qq}^\dag \omega\big)\, +\, 
\big(\bar{\omega}^\dag\cdot \qq  \, -\, 
\omega^\dag \cdot\bar{\qq}\big)\, +\, 
\big(\omega^\dag\cdot\pp  \,+\, 
\bar{\omega}^\dag\cdot\bar{\pp} \big).
\end{equation*}
\noindent
To simplify calculations, henceforth 
let $\omega$ be the normalized basis of abelian differentials such
that $\bA = \Id$ and $\bB=\Pi$.

If $D$ is a flat unitary connection, then $D = D_0 + \psi$, where
$\psi$ is purely imaginary.  By a gauge transformation, one
may assume that $\psi = \Psi - \bar{\Psi}$, where $\Psi\in\hH^{0,1}(X)$.
The corresponding vector $\qq\in\C^k$ is defined by
\begin{equation*}
\Psi = \sum_{j=1}^k q_j \bar\omega_j = \bar\omega^\dag \cdot \qq 
\end{equation*}
and
\begin{equation*}
\psi =  \bar\omega^\dag \cdot \qq  - \omega^\dag \cdot \bar\qq   =
2i \Im \big( \bar\omega^\dag \cdot \qq\big),
\end{equation*}
with corresponding periods
\begin{align*}
\int_{A_l} \psi & = 2 i \Im (q_l) \\ 
\int_{B_l} \psi & = 2 i \Im \big(\sum_{j=1}^k q_j \bar\Pi_{j,l}\big).
\end{align*}
In the notation
\begin{equation*}
\ba := 
\bmatrix  \int_{A_1} \psi \\ \vdots \\  \int_{A_k} \psi \endbmatrix, \qquad
\bbb := \bmatrix  \int_{B_1} \psi \\ \vdots \\  \int_{B_k} \psi \endbmatrix,
\end{equation*}
the periods are:
\begin{align*}
\ba & = 2 i \Im(\qq) \\ 
\bbb & = 2 i \Im(\bar\Pi\cdot\qq) = 
2 i \big(\Re(\Pi) \cdot\Im(\qq) -\Im(\Pi)\cdot\Re(\qq)\big).
\end{align*}
The lattice $L$ in $\C^k$ consists of $\qq$ such that the periods
$\ba,\bbb\in 2\pi i\Z$, that is,
\begin{equation*}
\Im(\qq)\in \pi\Z^k,\qquad  
\Re(\Pi)  \Im(\qq) -\Im(\Pi) \Re(\qq)\in \pi\Z^k
\end{equation*}
which is equivalent to:
\begin{equation*}
\qq \in \pi \Im(\Pi)^{-1} \big\{ \Z^k + \Pi\Z^k\big\}.
\end{equation*}
Thus $L$ is equivalent (by $\pi \Im(\Pi)^{-1}\in\GL(k,\C)$) to
the lattice $\Z^k + \Pi\Z^k$
spanned by the $2k$ columns of $\bmatrix \Id_k & \Pi \endbmatrix$.
Therefore $\Jac(X)$ is the quotient of $\C^k$ by 
$\Z^k + \Pi\Z^k$.
(See, for example, Gunning~\cite{Gunning1,Gunning3}, \S 2.5.)

\subsection{Higgs fields}
The space $\Hom(\pi,\R^+)$ of
positive real characters inherits structure from the conformal structure 
of $X$ through the period matrix $\Pi$.
For $\rho\in\Hom(\pi,\C^*)$, write
\begin{equation*}
\alpha_j := \frac12 \log \vert \rho(A_j)\vert, \qquad
\beta_j := \frac12 \log \vert \rho(B_j)\vert.
\end{equation*}
We determine the corresponding flat $\R^+$-connection with holonomy
$(\alpha,\beta)$.

After applying a gauge transformation, the $1$-form 
$\phi\in\hH^1(X;\R)$ decomposes as  $\phi = \Phi + \bar\Phi$,  
where $\Phi\in\hH^{1,0}(X)$ is a holomorphic $1$-form.
The vector $\pp\in\C^k$ corresponding to $\Phi$ is defined by
\begin{equation*}
\Phi = \sum_{j=1}^k p_j \omega_j = \omega^\dag   \pp
\end{equation*}
which has periods
\begin{align*}
\alpha_j & = \frac12 \int_{A_j} \phi = \Re \int_{A_j}\Phi = \Re(p_j) \\
\beta_j & = \frac12 \int_{B_j} \phi = \Re \int_{B_j}\Phi = 
\Re\big(\sum_{l=1}^k\Pi_{j,l} p_l)\big).
\end{align*}
Thus
\begin{align*}
\alpha & = \Re(\pp)  \\
\beta & = \Re(\Pi)  \Re(\pp) - \Im(\Pi)  \Im(\pp)
\end{align*}
with inverse mapping
\begin{align*}
\Re(\pp) &= \alpha  \\
\Im(\pp) & = \Im(\Pi)^{-1}\big( \Re(\Pi)  \alpha - \beta\big).
\end{align*}
As matrices,
\begin{equation*}
\bmatrix \Re(\pp) \\ \Im(\pp) \endbmatrix \longmapsto 
\bmatrix \Id_k & 0 \\ \Re(\Pi) & -\Im(\Pi) \endbmatrix 
\bmatrix \alpha \\ \beta \endbmatrix 
\end{equation*}
and 
\begin{equation*}
\bmatrix \alpha \\ \beta \endbmatrix \longmapsto 
\bmatrix \Id_k & 0 \\ \Im(\Pi)^{-1}\Re(\Pi) & -\Im(\Pi)^{-1} \endbmatrix 
\bmatrix \Re(\pp) \\ \Im(\pp) \endbmatrix. 
\end{equation*}

\subsection{The $\C^*$-action in terms of the period matrix}
We describe the $\C^*$-action explicitly as follows.
We break $\pp$ into its real and imaginary parts:
\begin{equation*}
\pp = \Re(\pp) + i \Im(\pp) 
\end{equation*}
Scalar multiplication 
$\pp \longmapsto i\pp$ is:
\begin{equation*}
\bmatrix \Re(\pp) \\ \Im(\pp) \endbmatrix \longmapsto 
\bmatrix \Im(\pp) \\ -\Re(\pp) \endbmatrix = 
\bmatrix 0 & \Id_k \\ -\Id_k & 0 \endbmatrix
\bmatrix \Re(\pp) \\ \Im(\pp) \endbmatrix.
\end{equation*}
Multiplying the Higgs field $\Phi$ by $i$ transforms the periods $\alpha,\beta$
by:
\begin{equation*}
\bmatrix \alpha \\ \beta \endbmatrix \longmapsto \J_\Pi
\bmatrix \alpha \\ \beta \endbmatrix 
\end{equation*}
where
\begin{align}\label{eq:JPi}
\J_\Pi & =  
\bmatrix \Id_k & 0 \\ \Re(\Pi) & -\Im(\Pi) \endbmatrix 
\bmatrix  0 & \Id_k \\ -\Id_k & 0 \endbmatrix 
\bmatrix \Id_k & 0 \\ \Im(\Pi)^{-1}\Re(\Pi) & -\Im(\Pi)^{-1} \endbmatrix 
\notag\\ \ & \ \notag \\ & =
\bmatrix -\Im(\Pi)^{-1}\Re(\Pi) & \quad\Im(\Pi)^{-1} \\
-\Re(\Pi)\Im(\Pi)^{-1}\Re(\Pi) -\Im(\Pi) 
& \quad\Re(\Pi)\Im(\Pi)^{-1} \endbmatrix.
\end{align}
Then the $\C^*$-action is given by
\begin{equation*}
\Re(\lambda) \Id_{2k} + \Im(\lambda) \J_\Pi 
\end{equation*}
for $\lambda\in\C^*$. 

\subsection{The $\C^*$-action and the real points}
The $\C^*$-action on $\Hom(\pi,\C^*)$ fixes the unitary component
and acts as above on the Higgs field. That is, $\lambda\in\C^*$ maps
\begin{equation*}
\rho = \big(\rho_u,\rho_\R\big) \longmapsto 
\bigg(\rho_u,\exp\big( 
(\Re(\lambda) \Id + \Im(\lambda) J_\Pi) \log (\rho_\R)\big)\bigg),
\end{equation*}
where $\J_\Pi$ is defined by \eqref{eq:JPi}.

The unitary projection $\Hom(\pi,\C^*)\longrightarrow \Hom(\pi,\Uo)$
associates to a Higgs line bundle the underlying holomorphic line bundle.
It can be described purely in terms of the $\C^*$-action as follows.
Let $\rho\in\Hom(\pi,\C^*)$. Then as $\lambda\longrightarrow 0$,
the orbit $\lambda\cdot \rho$ approaches $\rho_u$. In particular this
structure defines the foliation of $\Hom(\pi,\C^*)$ by copies of
$\Hom(\pi,\R^+) = \H^1(X;\R)$.

Moreover the complex-symplectic structure $(J,\Omega)$ defined by 
\eqref{eq:OmegaJ} arises from the complex-orthogonal structure
on $\C^*$. By \eqref{eq:OmegaRI} (compare also \eqref{eq:OmegaIK}),
the imaginary part of $\Omega$ is a real-symplectic structure $-\omega_K$
under which the decomposition
\begin{equation*}
\Hom(\pi,\C^*) = \Hom(\pi,\Uo)\times \Hom(\pi,\R^+)
\end{equation*}
is a Lagrangian product (a {\em real bipolarization\/}). 
The $\C^*$-action is the extra piece of structure to determine the
Dolbeault moduli space (expressed in terms of the complex structure $I$).
As $\Hom(\pi,\Uo)$
is $\omega_K$-Lagrangian, $\omega_K$ identifies 
the tangent space with the normal space,
which in turn identifies with $\H^1(X;\R)$. 
This real-symplectic space inherits a Hermitian structure from the 
complex structure $\J_\Pi$, which defines
a complex structure on $\Hom(\pi,\Uo)$.
With this complex structure,
$\Hom(\pi,\Uo)$ identifies with $\Jac(X)$.
Thus the full structure of the Dolbeault moduli space of Higgs bundles
arises from the dynamics of the $\C^*$-action and the natural structures
of the Betti moduli space which depend only on the 
fundamental group of $\Sigma$.

%\begin{thebibliography}{99}
\bibliographystyle{amsalpha}

\end{document}